\documentclass[11pt]{article}
\usepackage{amsmath}
\usepackage{amssymb}

\newtheorem{theorem}{Theorem}[section]
\newtheorem{lemma}{Lemma}[section]
\newtheorem{corollary}{Corollary}[section]

\begin{document}

\begin{center}
  {\Large\bf Pointwise decay in space and in time for incompressible viscous
     flow around a rigid body moving with constant velocity.
}
\vspace{1ex}

{\large\sl Paul Deuring}

Universit\'e du Littoral C\^ote d'Opale, Laboratoire de Math\'ematiques
Pures et Appliqu\'ees Joseph Liouville, F-62228 Calais, France.
\end{center}

\vspace{2ex}
\begin{abstract}
We present pointwise space-time decay estimates for the velocity part of
solutions to the time-dependent Oseen system in 3D, with
Dirichlet boundary conditions and vanishing velocity at infinity.
In addition, similar estimates are derived for solutions to the time-dependent incompressible
Navier-Stokes system with Oseen term, and for solutions to the stability problem
associated with the stationary incompressible Navier-Stokes system with Oseen term.



\end{abstract}

\section{Introduction}
Consider the incompressible time-dependent Navier-Stokes system
\begin{eqnarray} \label{50}&&
\partial _tv-\Delta_x v + \tau \, \partial _{1}v +\tau \,
(v \cdot  \nabla_x  ) v + \nabla _x\pi = f,
\;\;
\mbox{div}_x v = 0
\\&&\nonumber
\mbox{for}\;
t \in (0, \infty ),\;
x \in \overline{ \Omega }^c:= \mathbb{R}^3 \backslash \overline{ \Omega },
\end{eqnarray}
with the boundary conditions
\begin{eqnarray} \label{60a}
  v(t)| \partial \Omega =b(t), \quad v(x,t)\to (1,0,0) \; (|x|\to \infty )
  \quad \mbox{for}\;\; t \in (0, \infty ) ,
\end{eqnarray}
and the initial condition
\begin{eqnarray} \label{70a}
  v(0)=v_0,
\end{eqnarray}
where $\Omega $ is an open bounded set in $\mathbb{R}^3 $ with
connected Lipschitz boundary. Problem (\ref{50}) -- (\ref{70a})
is a mathematical model for the flow of a viscous incompressible fluid
around a rigid body that moves steadily and without rotation,
under the assumption that the underlying reference frame adheres to the
body, represented by the set $\Omega $. The ``exterior domain'' $\overline{ \Omega }^c$
is supposed to be filled with the fluid. The function $v$ stands for the unknown
velocity field and the function $\pi$ for the unknown pressure field of the fluid.
The real number $\tau >0$ (Reynolds number)
and the functions $f$ (volume force), $v_0$ (initial velocity) and $ b$
(velocity of the fluid particles on the surface of the body) are given.

Problem (\ref{50}) -- (\ref{70a}) is already normalized in the sense that the flow
is characterized by a single parameter -- the Reynolds number --  and the rigid
body moves with the constant velocity $(-1,0,0) $ with respect to an observer at rest.
This latter feature of the motion of the body, expressed  by the boundary condition at
infinity stated in (\ref{60a}), means in particular that the negative part of the
$x_1$-axis corresponds to the upstream and the positive part to the downstream direction
of the flow.

We will study the asymptotic behaviour of the fluid far from the rigid body.
These asymptotics are of interest because they may be interpreted as features of the
flow that may actually be observed. In this respect, particular attention
is directed  at the wake extending behind the rigid body, in our situation
around the positive $x_1$-axis. This wake should emerge in the asymptotics
provided by theory.

Since nonzero boundary conditions at infinity are inconvenient from a mathematical
point of view, we will not deal with equations (\ref{50}) -- (\ref{70a}) directly,
but with an equivalent problem. To this end, we introduce the new velocity
$u:=v-(1,0,0),$ which satisfies the Navier-Stokes system with Oseen term,
\begin{eqnarray} \label{100}
\partial _tu-\Delta_x u + \tau \, \partial _{1}u +\tau \,
(u \cdot  \nabla_x  ) u + \nabla _x\pi = f,
\;\;
\mbox{div}_x u = 0
\quad \mbox{in}\;
\overline{ \Omega }^c \times (0, \infty ),
\end{eqnarray}
as well as the side conditions
\begin{eqnarray} \label{110}&&
u(t)| \partial \Omega =(-1,0,0)+b(t), \quad u(x,t)\to 0 \; (|x|\to \infty )
\quad \mbox{for}\;\; t \in (0, \infty ) ,
\\&& \label{30}
u(0)=a,
\end{eqnarray}
with $a:=v_0-(1,0,0).$
In the work at hand, we are interested in temporal and spatial asymptotics of the velocity
part $u$ of solutions to (\ref{100}), (\ref{110}), (\ref{30}). But due to the stationary component
$(-1,0,0)$ of the
Dirichlet boundary data in (\ref{110}), the velocity $u$ cannot be expected to decay for
$t\to \infty ,$ even if $b(t)$ tends to zero. Therefore we modify our problem a second time.
To this end, we take a solution $(U,\Pi )$ of the stationary
Navier-Stokes system with Oseen term,
\begin{eqnarray} &&\label{170}
-\Delta U+\tau \, \partial _1U+\tau \, (U \cdot \nabla )U +\nabla \Pi =F,
\;\;
\mbox{div}\,U=0\;\; \mbox{in}\; \overline{ \Omega }^c,
\end{eqnarray}
under Dirichlet boundary conditions
\begin{eqnarray} \label{170b}
U| \partial \Omega =(-1,0,0)+B,\;\; U(x)\to 0\;(|x|\to \infty ),
\end{eqnarray}
with given functions $F: \overline{ \Omega }^c \mapsto \mathbb{R}^3 ,\; B: \partial \Omega \mapsto \mathbb{R}^3 $,
and then introduce the new unknowns
$\overline{ u}(x,t):= u(x,t)-U(x),\; \overline{ \pi}(x,t)-\Pi (x),$
and the new given functions
$\overline{ f}(x,t):=f(x,t)-F(x),\; \overline{ b}(x,t):= b(x,t)-B(x)$ and $ \overline{ a}:=a-U$.
For simplicity denoting these new quantities again by $u,\, \pi,\, f,\, b,\, a$, instead of
$\overline{ u},\, \overline{ \pi},\, \overline{ f},\, \overline{ b},\, \overline{ a}$, respectively,
we then arrive at the system
\begin{eqnarray}&& \label{180}
\partial _tu-\Delta_x u + \tau \, \partial _{1}u
+ \tau \, (u \cdot  \nabla_x  ) u
+ \tau \, (U \cdot  \nabla_x  ) u
+ \tau \, (u \cdot  \nabla ) U
+ \nabla _x \pi  = f,
\\&&\hspace{3em}\nonumber 
\mbox{div}_x u = 0
\quad \mbox{in}\;
\overline{ \Omega }^c \times (0, \infty ),
\end{eqnarray}
with the boundary conditions
\begin{eqnarray} \label{20b}
u(t)| \partial \Omega =b(t), \quad u(x,t)\to 0 \; (|x|\to \infty )
\quad \mbox{for}\;\; t \in (0, \infty ) ,
\end{eqnarray}
and initial condition (\ref{30}) (stability problem associated with (\ref{170}), (\ref{170b})).
According to \cite[Theorem 1.2]{DeSIMA2}, under suitable conditions on the data, the velocity part $u$
of a solution to (\ref{180}), (\ref{20b}), (\ref{30}) exhibits an asymptotic behaviour in space
described by the estimate
\begin{eqnarray} \label{130}
| \partial ^{\alpha }_xu(x,t)| \le \mbox{$\mathfrak C$} \, \bigl(\, |x|\, \nu (x) \,\bigr) ^{-1-| \alpha |/2}
\quad \mbox{for}\;\; x \in \mathbb{R}^3 \; \mbox{with}\; |x|\ge R,\; t \in (0, \infty ) ,
\end{eqnarray}
$\alpha \in \mathbb{N} _0^3$ with $ | \alpha |:=\alpha _1+\alpha _2+\alpha _3\le 1$,
where $R $ is some positive real so large that
$\overline{ \Omega }\subset B_{R}$.
The same estimate is shown in previous papers \cite{DeSIMA1}, \cite{DeDCDS-A}
for solutions to the linear system (''Oseen system'')
\begin{eqnarray} \label{10}
\partial _tu-\Delta_x u + \tau \, \partial _{1}u
+ \nabla _x \pi  = f,
\;\;
\mbox{div}_x u = 0
\quad \mbox{in}\;
\overline{ \Omega }^c \times (0, \infty ),
\end{eqnarray}
again under side conditions (\ref{20b}) and (\ref{30}).
The fact that inequality (\ref{130}) holds for
$\alpha \in \mathbb{N} _0^3$ with $ | \alpha |\le 1$
means that $u$ and the spatial gradient $\nabla _xu$ of $u$ are evaluated.
The factor $\nu(x)$ in (\ref{130}) is defined by
\begin{eqnarray} \label{101}
\nu (x):=1+|x|-x_1 \quad \mbox{for}\;\; x \in \mathbb{R}^3 .
\end{eqnarray}
Due to this factor, the right-hand side of (\ref{130})
decays less fast in the wake region around the positive $x_1$-axis than it does elsewhere;
see \cite[section VII.3]{Ga1} for more details.
For this reason the presence of the factor $\nu (x)$ in (\ref{130}) may be interpreted as a
mathematical manifestation of the wake.

Inequality (\ref{130}) deals only with spatial decay of the velocity.
However, for suitable data,
the velocity far from the body should decay
in space as well as in time.
It is the aim of the work at hand to make this idea more precise by determining
upper bounds of $|u(x,t)|$ and $|\nabla _xu(x,t)|$ reflecting this type of asymptotics.
In certain special situations, such estimates may be deduced immediately from (\ref{130})
and estimates of $\|u(t)\|_{\infty }$ 
available in literature, where $\|\;\|_{\infty }$ denotes the norm
of $L ^{ \infty } ( \overline{ \Omega }^c)^3$.
With respect to the linear problem
(\ref{10}), (\ref{20b}), (\ref{30}),
Enomoto, Shibata \cite{EnSh2} showed that if $f=0,\; b=0,\; \Omega $ smoothly bounded
and $a \in L^p( \overline{ \Omega }^c)^3 $ for some $p \in [1, \infty ),$
then
$
\|u(t)\|_{\infty }
\le
\mbox{$\mathfrak C$} \,\| a\|_p\, t^{-3/(2p)-| \alpha |/2}
$
for
$
t \in (0, \infty ) 
$
(\cite[Theorem 1.1]{EnSh2}).
This result combined with (\ref{130}) and the equation $d=d^{1-\epsilon }\,d^{\epsilon }\; (d\ge 0,\;
\epsilon \in [0,1])$
yields the pointwise space-time estimate
\begin{eqnarray*} 
| \partial ^{\alpha }_xu(x,t)| \le \mbox{$\mathfrak C$} \,
\bigl(\, |x|\, \nu (x) \,\bigr) ^{(-1-| \alpha |/2)\,(1-\epsilon )}
\, (1+t)^{(-3/(2p) -| \alpha |/2)\, \epsilon }
\end{eqnarray*}
for
$
x,\, t,\, \alpha
$
as in (\ref{130}) and for
$
\epsilon \in [0,1]
.
$
In particular, the highest rate of temporal decay that may be attained is $(1+t)^{-3/2-| \alpha |/2}$,
arising if $p=1$.

The theory of the linear problem (\ref{10}), (\ref{20b}), (\ref{30}) we present here
does not require $f$ or $b$ to vanish, nor $\Omega $ to be more regular than
Lipschitz bounded. The estimates we derive reflect the asymptotics of the data $f,\, a$ and $b$
in a rather precise way (Corollary \ref{corollaryC8.10}). In the best possible case, occuring
when $f$ and $a$ are both bounded and with compact support and $b$ decays sufficiently rapidly,
we show that for $ \zeta \in (0,1)$ arbitrary but fixed,
the inequality
\begin{eqnarray} \label{129}&&
| \partial ^{\alpha }_xu(x,t)|
\\&&\nonumber 
\le
\mbox{$\mathfrak C$} \, \bigl[\,
\bigl(\, |x|\, \nu (x) \,\bigr) ^{-1-| \alpha |/2}\, (1+t)^{-\zeta }
+
\bigl(\, |x|\, \nu (x) \,\bigr) ^{(-1-| \alpha |/2)\,(1-\epsilon )}
  \, (1+t)^{(-1 -| \alpha |/2)\, \epsilon }\,\bigr]
\end{eqnarray}
holds for $x,\,t,\, \alpha $ as in (\ref{130}) and for $\epsilon \in [0,1]$
(Corollary \ref{corollaryC8.20}).
Further below, when we sketch our method of proof, we will indicate why we do not
achieve the rate $(1+t)^{-3/2 -| \alpha |/2}$ of temporal decay.

Concerning the nonlinear problem (\ref{180}), (\ref{20b}), (\ref{30}),
decay estimates of $\|u(t)\|_{\infty }$ were provided by Masuda \cite{Ma},
Heywood \cite{Hey2}, Shibata \cite{Sh} and Enomoto, Shibata \cite{EnSh2}
for a non-normalized version of this problem, with \cite{Ma}, \cite{EnSh2} and \cite{Sh}
requiring that the viscosity equals $1$.
Masuda chose initial data in $L^2$ and assumed
$b=0,\; \partial \Omega $ smooth and $U$ small in a suitable sense.
Constructing $L^2$-weak solutions that become strong after a certain time, he obtained
$\|u(t)\|_{\infty }\le \mbox{$\mathfrak C$} \, t^{-1/8}$ for large $t$,
or instead $\|u(t)\|_{\infty }\to 0\;(t\to \infty )$, depending on the asymptotics of $f$
(\cite[p. 297, Theorem; p.298, Remark 1.3]{Ma}). Heywood \cite[p. 674-675]{Hey2},
\cite[Theorem 4]{Hey1} admitted nonvanishing $f$ and $b$, improving the decay rate of
$\|u(t)\|_{\infty }$ to $t^{-1/4}$ under various smallness conditions on $a,\, f$ and $b$.
Enomoto, Shibata \cite{EnSh2} worked with initial data in $L^3$,
constructing mild solutions under the assumptions that $f$ and $b$ vanish, $\Omega $ is smoothly
bounded and the initial data $a$ and the data of the stationary problem (\ref{170}), (\ref{170b}),
and hence $U$, are small. Within this framework, they showed that
$\|u(t)\|_{\infty }\le \mbox{$\mathfrak C$} \, t ^{-1/2} $
for $t \in (0, \infty ) $  (\cite[Theorem 1.3]{EnSh2}). In a previous article by Shibata \cite{Sh},
a similar but slightly weaker result is derived (\cite[Theorem 1.4]{Sh}).

As in the linear case, these decay estimates of $\|u(t)\|_{\infty }$ combined with (\ref{130})
yield pointwise space-time decay estimates of $u$, although not of $\nabla _x u$ because the
quantity  $\| \nabla _xu(t)\|_{\infty }$ is not considered in the references in question.
It is perhaps not astonishing that the assumptions in these references are restrictive.
After all, an algebraic rate of decay of $\|u(t)\|_{\infty }$ is a rather strong stability result.
Too strong for our purposes, we think, because it describes the behaviour of $|u(x,t)|$
for $t\to \infty $ at any  point $x \in \overline{ \Omega }^c$, whereas we are interested
in the asymptotics of $u(x,t)$ only at points $x$ with $|x|$ large.

We will show that any $L^2$-strong solution (see (\ref{8.50}), (\ref{8.60}))
to (\ref{180}), (\ref{20b}), (\ref{30}) satisfies the estimate
\begin{eqnarray} \label{131}
| \partial ^{\alpha }_xu(x,t)| \le \mbox{$\mathfrak C$}
\, \bigl(\, |x|\, \nu (x) \,\bigr) ^{(-1-| \alpha |/2)\,(1-\epsilon )}
  \, X(t)^{ \epsilon }
\end{eqnarray}
for
$
x,\, t,\, \alpha
$
as in (\ref{130}) and for
$
\epsilon \in [0,1],
$
where $X:(0, \infty ) \mapsto (0, \infty )$
is a bounded function with $X(t)\downarrow 0$ for $t\uparrow \infty $
(Theorem \ref{theoremT8.50}).
References on existence of $L^2$-strong solutions
to  (\ref{180}), (\ref{20b}), (\ref{30})
are listed in the passage following (\ref{8.60}).
We will take the point of view that such a solution is given.
Under this assumption, (\ref{131})
may be shown without any smallness condition.
A key role in our proofs
will be played by the property $\nabla _x u \in L^2 \bigl(\, \overline{ \Omega }^c
\times (0, \infty ) \,\bigr) ^3$ verified by the solutions under consideration.

Concerning algebraic decay of $| \partial ^{\alpha }_xu(x,t)|$ with respect to $t$,
we consider the somewhat simpler system (\ref{100}) (Navier-Stokes system with Oseen term),
again with the side conditions (\ref{20b}) and (\ref{30}).
$L^2$-strong solutions to this problem fulfilling the additional assumption
\begin{eqnarray} \label{8.30}
\| \nabla _xu(t)\|_2\le c\, t^{-\kappa _1} \;\; ( t \in (1, \infty ) )
\end{eqnarray}
for some constants $c,\, \kappa _1>0$
will be shown to satisfy an inequality that in the best possible case,
arising if $f$ and $a$ are bounded and with compact support and $b$ decays sufficiently fast,
takes the form
\begin{eqnarray} \label{132}
| \partial ^{\alpha }_xu(x,t)| \le \mbox{$\mathfrak C$} \, \bigl[\,
\bigl(\, |x|\, \nu (x) \,\bigr) ^{-1-| \alpha |/2}\, (1+t)^{-\zeta }
+
\bigl(\, |x|\, \nu (x) \,\bigr) ^{(-1-| \alpha |/2)\,(1-\epsilon )}
  \, (1+t)^{\widetilde{ \zeta }\, \epsilon }\,\bigr]
\end{eqnarray}
for $x,\,t,\, \alpha $ as in (\ref{130}) and $\epsilon \in [0,1]$,
where $\zeta $ and $\widetilde{ \zeta }$ are constants determined by the data
and the exponent $\kappa _1$ in (\ref{8.30}) (Theorem \ref{theoremT8.50}).
As in the case of (\ref{131}), no smallness condition is involved in the proof of
this estimate.

The interest of our theory may be illustrated by an existence result due to
Neustupa \cite{Neu2016}. According to \cite[Theorem 4.1]{Neu2016}
(also see \cite[Theorem 1]{Neu2009}), if $f$ and $b$ vanish, $\Omega $ is smoothly bounded,
$a$ is small with respect to the norm of $H^1( \overline{ \Omega }^c)^3$,
and all eigenvalues of a certain linear operator have negative real part and stay away
from the imaginary axis, then an $L^2$-strong solution to
(\ref{180}), (\ref{20b}), (\ref{30}) in the sense of (\ref{8.50}), (\ref{8.60})
exists. Since no smallness of $U$ (velocity part of a solution to (\ref{170}), (\ref{170b})) is required
and the viscosity is arbitrary, the results in
\cite{Ma}, \cite{Hey2} or \cite{EnSh2} cannot be applied to this solution.
The theory presented here, however,
yields that Neustupa's solution satisfies (\ref{131}) (see Theorem \ref{theoremT8.50}), provided, of course, that
the spatial asymptotics of $a$ are compatible with our assumptions, listed at the beginning of section 8.

Property (\ref{8.30}) is fulfilled by the solutions to
(\ref{180}), (\ref{20b}), (\ref{30})
constructed in \cite{Ma} and \cite{Hey2}; see
\cite[inequality (7)]{Ma}, \cite[p. 675]{Hey2}. It it true that as mentioned above, these references
additionally provide algebraic decay of the $L ^{ \infty } $-norm of $u(t)$. However,
this latter property is suspended on $H^2$-regularity of the Stokes operator
(\cite[p. 323; p. 299, Proposition 1]{Ma}, \cite[p. 675]{Hey2}),
and thus on smoothness of $\partial \Omega $.
Therefore the $L ^{ \infty } $-estimates of $u(t)$ from \cite{Ma} or \cite{Hey2} cannot be used in the proof
of (\ref{132}) if $\Omega $ is supposed to be only Lipschitz bounded.
And in any case, they do not yield an access to (\ref{132}) if $| \alpha |=1$.
So, as far as we know, inequality (\ref{132}) is new at least in the case $| \alpha |=1$,
even though it only relates to (\ref{100}) instead of (\ref{180}).

In order to prove our results, we will start from a representation formula
established in \cite{DeJMFM}
for
solutions to the linear problem (\ref{10}), (\ref{20b}), (\ref{30}), and stated as equation (\ref{T10.60.20})
below, in Theorem \ref{theoremT10.60}.
This formula consists of a sum involving two volume potentials -- one on $\mathbb{R}^3 \times (0, \infty ) $
and related to $f$, the other one on $\mathbb{R}^3 $ and linked to the initial data $a$ --,
as well as a single layer potential on $S_{\infty }:= \partial \Omega \times  (0, \infty ) $
whose weight function solves an integral equation on $S_{\infty }$ (equation (\ref{T10.60.10})).
We refer to section 4 for the definition of these potential functions. 

In order to solve the integral equation (\ref{T10.60.10}), we use an $L^2$-theory
developed by Shen \cite{Shen} for the Stokes system, and extended to the Oseen system
in \cite{DeEst}.
In the framework of this theory, the right-hand side of (\ref{T10.60.10}) must belong
to a space whose definition is rather complicated and thus gives rise to much of the
technicalities we have to grapple with in what follows. This space, denoted by $H_{\infty }$
in this work, is introduced in section 3. It is involved in the crucial part of
our argument, that is, in determining how the $L^2$-norm on
$S_{T, \infty }:=\partial \Omega \times (T, \infty )$
of the solution to (\ref{T10.60.10}) is bounded in terms of $T$. This point is settled
in Theorem \ref{theoremT3.10}.

In each of the sections 5, 6 and 7, we consider one of the three potential functions appearing
in the representation formula (\ref{T10.60.20}),
deriving a pointwise decay estimate in space and in time for the function in question, among other results.
Theorem \ref{theoremT3.10} is applied in this context in order to deal with the single layer potential
from (\ref{T10.60.20}) (Corollary \ref{corollaryC5.20}).
Once upper bounds of these potentials are available, the formula in (\ref{T10.60.20}) yields
an estimate of the solution to (\ref{10}), (\ref{20b}), (\ref{30}) (Theorem \ref{theoremT8.10}).

As concerns the nonlinear problems (\ref{180}), (\ref{20b}), (\ref{30}) and (\ref{100}), (\ref{20b}),
(\ref{30}), the idea is, of course, to replace $f$ by
$
f- \tau \, (u \cdot  \nabla_x  ) u
- \tau \, (U \cdot  \nabla_x  ) u
-\tau \, (u \cdot  \nabla ) U
$
and
$
f- \tau \, (u \cdot  \nabla_x  ) u
$,
respectively, and then apply our estimates of solutions to the linear problem (\ref{10}), (\ref{20b}), (\ref{30}).
However, we were not able to shift all difficulties into the theory of this linear problem.
In fact, we will need an intermediate result from \cite{DeSIMA2} -- stated as the first estimate
in Theorem \ref{theoremT8.40} -- whose proof exploits the interaction between the nonlinearity
and the kernel function of one of the potentials 
in (\ref{T10.60.20}).

It is mainly due to the integral equation (\ref{T10.60.10}) that we cannot deal with
weak solutions to (\ref{180}) or (\ref{100}).
For this type of solution, we are not able to show that the
right-hand side of (\ref{T10.60.10}) with $f$ replaced
in the way just mentioned
belongs to the function space $H_{\infty }$, as required by Theorem \ref{theoremT10.40},
on which the resolution theory of (\ref{T10.60.10}) is based. There are other aspects
of our results whose scope is limited by this resolution theory. For example,
since we may solve this equation only in an $L^2$-framework, but not in $L^p$ with
$p\ne 2$, we cannot admit values $\zeta \ge 1$ in (\ref{129}), nor
can we obtain an algebraic decay rate of the time variable in (\ref{131}).

Let us mention some further papers related to the work at hand. Knightly \cite{Kni2}
considered pointwise decay in space of strong solutions
to the nonlinear system (\ref{180}),
detecting  the wake phenomen, but he required various smallness conditions on the data
and restrictive assumptions on the asymptotics of the solution.
Mizumachi \cite{Miz} studied the spatial asymptotics of strong solutions of (\ref{100}), (\ref{20b}), (\ref{30}),
but still under rather restrictive assumptions.
The results of these two authors
were improved in the articles \cite{DeSIMA1}, \cite{DeDCDS-A} (linear case)
and \cite{DeSIMA2} (nonlinear problem (\ref{180}), (\ref{20b}), (\ref{30})),
with predecessor papers
\cite{De1}, \cite{De2}, \cite{De3}, \cite{DeEst}, \cite{DeJMFM}, \cite{De7}.
As concerns temporal decay of spatial $L^p$-norms of solutions to the Oseen system (\ref{10})
under side conditions (\ref{20b}), (\ref{30}) and with $ f=0$ and $b=0$,
a basic study is due to Kobayashi, Shibata \cite{KoSh}. Their theory was extended in various respects
in \cite{EnSh1} and \cite{EnSh2}. The $L ^{ \infty } $-estimate from \cite{EnSh2} mentioned above
is an example of such an extension.
A different approach was used
by Bae, Jin \cite{BJ}, who considered temporal decay of weighted $L^p$-norms of solutions
to (\ref{10}), (\ref{20b}), (\ref{30}) with $b=0, \;f=0$, where the
weight functions take account of the wake phenomenon.
This type of result was extended to the nonlinear problem (\ref{180}), (\ref{20b}), (\ref{30})
with $b=0,\; f=0$ by Bae, Roh \cite{BR}.

\section{Notation. Various auxiliary results.}
\setcounter{equation}{0}

As we may recall, the bounded open set $\Omega \subset \mathbb{R}^3 $
with connected Lipschitz boundary $\partial \Omega $ and the parameter
$\tau \in (0, \infty ) $ were fixed at the beginning of section 1.
We will write $n ^{( \Omega )} $ for the outward unit normal to $\Omega $.
The notations $S_{\infty }:=\partial  \Omega \times (0, \infty )$ for $T \in (0, \infty ]$
and $S_{T, \infty }:= \partial \Omega \times (T, \infty ) $ for $T \in [0, \infty )$
were also already introduced in section 1, as was the function $\nu$ (see (\ref{101})),
as well as the abbreviation $| \alpha |$ for the length $\alpha _1+\alpha _2+\alpha _3$
of a multiindex $\alpha \in \mathbb{N} _0^3$. The symbol $| \;|$ additionally denotes
the Euclidean norm in $\mathbb{R}^3$.

For $A \subset \mathbb{R}^3 $, we set $A^c:=\mathbb{R}^3 \backslash A$.
Moreover we abbreviate $e_1:=(1,0,0)$, and we put
$B_r(x):=\{y \in \mathbb{R}^3 \, :\;|y-x|<r\}$
for $x \in \mathbb{R}^3 $, and $B_r:=B_r(0)$, where $r \in (0, \infty ) $.
Let $A$ be a nonempty set. If $\varphi :A \mapsto \mathbb{R} $ is a function, we define
$| \varphi |_{\infty }:=\sup\{| \varphi (x)|\, :\, x \in A\}.$
Let $n \in \mathbb{N} $ and $B$ a vector space consisting of functions $f: A \mapsto \mathbb{R} .$
Suppose $B$ is equipped with a norm, denoted by $\|\;\|_B$.
Then we put
$
B^n:=\{F:A \mapsto \mathbb{R} ^n\, :\; F_j \in B\; \mbox{for}\; 1\le j\le n\},
$
and we equip $B^n$ with the norm
$
\|F\|_B ^{(n)} := (\sum_{j=1}^n\|F_j\|_B^2) ^{1/2} \; (F \in B^n)
$.
But instead of $\|\;\|_B ^{(n)}$, we will write $\|\;\|_B$ again.

Next we introduce a fractional derivative. Let
$A \subset \mathbb{R}^3 ,\; T \in (0, \infty ]$
and
$\psi:A \times (0,T) \mapsto  \mathbb{R} $
a function such that
$\psi (x,\, \cdot \,)$ is measurable and
$
\int_{ 0}^t(t-r) ^{-1/2} \,|\psi (x,r)|\, dr < \infty
$
for
$x \in A,\; t \in (0,T) $.
Define $W(x,t):=  \int_{ 0}^t (t-r) ^{-1/2} \, \psi(x,r)\; dr$
for these $x $ and $ t$.
If the derivative
$\partial _tW(x,t)$
exists for some such $x$ and $t$,
we put
$\partial _t ^{1/2} V(x,t):= \Gamma (1/2) ^{-1} \, \partial _tW(x,t).$
Here $\Gamma $ denotes the usual Gamma function. In the case that
$\partial _t W(x,t)$ exists for any $x \in A,\; t \in (0,T), $
we define $\partial _4 ^{1/2}V := \Gamma (1/2) ^{-1} \,\partial _4W.$

Let $A \subset \mathbb{R}^3 $ be open. For $p \in [1, \infty ],$
the norm of the Lebesgue space $L^p(A)$, defined with respect to the
Lebesgue measure on $R^3$, is denoted by $\|\;\|_p$.
The same notation is used for the norm of $L^p$-spaces on $\partial \Omega $ or
on subsets of $S_{\infty }$.
If $T \in (0, \infty ],$ put
\begin{eqnarray} \label{2.10}
L^2_n(S_{T }):=\{\psi \in L^2( S_{T})^3\, :\, \int_{ \partial \Omega }
n ^{( \Omega )} (x) \cdot \psi (x,t)\, do_x =0\; \mbox{for a. e.}\; t \in (0,T)\}.
\end{eqnarray}
For $a,b \in \mathbb{R} \cup \{\infty \},\; a<b,$
let $L^1_{loc}\bigl(\, [a,b) \,\bigr) $ stand for the set of all functions $g:(a,b) \mapsto \mathbb{R} $
such that $g|(a,c) \in L^1 \bigl(\, (a,c) \,\bigr) $ for $c \in (a,b)$.
For $p \in [1, \infty ),$ we write $W^{1,p}(A)$ for the usual Sobolev space of order $1$ and
exponent $p$. If $p=2$, we use the notation $H^1(A)$ instead of $W^{1,2}(A).$
For $s \in (0,1)$, let $H^s(A)$ be the Sobolev space defined via the intrinsic norm with
exponent 2 introduced in \cite[section 7.51]{Adams}. The symbol $H^s_{\sigma }(A) $
stands for the closure of the set $\{V \in C_0 ^{ \infty } (A)^3\, :\, \mbox{div}\, V =0\}$
with respect to the norm of $H^s(A)^3\; (s \in (0,1])$.
Moreover, if $A$ is again an open set in $\mathbb{R}^3 $,
then $W^{1,1}_{loc}(A)$ designates the set of all functions $V:A \mapsto \mathbb{R} $ such that
$V|K \in W^{1,1}(K)$ for any open, bounded set $K \subset \mathbb{R}^3 $ with $\overline{ K} \subset A$.
The Sobolev space $H^1( \partial \Omega )$ is to be defined in the standard way (see \cite[section III.6]{FJK},
for example). We write $\|\;\|_{H^1( \partial \Omega )^3}$ for the norm of $H^1( \partial \Omega )^3$,
and $\|\;\|_{H^1( \partial \Omega ) ^{\prime} }$ for the usual norm of the canonical dual space
$H^1( \partial \Omega ) ^{\prime} $ of $H^1( \partial \Omega )$. Frequently we will use the fact
that there is $c>0$ such that for functions $V \in C^1( \mathbb{R}^3 )^3,$
the inequality $\| V| \partial \Omega \|_{H^1( \partial \Omega )^3 }\le c \, (\|V| \partial \Omega \|_2+
\sum_{j = 1}^3 \| \partial _jV| \partial \Omega \|_2)$
holds.

If $J \subset \mathbb{R} $ is an interval, $A \subset \mathbb{R}^3 $ is open and $v:A \times J \mapsto \mathbb{R} $
is a function with suitable smoothness, then the notation $\Delta _xv,\; \nabla _xv,\; \mbox{div}_xv$
indicates that the differential operators in question refer to $v(x,t)$ as a function of $x \in A$. It
will be convenient to denote this function by $v(t)$, for $t \in J$. For a function $V:A \mapsto \mathbb{R} $,
we write $\Delta V,\; \nabla V$ and $\mbox{div}\, V$, respectively.

Let $B$ be a Banach space, $a,b \in \mathbb{R} \cup \{\infty \}$ with $a<b$, and $p \in [1, \infty ]$.
Then the norm of the space $L^p(a,b,B)$ is denoted by $\|\;\|_{L^p(a,b,B)}.$
If $B=L^q(A)^n$ for some $q \in [1, \infty ],\; n \in \{1,\, 3\},\; A \subset \mathbb{R}^3 $ open
or $A=\partial \Omega $, then we use the notation $\|\;\|_{q,p; b}$ instead of
$\|\;\|_{L^p \bigl(\, a,b,\, L^q(A)^n \,\bigr) }.$
We consider the spaces $L^p(a,b,\, L^p(A)^n)$ and $L^p \bigl(\, A \times (a,b) \,\bigr)^n $
as identical, and denote their norm by $\|\;\|_p$.
The term $L^p \bigl(\, [a,b),\, B \,\bigr) $ stands for the set of all functions $F:(a,b) \mapsto \mathbb{R} $
such that $F| (a,c) \in L^p(a,c,B)$ for any $c \in (a,b)$.

Let $T \in  (0, \infty ].$ For any $w \in S_T$, we introduce a mapping $F_w \in L^2 \bigl(\, 0,T,\,
H^1( \partial \Omega ) ^{\prime} \,\bigr) $ by setting
$F_w(t)(V):=\int_{ \partial \Omega } w(x,t)\, V(x)\, do_x$ for $V \in H^1( \partial \Omega )$ and for
a. e. $t \in (0,T)$. However, we will again write $w$ instead of $F_w$.

We write $C$ for numerical constants, and $C( \gamma _1,\, ...,\, \gamma _n)$ for constants
depending exclusively on parameters $\gamma _1,\, ...,\, \gamma _n  \in [0, \infty ) ,$
for some $n \in \mathbb{N} $. However, it will not be possible to specify all our constants
in such a precise way. So in most cases we will use a different symbol, namely \mbox{$\mathfrak C$},
for generic constants, assuming that their dependencies become clear from context. Occasionally
we will use expressions of the form $\mbox{$\mathfrak C$} (\gamma _1,\, ...,\, \gamma _n)$
in order to insist that the constant under consideration depends on
$\gamma _1,\, ...,\, \gamma _n \in  [0, \infty ) ,$ but it may additionally be a function of other quantities.

We will frequently use Minkowski's inequality for integrals. For the convenience of the reader, we
state a suitable version as
\begin{theorem}[\mbox{\cite[p. 271, Appendix A1]{Stein}}] \label{theoremT1.10}
Let $A_1,\, A_2$ be nonempty sets, $\mbox{$\mathfrak A$} _j$ a measure space on $A_j$ and $m_j$ a
$\sigma $-finite measure on $\mbox{$\mathfrak A$} _j$, for $j \in \{1,\, 2\} .$
Let $F:A_1 \times A_2 \mapsto \mathbb{R} $ be an $\mbox{$\mathfrak A$} _1 \otimes \mbox{$\mathfrak A$} _2$-measurable
function, and take $p \in [1, \infty ). $ Then
\begin{eqnarray*}
\Bigl( \int_{ A_1 }\Bigl( \int_{ A_2}|F(x,y)|\, dm_2(y) \Bigr)^p\, dm_1(x) \Bigr) ^{1/p}
\le
\int_{ A_2} \Bigl( \int_{ A_1}|F(x,y)|^p\, dm_1(x) \Bigr) ^{1/p}\, dm_2(y).
\end{eqnarray*}
\end{theorem}
Next we state a Sobolev inequality for certain functions in exterior domains.
\begin{theorem} \label{theoremT1.101}
Let $V \in W^{1,1}_{loc}( \overline{ \Omega }^c)$
with $V \in L ^{\kappa } ( \overline{ \Omega }^c)$
for some $\kappa \in [1, \infty )$
and $\nabla V \in L^2( \overline{ \Omega  }^c)^3.$
Then $V \in L^6( \overline{ \Omega }^c)$
and $\|V\|_6\le \mbox{$\mathfrak C$} \,\| \nabla V\|_2.$
\end{theorem}
{\bf Proof:}  \cite[Lemma 2.4]{DeDCDS-A}, which is a consequence of \cite[Theorem II.5.1]{Ga1}.
\hfill $\Box $

\vspace{1ex}
We note some technical details whose sense will become clear later on.
\begin{lemma} \label{lemmaL1.10}
Let $\epsilon \in (0,1)$. Put
$\varphi ( \epsilon ):= \epsilon -1/2$
if $\epsilon > 1/2$, and
$\varphi ( \epsilon ):=\min\{1/12,\; \epsilon /4\}$
if $\epsilon \le 1/2.$
Let $k \in \mathbb{N} $ with $k\, \epsilon \le 1< (k+1)\, \epsilon ,$
and let $j \in \{0,\, ...,\, k-1\}.$
Then
$$Z(j):= -1/2 +j/(2\, k) +1/k -(j+1)\, \epsilon /k\:\le\: -\varphi ( \epsilon ).$$
\end{lemma}
{\bf Proof:}
Suppose that $\epsilon > 1/2$. Then $k=1$ so that $j=0$. It follows that
$Z(j)=-1/2+1- \epsilon \le -\varphi ( \epsilon )$.
Next suppose that $\epsilon \in (1/3,\, 1/2],$ so $k=2$ and $j \in \{0,\, 1\} .$
If $j=0$, we have $Z(j)=-\epsilon /2\le -\varphi ( \epsilon ),$
and if $j=1,\; Z(j)=1/4-\epsilon \le -1/12\le - \varphi ( \epsilon ).$
Next consider the case $\epsilon \in (1/4,\, 1/3],$
so that $k=3$ and $j \in \{0,\,1,\, 2\} .$ If $j=0$, we find
$
Z(j)=-1/6-\epsilon /3\le -1/6\le -\varphi ( \epsilon ).
$
For $j=1$, we get
$Z(j)=-2\, \epsilon /3\le -1/6\le -\varphi ( \epsilon ),$
and for $j=2$,
$Z(j)=1/6 - \epsilon \le -1/12 \le -\varphi ( \epsilon )$.

Now we turn to the case $\epsilon \le 1/4$, which means that $k\ge 4.$
If $j \in [k/2,\, k-2]$, then
$Z(j)\le -1/2 +(k-2)/(2\, k) +1/k - \epsilon /2 = -\epsilon /2\le -\varphi ( \epsilon ).$
The value $j=k-1$ leads to the equation $Z(j)=1/(2\, k)-\epsilon ,$
hence by the choice of $k$,
$$
Z(j)
=
\bigl[\,(k+1)/k \,\bigr]\, \bigl(\, 2\, (k+1) \,\bigr) ^{-1} - \epsilon
\le
\bigl[\,(k+1)/k \,\bigr]\, \epsilon /2 - \epsilon
=
-\epsilon /2 +\epsilon /(2\, k)
\le
-3\, \epsilon /8,
$$
so that $Z(j)\le -\varphi ( \epsilon ).$
Next suppose that $k/4\le j\le k/2$ so that
\begin{eqnarray*}
Z(j)
\le
-1/2+1/4+1/k-j\, \epsilon /k
\le
-1/4+1/k -\epsilon /4
\le
-\epsilon /4
\le
-\varphi ( \epsilon ),
\end{eqnarray*}
where the second from last inequality holds because $k\ge 4$. If $j\le k/4$, we get
$
Z(j)\le -1/2+1/8+1/k\le -1/8\le -\varphi ( \epsilon ),
$
with the second from last inequality being a consequence of the relation $k\ge 4$.
Thus we have found in any case that $Z(j)\le -\varphi ( \epsilon ).$
\hfill $\Box $
\begin{lemma} \label{lemmaL1.20}
Let $T>0,\; n \in \mathbb{N} ,\;n\ge 2,\;
\phi \in L^2 \bigl(\, \partial \Omega \times (T/2,\, T) \,\bigr) ^3.$
For $i \in \{0,\, ...,\, n-1\} ,$ put
$I_i:= \bigl(\,  T\,(1+i/n)/2,\; T\, \bigl(\, 1+(i+1)/n \,\bigr) /2 \,\bigr) .$

Then there is $i_0 \in \{0,\, ...,\, n-1\} $ such that
$\| \phi  | \partial \Omega \times I_{i_0}\|_2 \le \| \phi \|_2\, n ^{-1/2} .$
\end{lemma}
{\bf Proof:} Suppose that
$\| \phi | \partial \Omega \times I_{i}\|_2 > \| \phi \|_2\, n ^{-1/2}$
for any 
$i \in \{0,\, ...,\, n-1\} $.
Then
\begin{eqnarray*}
\| \phi \|_2^2
=
\Bigl( \sum_{i = 0}^{n-1}\| \phi | \partial \Omega \times I_i\|_2^2 \Bigr) ^{1/2}
>
\| \phi \|_2\, n ^{-1/2} \, \Bigl( \sum_{i = 0}^{n-1}1 \Bigr) ^{1/2}
=
\| \phi \|_2.
\end{eqnarray*}
Since this is a contradiction, the lemma is proved.
\hfill $\Box $
\begin{lemma} \label{lemmaL10.10}
Let $T \in (0, \infty ) ,\; \phi \in L^1_{loc} \bigl(\, [0, \infty ) \,\bigr) $
with $\phi |(T, \infty ) \in C^1 \bigl(\, (T, \infty ) \,\bigr) .$
Define
$
H(t):= \int_{ 0}^t(t-r) ^{-1/2} \,\phi (r)\, dr
$
for $t \in (T, \infty ).$
Then $H \in C^1 \bigl(\, (T, \infty ) \,\bigr) $ and
\begin{eqnarray} \label{L10.10.10}&&
H ^{\prime} (t)
=
\sqrt{2}\, (t-T) ^{-1/2} \, \phi \bigl(\, (t+T)/2 \,\bigr)
- (1/2)\, \int_{ 0}^{(t+T)/2}(t-r) ^{-3/2} \, \phi (r)\, dr
\\&& \nonumber \hspace{3em}
+
\int_{(t+T)/2}^t(t-r) ^{-1/2} \, \phi ^{\prime} (r)\, dr
\quad \mbox{for}\;\; t \in (T, \infty ),
\end{eqnarray}
with all the preceding integrals existing as Lebesgue integrals.
\end{lemma}
{\bf Proof:} Let $T_1,\, T_2 \in (T, \infty )$ with $T_1<T_2$, and put
$T_0:=(T+T_1)/2,$
\begin{eqnarray*} 
F(t):= \int_{ 0}^{T_0}(t-r) ^{-1/2} \, \phi (r)\, dr,
\quad 
G(t):= \int_{ T_0}^{t}(t-r) ^{-1/2} \, \phi (r)\, dr
\quad \mbox{for}\;\; 
t \in [T_1, T_2],
\end{eqnarray*}
hence
$H|[T_1,T_2] =F+G.$
Lebesgue's theorem and the relation $T_0<T_1$ yield that $F \in C^1 ( [T_1,T_2] ) .$
Turning to $G$, we choose a function $\widetilde{ \phi}\in C^1( \mathbb{R} )$ such that
$\widetilde{ \phi }|[T_0, \infty ) = \phi |[T_0, \infty ),$
and then consider
\begin{eqnarray*} 
\widetilde{ G}( \kappa ,t):=\int_{ 0}^{\kappa } r ^{-1/2} \, \widetilde{ \phi}(t-r)\, dr 
\quad \mbox{for}\;\; 
(\kappa ,t)  \in B:= [T_1-T_0,\, T_2-T_0] \times  [T_1,T_2].
\end{eqnarray*}
Again by Lebesgue's theorem, the derivative $\partial _2 \widetilde{ G} $ exists and
belongs to $C^0(B)$. Moreover, since for any $t \in [T_1,T_2]$,
the function $r \mapsto r ^{-1/2} \, \widetilde{ \phi}(t-r),\; r \in [T_1-T_0, T_2-T_0],$
is continuous, the derivative $\partial _1 \widetilde{ G}$ exists, too, and
$\partial _1\widetilde{ G} \in C^0(B)$. But $\widetilde{ G}(t-T_0,\, t)=G(t)$ for $t \in [T_1,T_2],$
so we get $H|[T_1,T _2] =F+G \in C^1 ([T_1,T_2]).$

Next we differentiate the function $\widetilde{ G }(t-T_0,\, t),\; t \in [T_1,T_2],$
write the result as a sum of $(t-T_0) ^{-1/2} \, \phi (T_0)$
plus an integral from $T_0$ to $t$, split off an integral from $T_0$ to $(t+T)/2$
and integrate by parts in that latter integral. Equation (\ref{L10.10.10}) follows by
a differentiation of $F$ and because
$H|[T_1,T _2] =F+G$. \hfill $\Box $
\begin{corollary} \label{corollaryC10.10}
Let $T \in (0, \infty ) ,\; \Psi  \in L^1_{loc} \bigl(\, [0, \infty ),\, L^2( \partial \Omega )^3 \,\bigr) ,\;
v \in C^1 \bigl(\, \mathbb{R}^3 \times (T, \infty ) \,\bigr) ^3$
with $\Psi |S_{T, \infty }=v|S_{T, \infty }.$

Then there is a set $N \subset \partial \Omega $ of measure zero such that for any $x \in \partial \Omega
\backslash N,\; t \in (T, \infty ),$ we have
$\int_{ 0}^t(t-r) ^{-1/2} \,|\Psi (x,r)|\, dr < \infty .$
Define
$W(x,t):=\int_{ 0}^t(t-r) ^{-1/2} \,\Psi (x,r)\, dr$
for such $x$ and $t$.
Then, for $x \in \partial \Omega \backslash N$,
we have $W(x,\, \cdot \, ) \in C^1 \bigl(\, (T, \infty ) \,\bigr) ^3,$
and equation (\ref{L10.10.10}) holds with $\Psi (x,\, \cdot \, )$ in the role of $\phi$, with all
integrals existing in the Lebesgue sense. In particular, the fractional derivative
$\partial ^{1/2} _t\Psi (x,t)$ exists for $x, \, t$ as above.
\end{corollary}
{\bf Proof:}
Let $T_0 \in (T, \infty ).$ Then $\Psi |S_{T_0} \in L^1 \bigl(\, 0,T_0,\, L^2 ( \partial \Omega )^3 \,\bigr) $
by our assumptions on $\Psi $, so there is a set $N \subset \partial \Omega $ of measure zero such that
$\Psi (x,\, \cdot \, )|(0,T_0) \in L^1 \bigl(\, (0,T_0) \,\bigr) ^3$ for $x \in \partial \Omega \backslash N$.
Since $\Psi (x,\, \cdot \, )|(T, \infty ) $ is a $C^1$-function for any $x \in \partial \Omega $,
we thus get that $\Psi (x,\, \cdot \, ) \in L^1_{loc} \bigl(\, [0, \infty ) \,\bigr) ^3$
for $x \in \partial \Omega \backslash N$.
Now the corollary follows from Lemma \ref{lemmaL10.10}.
\hfill $\Box $
\vspace{1ex}

In order to handle the term $\nu(x)$ defined in (\ref{101}),
we need the ensuing three lemmas, all of them well known.
\begin{lemma}[\mbox{\cite[Lemma 4.8]{DeKr2}}] \label{lemmaL10.1}
$\nu (x-y) ^{-1} \le C\, |y|\, \nu(x) ^{-1} $ for $x,y \in \mathbb{R}^3 $.
\end{lemma}
\begin{lemma}[\mbox{\cite[Lemma 2.3]{Farwig}}] \label{lemmaL10.2}
Let $\beta \in (1, \infty ).$ Then
$\int_{ \partial B_r}\nu (x) ^{-\beta }\; do_x\le C( \beta )\, r$ for $r \in (0, \infty ) $.
\end{lemma}
\begin{lemma}[\mbox{\cite[Lemma 2]{De1}}] \label{lemmaL10.3}
Let $K \in (0, \infty ) $. Then for $x \in \mathbb{R}^3 ,\; t \in (0, \infty ) ,$
\begin{eqnarray*}
|x- \tau \, t \, e_1|^2+t \ge C(K, \tau )\, \bigl(\, \chi_{[0,K]}(|x|)\, (|x|^2+t)
+
\chi_{(K, \infty )}(|x|)\, (|x|\, \nu (x) +t) \,\bigr) .
\end{eqnarray*}
\end{lemma}
\begin{corollary} \label{corollaryC1.10}
Let $ \beta \in (2, \infty ).$ Then
$
\int_{ B_R^c}\bigl(\, |x|\,\nu (x) \,\bigr) ^{- \beta }\, do_x
\le
C( \beta )\, R^{-\beta +2}
$
for $R \in (0, \infty ) $.
\end{corollary}
{\bf Proof:}
By Lemma \ref{lemmaL10.2},
$$
\int_{ B_R^c}\bigl(\, |x|\, \nu (x) \,\bigr) ^{-\beta }\, do_x
=
\int_{ R} ^{ \infty } r ^{-\beta }\, \int_{ \partial B_r}\nu (x)^{- \beta }\, do_x\, dr
\le
C( \beta )\, \int_{ R} ^{ \infty } r^{-\beta +1}\, dr
=
C( \beta )\, R^{-\beta +2}.
$$
\hfill $\Box $

\section{Definition of the space $H_{\infty }$ and its norm. Some properties of this space.}
\setcounter{equation}{0}

In this section, we introduce the key function space of our theory, denoted by $H_{\infty }$
and taken from \cite{Shen}. We start by fixing $T_0 \in (0, \infty ]$ and defining
\begin{eqnarray*}
\widetilde{ H}_{T_0}:= \{ w|S_{T_0}\, :\, w \in C ^{ \infty } _0( \mathbb{R}^4)^3,\;
w| \mathbb{R}^3 \times (-\infty , 0) =0\}.
\end{eqnarray*}
For $\Psi \in \widetilde{ H}_{T_0}$,
we set
\begin{eqnarray*}
\|\Psi \|_{H_{T_0}}:= \biggl(\,  \int^{ T_0 }_0 \Bigl( \|\Psi (t)\|_{H^1( \partial \Omega )^3}^2
+
\| \partial ^{1/2} _4\Psi (t)\|_{2}^2
+
\| n^{ (\Omega )}\cdot  \partial _4\Psi (t)\|_{H^1( \partial \Omega ) ^{\prime} }^2 \Bigr) \, dt \,\biggr) ^{1/2} .
\end{eqnarray*}
The mapping $\|\; \|_{H_{T_0}}$ is a norm on  $\widetilde{ H}_{T_0}$.
In the case $T_0=\infty $, we further put
$$
\|\Psi |S_{T, \infty  }\|_{H_{T, \infty }}:= \biggl(\,  \int^{ \infty }_T \Bigl( \|\Psi (t)\|_{H^1( \partial \Omega )^3}^2
+
\| \partial ^{1/2} _4\Psi (t)\|_{2}^2
+
\| n^{ (\Omega )} \cdot  \partial _4\Psi (t)\|_{H^1( \partial \Omega ) ^{\prime} }^2 \Bigr) \, dt \,\biggr) ^{1/2} 
$$
for $T \in [0, \infty ),\; \Psi \in \widetilde{ H}_{\infty }$.
Note that the term $\|\Psi |S_{T, \infty  }\|_{H_{T, \infty }}$ depends on
$\Psi |S_T$, via the fractional derivative $ \partial ^{1/2} _4\Psi$. Further note
that in the case $T=0$, we have
$\|\Psi |S_{T, \infty  }\|_{H_{T, \infty }}=\|\Psi\|_{H_{\infty }}.$

In order to construct a completion of $\widetilde{ H}_{T_0}$, we consider sequences
$(\Psi_n)$ in $\widetilde{ H}_{T_0}$ which are Cauchy sequences with respect to the norm
$\|\;\|_{H_{T_0}}$, with the additional property that there is $\Psi \in L^2_n( S_{T_0})$
with $\|\Psi_n-\Psi\|_{2}\to 0$. In such a situation it is not clear a priori how $\Psi$ is related
to the limit of $(\partial ^{1/2} _4\Psi_n) $ in $L^2(S_{T_0})$
and to the limit of $( n^{ (\Omega )} \cdot  \partial _4\Psi_n)$
in $L^2 \bigl(\, 0,T_0,\, H^1( \partial \Omega ) ^{\prime} ).$
But it turns out that if $\Psi=0$, these limits vanish as well. This fact and some other, related ones
are the subject of
\begin{lemma} \label{lemmaL10.30}
Let $T_0 \in (0, \infty ],\; \Psi \in L^2(S_{T_0})^3,\; v_n \in C ^{ \infty } _0( \mathbb{R} ^4)^3$
with $v_n| \mathbb{R}^3 \times (-\infty , 0) =0$ for $n \in \mathbb{N} $ such that
$\|v_n-\Psi \|_2\to 0$ and $(v_n|S_{T_0})$ is a Cauchy sequence with respect to the norm
$\|\;\|_{H_{T_0}}$ on $\widetilde{ H}_{T_0}$.

Then the sequence $(\|v_n|S_{T_0}\|_{H_{T_0}})$  converges.

If $\Psi = 0$, then $\lim\limits_{n\to \infty }\|v_n|S_{T_0}\|_{H_{T_0}}=0$.

If $T_0=\infty ,\; T \in (0, \infty ) $ and $\Psi |S_T =0$,
then
$\lim\limits_{n\to \infty }\|v_n|S_{T, \infty }\|_{H_{T, \infty }}=\lim\limits_{n\to \infty }\|v_n|S_{\infty }\|_{H_{\infty }}.$

If $T_0= \infty ,\; \widetilde{ T},\, T \in (0, \infty ) $ with $\widetilde{ T}\le T,$
then $\lim\limits_{n\to \infty }\|v_n|S_{T, \infty }\|_{H_{T, \infty }}\le
\lim\limits_{n\to \infty }\|v_n|S_{\widetilde{ T} , \infty }\|_{H_{\widetilde{ T}, \infty }}.$
\end{lemma}
{\bf Proof:}
Since $(v_n|S_{T_0})$ is a Cauchy sequence with respect to the norm $\|\;\|_{H_{T_0}}$,
the sequence $(\|v_n|S_{T_0}\|_{H_{T_0}})$ is a Cauchy sequence in $\mathbb{R} $ and thus converges.
It further follows there are functions
$\Psi ^{(1)} \in L^2 \bigl(\, 0, T_0 ,\, H^1( \partial \Omega )^3 \,\bigr) ,\; \Psi ^{(2)} \in L^2(S_{T_0})^3$
and
$\Psi ^{(3)} \in L^2 \bigl(\, 0, T_0 ,\, H^1( \partial \Omega ) ^{\prime}  \,\bigr)$
such that
\begin{eqnarray} \label{L10.30.10}
\|v_n-\Psi ^{(1)} \|_{L^2 \bigl(\, 0, T_0 ,\, H^1( \partial \Omega )^3 \,\bigr)}
+
\| \partial ^{1/2} _4v_n -\Psi ^{(2)} \|_2
+
\| \overline{ v}_n-\Psi ^{(3)} \|_{L^2 \bigl(\, 0, T_0 ,\, H^1( \partial \Omega ) ^{\prime}  \,\bigr)}
\longrightarrow 0,
\end{eqnarray}
where $\overline{ v}_n(x,t):= n^{ (\Omega )}(x) \cdot v_n(x,t)$ for $n \in \mathbb{N} ,\; x \in \partial
\Omega ,\; t \in (0,T_0).$
By the definition of the mapping $\|\;\|_{H_{T, \infty }}$,
this means in the case $T_0=\infty $ that
\begin{eqnarray} \label{L10.30.20}
\lim_{n\to \infty }\|v_n|S_{\widetilde{ T} , \infty }\|_{H_{\widetilde{ T}, \infty }}^2
=
\int^{ \infty }_{\widetilde{ T}} \Bigl( \|\Psi ^{(1)}  (t)\|_{H^1( \partial \Omega )^3}^2
+
\| \Psi ^{(2)}  (t)\|_{2}^2
+
\| \Psi ^{(3)}  (t)\|_{H^1( \partial \Omega ) ^{\prime} }^2 \Bigr) \, dt 
\end{eqnarray}
for any $\widetilde{ T} \in [0, \infty ).$ Now suppose that $T_0=\infty ,\; T \in (0, \infty ) ,$
and $\Psi |S_T =0.$ We will show that the limits
$\lim\limits_{n\to \infty }\|v_n|S_{T, \infty }\|_{H_{T, \infty }}$ and $\lim\limits_{n\to \infty }\|v_n|S_{\infty }\|_{H_{\infty }}$
coincide. Since $\|\Psi -v_n\|_2\to 0$, we obtain from (\ref{L10.30.10}) that $\Psi = \Psi ^{(1)} $.
But $\Psi |S_T =0,$ so we may conclude that $\Psi ^{(1)} |S_T =0$.
Moreover, for $V \in L^2( \partial \Omega )^3,\; \phi \in C ^{ \infty } _0 \bigl(\, (0,T) \,\bigr) ,$
we get from (\ref{L10.30.10}) that
\begin{eqnarray} \label{L10.30.30}
\int_{ \partial \Omega }\int_{ 0}^T \Psi ^{(2)} (x,t) \cdot V(x) \, \phi (t)\, dt\, do_x
=
\lim_{n\to \infty }
\int_{ \partial \Omega }\int_{ 0}^T \partial ^{1/2} _4v_n (x,t) \cdot V(x) \, \phi (t)\, dt\, do_x.
\end{eqnarray}
Put $w_n(x,t):=\int_{ 0}^t (t-r) ^{-1/2} \, v_n(x,r)\, dr$
for $n \in \mathbb{N} ,\; x \in \partial \Omega ,\; t \in (0, \infty ) $.
By Lemma \ref{lemmaL10.10}, we have $w_n(x,\, \cdot \,) \in C^1 \bigl(\, (0, \infty ) \,\bigr) ^3$
for any $x \in \partial \Omega ,\; n \in \mathbb{N} .$ Since $\partial ^{1/2} _4v_n (x,t)=\partial _tw_n(x,t)$
for $n,\,x,\, t$ as before, we get
\begin{eqnarray} \label{L10.30.35}
\int_{ \partial \Omega }\int_{ 0}^T  \partial ^{1/2} _4v_n (x,t) \cdot V(x) \, \phi (t)\, dt\, do_x
=
-
\int_{ \partial \Omega }\int_{ 0}^T w_n (x,t) \cdot V(x) \, \phi ^{\prime} (t)\, dt\, do_x
\end{eqnarray}
for $n \in \mathbb{N} $. But
\begin{eqnarray} \label{L10.30.40}&&
\Bigl| \int_{ \partial \Omega }\int_{ 0}^T w_n (x,t) \cdot V(x) \, \phi ^{\prime} (t)\, dt\, do_x \Bigr|
\\&&\nonumber
=
\Bigl| \int_{ \partial \Omega }\int_{ 0}^T \int_{ r}^T
(t-r) ^{-1/2} \, \phi ^{\prime} (t)\, dt\, v_n(x,r) \cdot V(x) \, dr\, do_x \Bigr|
\\&&\nonumber
\le 
\mbox{$\mathfrak C$} (T)\, | \phi ^{\prime} |_{\infty } \, \int_{ 0}^T \|v_n(r)| \partial \Omega \|_2\, dr
\: \|V\|_2
\le
\mbox{$\mathfrak C$} (T)\, | \phi ^{\prime} |_{\infty } \, \|V\|_2 \, \|v_n|S_T \|_2
\;\; (n \in \mathbb{N} ).
\end{eqnarray}
Since $\|v_n|S_T\|_2 =\|v_n-\Psi|S_T\|_2\; (n \in \mathbb{N} )$
and $\|\Psi -v_n\|_2\to 0$, we may conclude from (\ref{L10.30.30}) -- (\ref{L10.30.40})
that $\Psi ^{(2)} |S_T=0.$ A similar
reasoning based on (\ref{L10.30.10}) yields $\Psi ^{(3)} |S_T=0$.
Now we may conclude from (\ref{L10.30.20}) with $\widetilde{ T}=0$ and $\widetilde{ T}=T$
that
$\lim\limits_{n\to \infty }\|v_n|S_{T, \infty }\|_{H_{T, \infty }}=\lim\limits_{n\to \infty }\|v_n|S_{\infty }\|_{H_{\infty }}.$

If $T_0 \in (0, \infty ]$ and $ \Psi =0$, the same reasoning yields $\Psi ^{(i)} =0$ for 
$i \in \{1,\, 2,\, 3\} .$ In this way we deduce
from (\ref{L10.30.10}) that $\lim_{n\to \infty }\|v_n|S_{T_0}\|_{H_{T_0}}=0.$

The last claim of the lemma is an immediate consequence of (\ref{L10.30.20}).
\hfill $\Box $
\begin{corollary} \label{corollaryC10.30}
Let $T_0 \in (0, \infty ],\; \Psi \in L^2(S_T)^3,\; (\Psi_n) $ and $( \widetilde{ \Psi}_n)$
sequences in $\widetilde{ H}_{T_0}$ with
$\|\Psi-\Psi_n\|_2\to 0,\; \| \Psi -\widetilde{ \Psi}_n\|_2 \to 0$
and such that
$(\Psi_n)$ and $( \widetilde{ \Psi}_n)$ are Cauchy sequences with respect to the norm $\|\;\|_{H_{T_0}}$
on $\widetilde{ H}_{T_0}$.
Then $\lim\limits_{n\to \infty }\|\Psi_n\|_{H_{T_0}}=\lim\limits_{n\to \infty }\| \widetilde{ \Psi}_n\|_{H_{T_0}}.$
If $T_0=\infty$ and $T \in (0, \infty ) $, it further follows that
$\lim\limits_{n\to \infty }\|\Psi_n|S_{T, \infty }\|_{H_{T, \infty }}
=\lim\limits_{n\to \infty }\| \widetilde{ \Psi}_n|S_{T, \infty }\|_{H_{T, \infty }}.$
\end{corollary}
{\bf Proof:} Lemma \ref{lemmaL10.30} \hfill $\Box $
\vspace{1ex}

Let $T_0 \in (0, \infty ].$
We define the space
$H_{T_0}$ as the set of all functions $\Psi \in L^2_n(S_{T_0})$ for which there is a sequence
$(\Psi_n)$ in $\widetilde{ H}_{T_0}$ such that $\|\Psi-\Psi_n\|_2\to 0 $ and $(\Psi_n)$ is a Cauchy
sequence with respect to the norm $\|\;\|_{H_{T_0}}$. This means in particular that $H_{T_0}\subset
L^2_n(S_{T_0})$. (See (\ref{2.10}) for the definition of $L^2_n(S_{T_0})$.)

For $\Psi \in H_{T_0}$ and a sequence $(\Psi_n)$ as above, we set
$\|\Psi\|_{H_{T_0}}:=\lim\limits_{n\to \infty }\|\Psi_n\|_{H_{T_0}}.$
Corollary \ref{corollaryC10.30} implies that the mapping $\|\;\|_{H_{T_0}}$ is well defined on
$H_{T_0}$. The space $H_{T_0}$ equipped with this mapping
is a Banach space, a fact that we will not need in the following.

Suppose that $T_0=\infty $. We set
$\|\Psi |S_{T, \infty }\|_{H_{T, \infty }}:=\lim\limits_{n\to \infty }\|\Psi_n | S_{T, \infty }\|_{H_{T, \infty }}$
for $\Psi \in H _{ \infty } ,\; T \in (0, \infty ) ,$
where $\Psi_n \in \widetilde{ H}_{\infty }$ for $n \in \mathbb{N} $ with the properties that
$\|\Psi-\Psi_n\|_2 \to 0$
and $(\Psi_n)$ is a Cauchy sequence with respect to the norm $\|\;\|_{H_{\infty }}.$
It again follows from Corollary \ref{corollaryC10.30} that the mapping $\|\;\|_{H_{T, \infty} }$
is well defined.
\begin{corollary} \label{corollaryC10.40}
Let $T \in (0, \infty ) $ and $\Psi \in H_{\infty}.$
If $\Psi|S_T =0,$ then
$\|\Psi |S_{T, \infty }\|_{H_{T,  \infty }} =\|\Psi \|_{H_{\infty }}$.
If $\widetilde{ T} \in [0,T],$ then
$\| \Psi | S_{T, \infty }\|_{H_{ T, \infty }}
\le
\| \Psi | S_{\widetilde{ T}, \infty }\|_{H_{\widetilde{ T}, \infty }}.$
\end{corollary}
{\bf Proof:} Lemma \ref{lemmaL10.30}.
\begin{lemma} \label{lemmaL10.20}
Let $T \in (0, \infty ) ,\;\Psi \in H_{\infty },$ and suppose there is
$v \in C^1 \bigl(\, \mathbb{R}^3 \times (T, \infty ) \,\bigr) ^3$
with $\Psi |S_{T, \infty }=v|S_{T, \infty }$.

Then $\int_{ 0}^t(t-r) ^{-1/2} \, |\Psi (x,r)|\, dr < \infty $
for $t \in (T , \infty ) $ and for a. e. $x \in \partial \Omega $.
Put
$W(x,t):= \int_{ 0}^t(t-r) ^{-1/2} \, \Psi (x,r)\, dr$
for these $t$ and $x$. Then $W(x,\, \cdot \,) \in C^1 \bigl(\, (T,  \infty ) \,\bigr) ^3$
and equation (\ref{L10.10.10}) holds with $\phi  $ replaced by $\Psi (x,\, \cdot \, )$,
for a. e. $x \in \partial \Omega $.
Therefore the fractional derivative $\partial _t ^{1/2} \Psi (x,t)$ is well defined and equals
$\partial _4W(x,t)$ for $t \in (T, \infty )$ and a. e. $x \in \partial \Omega $. Moreover
\begin{eqnarray} \label{L10.20.05}
\|\Psi |S_{T, \infty  }\|_{H_{T, \infty }}^2= \int^{ \infty }_T \Bigl( \|\Psi (t)\|_{H^1( \partial \Omega )^3}^2
+
\| \partial _4W (t)\|_{2}^2
+
\| n^{ (\Omega )} \cdot  \partial _4\Psi (t)\|_{H^1( \partial \Omega ) ^{\prime} }^2 \Bigr) \, dt
.
\end{eqnarray}
\end{lemma}
{\bf Proof:} Since $\Psi \in H_{\infty }$, we have in particular $\Psi \in L^2(S_{\infty })^3,$
so $\Psi \in L^1_{loc} \bigl(\, [0, \infty ) \, L^2( \partial \Omega )^3 \,\bigr) .$
In addition $\Psi |S_{T, \infty }=v|S_{T, \infty }$ with
$v \in C^1 \bigl(\, \mathbb{R}^3 \times (T, \infty ) \,\bigr) ^3$,
so the claims about $W$ are valid according to Corollary \ref{corollaryC10.10}. This leaves us to show
(\ref{L10.20.05}). Again since $\Psi \in H_{\infty }$, there are functions
$v_n \in C ^{ \infty }  ( \mathbb{R}^4 ) ^3$ with $v_n| \mathbb{R}^3 \times (-\infty ,0)=0$
for $n \in \mathbb{N}  $ such that $\|v_n-\Psi\|_2 \to 0$ and such that $(v_n|S_{\infty })$
is a Cauchy sequence with respect to the norm $\|\;\|_{H_{\infty }}$ on $\widetilde{ H}_{T_0}$.
The latter property implies there are mappings
$\Psi ^{(1)} \in L^2 \bigl(\, 0, \infty  ,\, H^1( \partial \Omega )^3 \,\bigr) ,\; \Psi ^{(2)} \in L^2(S_{\infty })^3$
and
$\Psi ^{(3)} \in L^2 \bigl(\, 0, \infty  ,\, H^1( \partial \Omega ) ^{\prime}  \,\bigr)$
as in (\ref{L10.30.10}) with $T_0=\infty $.
Thus
(\ref{L10.30.20}) holds, so by the definition of $\|\;\|_{H_{T,\infty} }$,
\begin{eqnarray} \label{L10.20.20}
\|\Psi |S_{ T , \infty }\|_{H_{ T, \infty }}^2
=
\int^{ \infty }_{ T} \Bigl( \|\Psi ^{(1)}  (t)\|_{H^1( \partial \Omega )^3}^2
+
\| \Psi ^{(2)}  (t)\|_{2}^2
+
\| \Psi ^{(3)}  (t)\|_{H^1( \partial \Omega ) ^{\prime} }^2 \Bigr) \, dt 
\end{eqnarray}
We further conclude from (\ref{L10.30.10}) and the choice of the sequence $(v_n)$ that
$\Psi ^{(1)} =\Psi $. Let us show that $\Psi ^{(2)} |S_{T, \infty }=\partial _4W.$
To this end, we proceed as in the proof of Lemma \ref{lemmaL10.30}. We take
$V \in L^2( \partial \Omega )^3,\; \phi \in C ^{ \infty } _0 \bigl(\, (T, \infty ) \,\bigr) $,
and we put
$w_n(x,t):=\int_{ 0}^t (t-r) ^{-1/2} \, v_n(x,r)\, dr$
for $n \in \mathbb{N} ,\; x \in \partial \Omega ,\; t \in (0, \infty ) $.
By Lemma \ref{lemmaL10.10}, we have $w_n(x,\, \cdot \,) \in C^1 \bigl(\, (0, \infty ) \,\bigr) ^3$
for any $x \in \partial \Omega $, and $\partial ^{1/2} _4v_n =\partial _4w_n$ for $n \in \mathbb{N} $.
Similarly to (\ref{L10.30.30}) and (\ref{L10.30.35}), we get
\begin{eqnarray} \label{L10.20.40}
\int_{ \partial \Omega }\int_{ T} ^{\infty}  \Psi ^{(2)} (x,t) \cdot V(x) \, \phi (t)\, dt\, do_x
=
-\lim_{n\to \infty }
\int_{ \partial \Omega }\int_{ T} ^{ \infty }  w_n (x,t) \cdot V(x) \, \phi ^{\prime} (t)\, dt\, do_x.
\end{eqnarray}
There is $T_1 \in (T, \infty )$ such that $supp( \phi ) \subset (T,T_1).$
Thus, with $W$ introduced in the lemma,
\begin{eqnarray*} &&
\Bigl| \int_{ \partial \Omega }\int_{ T} ^{ \infty }  (w_n-W) (x,t) \cdot V(x) \, \phi ^{\prime} (t)\, dt\, do_x
\Bigr|
\\ &&
=
\Bigl| \int_{ \partial \Omega }\int_{ 0} ^{ T_1 } \int_{ 0}^t (t-r) ^{-1/2} \, (v_n-\Psi) (x,r)
\cdot V(x) \, \phi ^{\prime} (t)\, dr \, dt\, do_x \Bigr|
\\&&
=
\Bigl| \int_{ \partial \Omega }\int_{ 0} ^{ T_1 } \int_{ r}^{T_1} (t-r) ^{-1/2} \,\phi ^{\prime} (t)\, dt
\, (v_n-\Psi) (x,r) \cdot  V(x) \,  dr \,  do_x \Bigr| 
\\&&
\le
\mbox{$\mathfrak C$} (T_1)\, | \phi ^{\prime} |_{\infty }\, 
\int_{ \partial \Omega } \|(v_n-\Psi )(x,\,  \cdot \,) \|_2\, | V(x)|\; do_x
\le
\mbox{$\mathfrak C$} (T_1)\, | \phi ^{\prime} |_{\infty } \,\|V\|_2\, \|v_n-\Psi \|_2
\end{eqnarray*}
Since $ \|v_n-\Psi \|_2\to 0$, we thus get
\begin{eqnarray} \label{L10.20.50}
-\int_{ \partial \Omega }\int_{ T}^{ \infty } W(x,t) \cdot  V(x)\, \phi ^{\prime} (t)\, dt\, do_x
=
-\lim\limits_{n\to \infty }\int_{ \partial \Omega }\int_{ T}^{ \infty } w_n(x,t) \cdot  V(x)\, \phi ^{\prime} (t)
\, dt\, do_x.
\end{eqnarray}
On the other hand, since $W(x,\, \cdot \, ) \in C^1 \bigl(\, (t, \infty ) \,\bigr) $
for a. e. $x \in \partial \Omega $,
\begin{eqnarray} \label{L10.20.60}
-\int_{ \partial \Omega }\int_{ T}^{ \infty } W(x,t) \cdot  V(x)\, \phi ^{\prime} (t)\, dt\, do_x
=
\int_{ \partial \Omega }\int_{ T}^{ \infty } \partial _4W(x,t) \cdot  V(x)\, \phi (t)\, dt\, do_x.
\end{eqnarray}
Combining (\ref{L10.20.40}) -- (\ref{L10.20.60}), we may conclude that
$\Psi ^{(2)} |S_{T, \infty }=\partial _4W$. An analogous but simpler reasoning yields that
$\Psi ^{(3)}(x,t)=n^{ (\Omega )} (x) \cdot  \partial _4\Psi (x,t)$
for $t \in (T, \infty )$ and for a. e. $x \in \partial \Omega $. Now equation (\ref{L10.20.05}) follows
from (\ref{L10.20.20}).
\hfill $\Box $

\section{Some fundamental solutions and potential functions. A representation formula
for the velocity part of solutions to the Oseen system with Dirichlet boundary conditions.}
\setcounter{equation}{0}

We define three potential functions, denoted by $\mbox{$\mathfrak R$} ^{( \tau )} (f),\, \mbox{$\mathfrak I$}
^{( \tau )} (a)$ and $\mbox{$\mathfrak V$} ^{( \tau )} ( \phi ), $ respectively. The first is related to the
right-hand side $f$ in (\ref{10}), the second to the initial data $a$ in (\ref{30}), and the third to
the Dirichlet boundary data $b$ in (\ref{20b}). We begin by introducing fundamental solutions to
the heat equation, the time-dependent Stokes system and
the Oseen system (\ref{10}), respectively.
We write \mbox{$\mathfrak H$} for the fundamental solution
of the heat equation in $\mathbb{R}^3 $, that is,
\begin{eqnarray*}
\mbox{$\mathfrak H$} (z,t):=(4\,\pi\, t) ^{-3/2} \,e^{-|z|^2/(4\, t)}
\quad \mbox{for}\;\; (z,t) \in \mathbb{R}^3 \times (0, \infty ).
\end{eqnarray*}
As a fundamental solution of the time-dependent Stokes system, we choose the same function
$\Gamma $ as in \cite{Shen}, that is,
\begin{eqnarray*}
\Gamma _{jk} (z,t):= \delta _{jk} \, \mbox{$\mathfrak H$} (z,t)+ \int_{ t} ^{ \infty } \partial _j \partial _k
\mbox{$\mathfrak H$} (z,s)\, ds
\quad \mbox{for}\;\; (z,t) \in \mathbb{R}^3 \times (0, \infty ),\; 1\le j,k\le 3.
\end{eqnarray*}
Actually this is the velocity part of the fundamental solution in question; we will not need
the  pressure part associated with. Finally we introduce the velocity part of the looked-for fundamental
solution of the time-dependent Oseen system (\ref{10}), setting
\begin{eqnarray*} 
\Lambda  _{jk} (z,t, \tau ):= \Gamma _{jk} (z-t\, \tau \,e_1,\, t)
\quad \mbox{for}\;\; z,\, t,\, j,\, k\; \mbox{as before}.
\end{eqnarray*}
Note that in what follows, $\Gamma $ does not stand for the usual Gamma function.
We state some properties of $\mbox{$\mathfrak H$} ,\; \Gamma =(\Gamma _{jk} )_{1\le j,k\le 3}$
and $\Lambda = ( \Lambda _{jk} )_{1\le j,k\le 3}$.
\begin{lemma} \label{lemmaL10.40}
$\mbox{$\mathfrak H$} \in C ^{ \infty } \bigl(\, \mathbb{R}^3 \times (0, \infty ) \,\bigr) ,
\; \int_{ \mathbb{R}^3 }\mbox{$\mathfrak H$} (z,t)\, dz = 1$
for
$ t \in (0, \infty ) ,$
and
$| \partial _t^l \partial _z ^{\alpha }\mbox{$\mathfrak H$} (z,t)|\le
C \, (|z|^2+t)^{-3/2-| \alpha |/2-l}$
for
$z \in \mathbb{R}^3 ,\; t \in (0, \infty ) ,\; \alpha \in \mathbb{N} _0^3,\; l \in \mathbb{N} _0$
with $| \alpha |+l\le 1.$
\end{lemma}
{\bf Proof:}
For a proof of the estimate at the end of the lemma, we refer to
\cite{Sol1}. An estimate of this kind holds for any $\alpha \in \mathbb{N} _0^3,\; l \in \mathbb{N} _0,$
but we will need it only in the case $| \alpha |+l\le 1$.
\hfill $\Box $
\begin{lemma}[\mbox{\cite[Proposition 2.19]{Shen}}] \label{lemmaL10.50}
$\Gamma _{jk} ,\: \Lambda _{jk}(\, \cdot \, ,\, \cdot \, , \tau )
\in C ^{ \infty } \bigl(\, \mathbb{R}^3 \times (0, \infty ) \,\bigr) ,$
\begin{eqnarray*} &&
| \partial _t^l \partial _z ^{\alpha }\Gamma _{jk}  (z,t)|\le
C \, (|z|^2+t)^{-3/2-| \alpha |/2-l},
\\&&
| \partial _t^l \partial _z ^{\alpha }\Lambda  _{jk}  (z,t)|
\le
C( \tau )\, \bigl[\, (|z- t\, \tau \, e _1 |^2+t)^{-3/2-| \alpha |/2-l}
+ \delta _{l1}\, (|z- t\, \tau \, e_1|^2+t) ^{-2}\,\bigr] 
\end{eqnarray*}
for
$z \in \mathbb{R}^3 ,\; t \in (0, \infty ) ,\; \alpha \in \mathbb{N} _0^3,\; l \in \mathbb{N} _0$
with $| \alpha |+l\le 1,\; j,k \in \{1,\, 2,\, 3\} .$
\end{lemma}
By combining Lemma \ref{lemmaL10.3} and the preceding lemma, we get
\begin{corollary} \label{corollaryC10.60}
Let $K \in (0, \infty ) $ and put $\gamma _K(z):=|z|^2$ for $z \in  B_K, \;
\gamma _K(z):=|z|\, \nu (z)$ for $z \in B_K^c$.
Then
$
| \partial _t^l \partial _z ^{\alpha }\Lambda _{jk}    (z,t)|
\le
C( \tau ,K )\, \bigl[\, ( \gamma _K(z)+t)^{-3/2-| \alpha |/2-l} + \delta _{l1}\, ( \gamma _K(z)+t) ^{-2} \,\bigr]
$
for $z,\, t,\, \alpha ,\, l,\, j,\,k $ as in the preceding lemma.
\end{corollary}
The ensuing theorem provides an estimate of convolutions of $\Lambda $.
\begin{theorem} \label{theoremT4.30}
Let $ q \in [1, \infty ),\; \varrho  \in (1, \infty ],\;
s \in [1, \infty ]$ with $  s\le \varrho ,\;
M \in (0, \infty ) ,\; j,k \in \{1,\, 2,\, 3\} ,\; \alpha  \in \mathbb{N} _0^3$
with $|\alpha |  \le 1$. Take $h \in L^s\bigl(\, 0, \infty ,\, L^q( \mathbb{R}^3 ) \,\bigr) $.

Then the ensuing inequality holds
for $W= (0,M)$ if $\; 1-|\alpha | /2   - 3/(2\, q) - 1/s +1/ \varrho >0$,
and for $W= (M, \infty )$ if $\; 1-|\alpha | /2  -3/(2\, q) -1/s+1/ \varrho <0$:
\begin{eqnarray*} &&
\Bigl( \int_{ 0}^{ \infty } 
\Bigl( \int_{ 0}^{ \infty } \int_{ \mathbb{R}^3 }
\mbox{\Large $\chi$}_W
(t- \sigma ) \, 
| \partial _x^{\alpha}  \Lambda  _{jk}(x-y,t- \sigma , \tau )| 
\,
|h(y, \sigma )| \; dy \; d \sigma \Bigr) ^{\varrho } \; dt \Bigr) ^{1/ \varrho }
\\&& \nonumber
\le 
C( \tau ,q, \varrho ,s) \, M^{ 1-| \alpha | /2 -3/(2\, q) -1/s +1 / \varrho } \, \|h\|_{q,s; \infty }
\end{eqnarray*}
for $x \in \mathbb{R}^3 $
if $\varrho < \infty $, and
\begin{eqnarray*} &&
\int_{ 0}^{ \infty } \int_{ \mathbb{R}^3 }
\mbox{\Large $\chi$}_{W}
(t- \sigma ) \, 
| \partial _x^{\alpha}  \Lambda  _{jk}(x-y,t- \sigma , \tau )| 
\,
|h(y, \sigma )| \; dy \; d \sigma 
\\ &&
\le
C( \tau ,q,s) \, M^{ 1-| \alpha | /2 -3/(2\, q)-1/s } \, \|h\|_{q, s ; \infty }
\end{eqnarray*}
for $x \in \mathbb{R}^3 ,\; t \in (0, \infty ) $ if $\varrho =\infty $.
\end{theorem}
{\bf Proof:} Theorem \ref{theoremT4.30} follows from \cite[Lemma 2.7]{DeSIMA1} and
\cite[Theorem 2.8]{DeSIMA2} if the parameter $p$ in those references is chosen as $p= \infty $.
As becomes apparent from the proof of these references, not only $L ^{ \infty } $-estimates are provided
if $p= \infty $ or $\varrho = \infty $, but pointwise estimates not involving any sets of measure zero.
\hfill $\Box $

\vspace{1ex}
The ensuing lemma is the basis of the definition of our first potential function.
\begin{lemma} \label{lemmaL10.60}
Let $q,\, s \in [1, \infty )$ and $h \in L^s \bigl(\, 0, \infty ,\, L^q( \mathbb{R}^3 )^3 \,\bigr) .$
Then, for a. e. $(x,t) \in \mathbb{R}^3 \times (0, \infty ),\; \alpha \in \mathbb{N} _0^3$ with
$| \alpha |\le 1,$ we have
$
\int_{ 0}^{ t } \int_{ \mathbb{R}^3 }
| \partial _x^{\alpha}  \Lambda (x-y,t- \sigma , \tau )  
\cdot 
h(y, \sigma )| \, dy \, d \sigma 
< \infty .
$
\end{lemma}
{\bf Proof:}
The lemma follows from a more general version of Theorem \ref{theoremT4.30} (\cite[Lemma 2.7]{DeSIMA1});
see the remarks in \cite[p. 898 below]{DeSIMA1}. \hfill $\Box $

\vspace{1ex}
Due to the preceding lemma, we may define a function $\mbox{$\mathfrak R$} ^{( \tau )}(h) : \mathbb{R}^3 \times
[0, \infty ) \mapsto \mathbb{R}^3 $ for any $h \in L^s \bigl(\, 0, T,\, L^q(A)^3 \,\bigr) $
with $q,s \in [1, \infty ),\; A \subset \mathbb{R}^3 $ measurable, $T \in (0, \infty ],$
by setting
\begin{eqnarray*}
\mbox{$\mathfrak R$} ^{( \tau )}(x,t):= \int_{ 0}^t \int_{ \mathbb{R}^3 } \Lambda (x-y,t- \sigma , \tau ) \cdot 
\widetilde{ h}(y, \sigma )\, dy\, d \sigma
\quad \mbox{for a. e.}\; (x,t) \in \mathbb{R}^3 \times [0, \infty ),
\end{eqnarray*}
where $\widetilde{ h}$ denotes the zero extension of $h$ to $\mathbb{R}^3 \times (0, \infty ).$
\begin{lemma}[\mbox{\cite[Lemma 2.11]{DeSIMA1}}] \label{lemmaL10.70}
Let $q,\, s \in [1, \infty )$ and $h \in L^s \bigl(\, 0, \infty ,\, L^q( \mathbb{R}^3 )^3 \,\bigr) .$
Then the weak derivativee $\partial _l \mbox{$\mathfrak R$} ^{( \tau )} (h)$ exists for
$1\le l\le 3$, in particular $\mbox{$\mathfrak R$} ^{( \tau )} (h)(t) \in W^{1,1}_{loc}( \mathbb{R}^3 )^{3 \times 3}$
for a. e. $t \in (0, \infty ).$
Moreover
$
\partial _l\mbox{$\mathfrak R$} ^{( \tau )}(x,t)= \int_{ 0}^t \int_{ \mathbb{R}^3 }
\partial _l\Lambda (x-y,t- \sigma , \tau ) \cdot 
h(y, \sigma )\, dy\, d \sigma
$
for a. e. $x \in \mathbb{R}^3 ,$ a. e. $t \in  (0, \infty ),\; 1\le l\le 3.$
In particular the trace of $\mbox{$\mathfrak R$} ^{( \tau )} (h)(t)$ on $\partial \Omega $
is well defined for a. e. $t \in (0, \infty ) $.
\end{lemma}
\begin{lemma}[\mbox{\cite[Lemma 2.12, Corollary 2.13]{DeSIMA1}}] \label{lemmaL10.80}
Let $q,\, s \in [1, \infty )$ and let $h$ be a function belonging to
$L^s \bigl(\, 0, \infty ,\, L^q( \mathbb{R}^3 )^3 \,\bigr) .$
Then 
$
\int_{ 0}^t \int_{ \mathbb{R}^3 }|\Lambda (x-y,t- \sigma , \tau ) \cdot 
h(y, \sigma )|\, dy\, d \sigma < \infty 
$
for a. e. $(x,t) \in \partial \Omega \times (0, \infty ).$
Thus $\mbox{$\mathfrak R$} ^{( \tau )} (h)$ is well defined also as a function on $S_{\infty }$.
Moreover, $\mbox{$\mathfrak R$} ^{( \tau )} (h)(t)$ as a function on $\partial \Omega $ is the
trace of $\mbox{$\mathfrak R$} ^{( \tau )} (h)(t)$ as a function on $\mathbb{R}^3 $, for a. e.
$t \in (0, \infty ) $.
\end{lemma}
\begin{lemma} \label{lemmaL10.85}
Let $T \in (0, \infty ) ,\; q \in [1,4], \; s \in [1, \infty ),\;
f \in L^s \bigl(\, 0,T,\;L^q( \mathbb{R}^3 )^3 \,\bigr) .$
Then $\mbox{$\mathfrak R$} ^{( \tau )} (f)| \mathbb{R}^3 \times (T, \infty ) \in C^1 \bigl(\,
\mathbb{R}^3 \times (T, \infty ) \,\bigr) ^3$ and
\begin{eqnarray} \label{L10.85.10}
\partial _t^l \partial ^{\alpha }_x \mbox{$\mathfrak R$} ^{( \tau )} (f)(x,t)
=
\int_{ 0}^T \int_{ \mathbb{R}^3 } \partial _t^l \partial _x^{\alpha } \Lambda (x-y,t-\sigma , \tau ) \cdot
f(y, \sigma )\, dy\, d \sigma
\end{eqnarray}
for $x \in \mathbb{R}^3 ,\; t \in (T, \infty ) ,\; l \in \mathbb{N} _0,\; \alpha \in \mathbb{N} _0^3$
with $l+| \alpha |\le 1,$ where the preceding integral exists as a Lebesgue integral.
\end{lemma}
{\bf Proof:}
The lemma follows from Lebesgue's theorem. However, since $f $ is not required to belong
to $L^1 \bigl(\, \mathbb{R}^3 \times (0,T) \,\bigr) ^3$, it is perhaps not completely obvious
how to apply that theorem. In particular, the reason for the condition $q\le 4$
may not be clear. So we indicate a proof.
Let $R \in (0, \infty ) ,\; T_0 \in (T, \infty ).$ It is enough to show that
$\mbox{$\mathfrak R$} ^{( \tau )} (f)|B_R \times (T_0, \infty )$ is a $C^1$-function and
equation (\ref{L10.85.10}) holds for $(x,t) \in B_R \times (T_0,  \infty ).$ Put
$f ^{(1)} := \chi_{B_{2 R}^c \times (0,T)}\, f.$
Then, by Corollary \ref{corollaryC10.60} with $K=R$ and Lemma \ref{lemmaL10.1}, we get
\begin{eqnarray} \label{L10.85.20}
|\partial _t^l \partial _x^{\alpha } \Lambda (x-y,t-\sigma , \tau ) \cdot f ^{(1)} (y, \sigma )|
\le
C(R, \tau )\, g_{l, \alpha }(y)\,|f(y, \sigma )|
\end{eqnarray}
for $x \in B_R,\; y \in \mathbb{R}^3 ,\; t \in (T_0, \infty ),\; \sigma \in (0,T),\; \alpha ,\, l$
as in (\ref{L10.85.10}), where
$$
g_{l, \alpha }(y):= \chi_{(2R,\, \infty )}(|y|)\, \bigl[\,  \bigl(\, |y|\, \nu (y) \,\bigr) ^{-3/2-| \alpha |/2-l}
+ \delta _{l1}\, \bigl(\, |y|\,\nu (y) \,\bigr) ^{-2} \,\bigr]
\quad \mbox{for}\;\; y \in \mathbb{R}^3 .
$$
Since $q<4$, we have $3\, q ^{\prime} /2>2,$ so by Corollary \ref{corollaryC1.10},
\begin{eqnarray*} &&
\int_{ 0}^T \int_{ \mathbb{R}^3 } g_{l, \alpha }(y)\, |f(y, \sigma )|\, dy\, d \sigma
\le
\int_{ 0}^T \Bigl( \int_{ \mathbb{R}^3 } g_{l, \alpha }(y)^{q ^{\prime} }\, dy \Bigr) ^{1/p ^{\prime} }\,
\|f( \sigma )\|_q\, d \sigma
\\&&
\le
C(R,q,s, \tau )\, T^{1/s ^{\prime} }\,\|f\|_{q,s;T}
<
\infty ,
\end{eqnarray*}
for $l,\, \alpha $ as in (\ref{L10.85.10}). Therefore from (\ref{L10.85.20}) and Lesbesgue's theorem,
we may conclude that $\mbox{$\mathfrak R$} ^{( \tau )} (f ^{(1)} )$ is a $C^1$-function on
$B_R \times (T_0, \infty ),$
and equation (\ref{L10.85.10}) is valid for $x \in B_R,\; t \in (T_0, \infty )$ and
with $f$ replaced by $f ^{(1)} $. Put $f ^{(2)} :=\chi_{B_{2 R}\times (0,T)}\, f$.
Then $f ^{(2)} \in L^1 \bigl(\, \mathbb{R}^3 \times (0,T) \,\bigr) ^3,$
and due to Lemma \ref{lemmaL10.50}, we have
\begin{eqnarray*}
|\partial _t^l \partial _x^{\alpha } \Lambda _{jk}  (x-y,t-\sigma , \tau ) |
\le
C(\tau )\, \bigl(\, (T_0-T)^{-3/2-| \alpha |/2-l}+ \delta _{1l}\, (T_0-T) ^{-2} \,\bigr) 
\end{eqnarray*}
for $x , \, y,\, t,\, \sigma $ as in (\ref{L10.85.20}), $1\le j,k\le 3$ and  $\alpha ,\, l$
as in (\ref{L10.85.10}). Therefore Lebesgue's theorem yields that
$\mbox{$\mathfrak R$} ^{( \tau )} ( f ^{(2)} )|B_R \times (T_0, \infty )$
is a $C^1$-function and (\ref{L10.85.10}) holds on $B_R \times (T_0, \infty )$
if $f$ is replaced by $f ^{(2)} $. Since $f=f ^{(1)} +f ^{(2)} $, the lemma is proved.
\hfill $\Box $

\vspace{1ex}
The next theorem presents a criterion on $f$ implying $\mbox{$\mathfrak R$} ^{( \tau )} (f)|S_{\infty }
\in H_{\infty }.$
\begin{theorem}[\mbox{\cite[Theorem 8.1]{DeJMFM}}] \label{theoremT10.5}
Let
$f \in L^2 \bigl(\, 0, \infty ,\;L^{\varrho }( \mathbb{R}^3 )^3 \,\bigr)$
for $\varrho =3/2$ and for some $\varrho  \in [1,\, 3/2).$
Then $\mbox{$\mathfrak R$} ^{( \tau )} (f)|S_{\infty } \in H_{\infty }$ and
$
\| \mbox{$\mathfrak R$} ^{( \tau )} (f)|S_{\infty }\|_{H_{\infty }}
\le
\mbox{$\mathfrak C$} \ (\|f\|_{3/2,2; \infty }+\|f\|_{q,2; \infty }).
$
\end{theorem}
\begin{corollary} \label{corollaryC10.65}
Let $T,\, q,\, f$ be given as in Lemma \ref{lemmaL10.85}, and abbreviate
$\Psi := \mbox{$\mathfrak R$} ^{( \tau )} (f)|S_{\infty }$.
For a. e. $x \in \partial \Omega $ we then have
$\Psi (x,\, \cdot \,)|(T, \infty ) \in C^1 \bigl(\, (T, \infty ) \,\bigr)^3 ,\;
\int_{ 0}^t (t-r) ^{-1/2} \, |\Psi(x,r)|\, dr < \infty
$
for $t \in (T, \infty ),$
and the function
$t \mapsto \int_{ 0}^t (t-r) ^{-1/2} \, \Psi (x,r)\, dr,\; t \in (T, \infty ),$
belongs to $C^1 \bigl(\, (T, \infty ) \,\bigr) ^3$, so that
$\partial ^{1/2} _t\Psi (x,t)$ exists for any $t \in (T, \infty )$. In addition,
equation (\ref{L10.10.10}) holds with $\Psi (x,\, \cdot \, )$ in the place of $\phi $,
for a. e. $x \in \partial \Omega .$
Finally,
\begin{eqnarray}&& \label{C10.65.10}
\|\Psi |S_{T, \infty }\|_{H_{T, \infty }}^2
\\&&\nonumber 
=
\int_{ T} ^{ \infty } \Bigl( \|\Psi (t)| \partial \Omega \|_{H^1( \partial \Omega )^3}^2
+
\int_{ \partial \Omega }| \partial _t ^{1/2} \Psi(x,t)|^2\, do_x
+
\int_{ \partial \Omega }| n^{ (\Omega )} (x) \cdot \partial _t\Psi (x,t)|^2\, do_x \Bigr) \, dt
\end{eqnarray}
\end{corollary}
{\bf Proof:} Theorem \ref{theoremT10.5}, Lemma \ref{lemmaL10.85} and \ref{lemmaL10.20}. \hfill $\Box $

\vspace{1ex}
We will need the following estimate on pointwise spatial decay of $\mbox{$\mathfrak R$} ^{( \tau )} (f)$.
\begin{theorem}[\mbox{\cite[Theorem 3.1]{DeDCDS-A}}]
\label{theoremT4.10}
Let $A \in (2, \infty ),\; B \in [0,\, 3/2]$ with $A+\min\{1,B\} > 3,\; A+B\ge 7/2,\;
\varrho _0 \in (2, \infty ),\; R_0 \in (0, \infty ) ,\; \gamma \in L^2 \bigl(\, (0, \infty ) \,\bigr)
\cap L^{\varrho _0} \bigl(\, (0, \infty ) \,\bigr) ,\; f: \mathbb{R}^3 \times (0,  \infty ) \mapsto \mathbb{R}^3 $
measurable with
$|f(y, \sigma )|\le \gamma ( \sigma ) \, |y| ^{-A}\, \nu (y) ^{-B}$
for $y \in B_{R_0}^c,\; \sigma \in (0, \infty ) .$
Further suppose that $f|B_{R_0} \times (0, \infty ) \in L^2 \bigl(\, B_{R_0} \times (0, \infty ) \,\bigr) ^3.$
Let $R \in (R_0, \infty )$.

Then, for $x \in B_R^c,\; t \in (0, \infty ) ,\; \alpha \in \mathbb{N} _0^3$
with $ | \alpha |\le 1,$
inequality (\ref{130}) holds.

(Note that by Lemma \ref{lemmaL10.60} and \ref{lemmaL10.70} with $q=s=2$, and because
$\gamma \in L^2 \bigl(\, (0, \infty ) \,\bigr) $ and $A >2$, the term
$\partial _x^{\alpha } \mbox{$\mathfrak R$} ^{( \tau )} (f)(x,t)$ in (\ref{130}) is well defined.)
\end{theorem}
The next lemma allows to define the potential $\mbox{$\mathfrak I$} ^{( \tau )} (a)$ further below and
yields a pointwise temporal estimate of this potential.
\begin{lemma} \label{lemmaL10.90}
Let $p \in [1, \infty ], \; a \in L^p( \mathbb{R}^3 )^3,\; l \in \mathbb{N} _0,\; \alpha \in \mathbb{N} _0^3$
with $l+| \alpha |\le 1.$ Then
$$
\int_{ \mathbb{R}^3 }| \partial _t^l \partial _x^{\alpha }
\bigl(\, \mbox{$\mathfrak H$}  (x-y -\tau \, t\, e_1,\, t ) \,\bigr) \, a(y)|\, dy
\le
C(p)\,  \|a\|_p\,
(t^{-3/(2\, p)-| \alpha |/2-l}+ \delta _{1l}\, t^{-3/(2 \, p)-1/2})
$$
for $x \in \mathbb{R}^3 ,\; t \in (0, \infty ) ,$ where
$3/(2\, p):=0 $ if $p=\infty $.
\end{lemma}
{\bf Proof:}
Take $x,\, t$ as in the lemma, and let the left-hand side of the estimate
in Lemma \ref{lemmaL10.90} be denoted by \mbox{$\mathfrak A$}.
In the case $p \in (1, \infty )$, we get with Lemma \ref{lemmaL10.40} that
\begin{eqnarray*} &&
\mbox{$\mathfrak A$}
=
\int_{ \mathbb{R}^3 } | \partial ^l_4 \partial ^{\alpha }_x \mbox{$\mathfrak H$} (x-y -\tau \, t\, e_1,\, t)
-\delta _{1l}\, \tau \, \partial _1 \mbox{$\mathfrak H$} (x-y -\tau \, t\, e_1,\, t)|\, |a(y)|\, dy
\\&& 
\le
C\, \|a\|_p\, \Bigl( \int_{ \mathbb{R}^3 }\Bigl[ (|x-y-\tau \, t\, e_1|+t ^{1/2} )^{(-3-| \alpha |-2\, l)\,p ^{\prime} }
\\&& \hspace{15em}
+ \delta _{1l}\, (|x-y-\tau \, t\, e_1|+t ^{1/2} )^{-4\,p ^{\prime} }\Bigr]\, dy \Bigr) ^{1/p ^{\prime} }
\\&&
\le
C\, \|a\|_p\, \Bigl( \int_{ \mathbb{R}^3 }\Bigl[ (|z|+t ^{1/2} )^{(-3-| \alpha |-2\, l)\,p ^{\prime} }
+ \delta _{1l}\, (|z|+t ^{1/2} )^{-4\,p ^{\prime} }\Bigr]\, dz \Bigr) ^{1/p ^{\prime} }
\\&&
\le
C\, \|a\|_p\, ( t^{-3/2-| \alpha |/2- l +3/(2p ^{\prime} ) }
+ \delta _{1l}\, t ^{-2+3/(2p ^{\prime}) })
\\&&
\le
C\, \|a\|_p\, ( t^{-3/(2p)-| \alpha |/2- l }
+ \delta _{1l}\, t^{-3/(2p)-1/2 })
.
\end{eqnarray*}
If $p=1$, the estimate
$
\mbox{$\mathfrak A$} \le C\, \|a\|_p\, ( t^{-3/(2p)-| \alpha |/2- l }
+ \delta _{1l}\, t^{-3/(2p)-1/2 })
$
follows almost immediately from Lemma \ref{lemmaL10.40}.
Suppose $p=\infty $. Then, in the case $\alpha =0,\; l=0$,
\begin{eqnarray*}
\mbox{$\mathfrak A$}
\le
\|a\|_p\, \int_{ \mathbb{R}^3 }\mbox{$\mathfrak H$} (x-y-\tau \, t\, e_1,\, t)\, dy
=
\|a\|_p\, \int_{ \mathbb{R}^3 }\mbox{$\mathfrak H$} (z,t)\, dz =\|a\|_p
\end{eqnarray*}
by Lemma \ref{lemmaL10.40}, and in the case $| \alpha |+l>0$ by the same reference,
\begin{eqnarray*}&&
\mbox{$\mathfrak A$}
\le
C\, \|a\|_p\, \int_{ \mathbb{R}^3 }\bigl(\, (|z|+ t ^{1/2} )^{-3-| \alpha |-2l}+\delta _{l1}\,(|z|+t ^{1/2} ) ^{-4} \,\bigr)
\, dz
\\&&
\le
C\, \|a\|_p\,(t^{-| \alpha |/2-l}+\delta _{1l}\, t ^{-1/2} ).
\end{eqnarray*}
This completes the proof of the lemma.
\hfill $\Box $

\vspace{1ex}
In view of Lemma \ref{lemmaL10.90}, we may define
\begin{eqnarray*}
\mbox{$\mathfrak I$} ^{( \tau )} (a)(x,t):= \int_{ \mathbb{R}^3 } \mbox{$\mathfrak H$} (x-y-\tau \,t\, e_1,\, t)\,
\widetilde{ a}(y)\, dy
\quad 
\bigl(\, x \in \mathbb{R}^3 ,\; t \in (0, \infty ) \,\bigr) ,
\end{eqnarray*}
where
$p \in [1, \infty ],\; a \in L^p(A)^3,$ for some
measurable subset $A$ of $ \mathbb{R}^3 $, with $\widetilde{ a}$ denoting the zero
extension of $a$ to $\mathbb{R}^3 $.
\begin{lemma}
\label{lemmaL10.100}
Let $p \in [1, \infty ]$ and $ a \in L^p( \mathbb{R}^3 )^3.$
Then $\mbox{$\mathfrak I$} ^{( \tau )} (a) \in C^1 \bigl(\, \mathbb{R}^3 \times (0, \infty ) \,\bigr) ^3$
and
\begin{eqnarray} \label{L10.100.10}
\partial _t^l \partial _x^{\alpha }\mbox{$\mathfrak I$} ^{( \tau )} (a)(x,t)
=
\int_{ \mathbb{R}^3 } \partial _t^l \partial _x^{\alpha } \bigl(\, \mbox{$\mathfrak H$}  (x-y- \tau \, t\, e_1,\, t)
\,\bigr) \, a(y)\, dy
\end{eqnarray}
for
$
x \in \mathbb{R}^3 ,\; t \in (0, \infty ),\; \alpha \in \mathbb{N} _0^3,\; l \in \mathbb{N} _0\;
\mbox{with}\; | \alpha |+l\le 1.
$
\end{lemma}
{\bf Proof:} See \cite[Lemma 2.3]{De7} and its proof.

\vspace{1ex}
Actually $\mbox{$\mathfrak I$} ^{( \tau )} (a)$ is a $C ^{ \infty } $-function, but we will not need
this fact.
\begin{theorem}[\mbox{\cite[Theorem 3.1]{DeSIMA1}}] \label{theoremT10.6}
Let $\epsilon _0 \in (0,1/2]$ and $a \in H^{1/2+\epsilon _0}_{\sigma }( \overline{ \Omega }^c).$
Then the function $\mbox{$\mathfrak I$} ^{( \tau )} (a)|S_{\infty }$ belongs to $ H_{\infty }.$
\end{theorem}
\begin{corollary} \label{corollaryC10.66}
Let $\epsilon _0$ and $a$ be given as in the preceding theorem. Take $T \in (0, \infty ) .$
Then all the conclusions stated in Corollary \ref{corollaryC10.65} for the function $\Psi$ introduced
there hold for $\Psi:=\mbox{$\mathfrak I$} ^{( \tau )} (a)|S_{\infty }$ as well.
\end{corollary}
{\bf Proof:} Theorem \ref{theoremT10.6}, Lemma \ref{lemmaL10.100} and \ref{lemmaL10.20}. \hfill $\Box $

\vspace{1ex}
We will need the following pointwise spatial decay estimate of $\mbox{$\mathfrak I$} ^{( \tau )} (a)$.
\begin{theorem} \label{theoremT4.60}
Let $R_0,\, \delta _0 \in (0, \infty ) ,\; \kappa _0 \in (0,1],\; a \in L_{loc}^1( \mathbb{R}^3 )^3$
such that
$a| \overline{ B_{R_0}}^c \in W^{1,1}_{loc}(  \overline{ B_{R_0}}^c)^3,\;
| \partial ^{\alpha }a(y)|\le \delta _0\, \bigl(\, |y|\, \nu (y) \,\bigr) ^{-1-| \alpha |/2-\kappa _0}
$
for $y \in  \overline{ B_{R_0}}^c,\; \alpha \in \mathbb{N} _0^3$ with $| \alpha |\le 1.$

Let $R \in (R_0, \infty )$ and take $\alpha $ as before. Then, for $x \in B_R^c$ and $t \in (0, \infty ) ,$
\begin{eqnarray*}
| \partial _x^{\alpha }\mbox{$\mathfrak I$} ^{( \tau )} (a)(x,t)|
\le
\mbox{$\mathfrak C$} \, ( \delta _0+\|a|B_{R_0}\|_1)\, \bigl(\, |x|\, \nu (x) \,\bigr) ^{-1-| \alpha |/2}.
\end{eqnarray*}
\end{theorem}
{\bf Proof:} See \cite[Theorem 1.1]{De7}. Note that in \cite[inequality (1.10)]{De7}, it should read $| \beta |/2$
instead of $| \beta |$. \hfill $\Box $

\vspace{1ex}
We turn to a single layer potential whose definition is based on the ensuing
\begin{lemma} \label{lemmaL10.110}
Let $\phi \in L^2(S_{\infty })^3 ,\; \alpha \in \mathbb{N} _0^3$ with $| \alpha |\le 1.$
Then, for $x \in \mathbb{R}^3 \backslash \partial \Omega ,\;
t \in (0, \infty ) $, we have
$
\int_{ 0}^t \int_{ \partial \Omega }| \partial _x^{\alpha } \Lambda (x-y,t-\sigma , \tau ) \cdot \phi (y, \sigma )|\,
do_y\, d \sigma < \infty .
$
In addition, the preceding relation holds for a. e. $x \in \partial \Omega $ and a. e. $t \in (0, \infty ) $.
\end{lemma}
{\bf Proof:} 
If $x \in \mathbb{R}^3 \backslash \partial \Omega ,$ we have
$|x-y|\ge dist(x, \partial \Omega )>0$ for $y \in \partial \Omega $, so
Corollary \ref{corollaryC10.60} with
$K=dist(x, \partial \Omega )/2$
yields
$| \partial _x^{\alpha } \Lambda _{jk} (x-y,t-\sigma , \tau )|\le \mbox{$\mathfrak C$} \,
\bigl(\, dist(x, \partial \Omega ) \,\bigr) ^{-3/2-| \alpha |/2}>0$
for
$y$
as before,
$t \in (0, \infty ) $ and $\sigma \in (0,t)$.
This estimate and H\"older's inequality imply the first claim of the lemma.
The second holds according to \cite[Lemma 2.19]{DeSIMA1}.
\hfill $\Box $

\vspace{1ex}
In view of the preceding lemma, we may define
\begin{eqnarray*}
\mbox{$\mathfrak V$}^{( \tau )} ( \phi )(x,t)
:=
\int_{ 0}^t \int_{ \partial \Omega } \Lambda (x-y,t-\sigma , \tau ) \cdot \widetilde{ \phi }(y, \sigma ) \,
do_y\, d \sigma
\end{eqnarray*}
for $T \in (0, \infty ],\; \phi \in L^2( S_T)^3,\; x \in \mathbb{R}^3 \backslash \partial \Omega ,\;
t \in (0, \infty ),$
and for a. e. $x \in \partial \Omega $ and a. e. $t \in (0, \infty ) $,
where $\widetilde{ \phi }$ denotes the zero extension of $\phi $ to $S_{\infty }$.

Note that $\mbox{$\mathfrak V$} ^{( \tau )} ( \phi )$ is defined as a function on $( \mathbb{R}^3 \backslash
\partial \Omega ) \times (0, \infty )$ and also as a function on $\partial \Omega \times (0, \infty )$.
The second of these functions is the trace of the first:
\begin{theorem} \label{theoremT10.10}
Let $\phi \in L^2(S_{\infty } )^3$ and abbreviate $v:= \mbox{$\mathfrak V$} ^{( \tau )} ( \phi )
| \overline{ \Omega }^c \times (0, \infty ).$
Then $v \in W_{loc}^{1,1} \bigl(\,  \overline{ \Omega }^c \times (0, \infty ) \,\bigr) ^3
\cap
C^0 \bigl(\,  \overline{ \Omega }^c \times  [0, \infty ) \,\bigr) ^3,\;
v(t) \in C ^{ \infty } ( \overline{ \Omega }^c)^3$ and
\begin{eqnarray} \label{T10.10.10}
\partial _lv(x,t)
=
\int_{ 0}^t \int_{ \partial \Omega } \partial _l\Lambda (x-y,t-\sigma , \tau ) \cdot  \phi (y, \sigma ) \,
do_y\, d \sigma
\end{eqnarray}
for
$t \in (0, \infty ) ,\; x \in \overline{ \Omega }^c,\; 1\le l\le 3.$
Moreover
$
v \in L^{\infty } \bigl(\,  0, \infty ,\, L^2( \overline{ \Omega  }^c)^3 \,\bigr) ,\;
\nabla _xv \in L^{2 } \bigl(\, \overline{ \Omega  }^c \times (0, \infty )\,\bigr)^9,
$
and for $t \in (0, \infty ) $, the trace of $v(t)$ coincides with $\mbox{$\mathfrak V$} ^{( \tau )} ( \phi )(t)$
considered as a function on $\partial \Omega $.
\end{theorem}
{\bf Proof:}
For the first part of the lemma, up to and including (\ref{T10.10.10}), we refer to
\cite[Lemma 2.21]{DeSIMA1}. The second part holds according to \cite[Theorem 2.3, 2.4]{De3}.
\hfill $\Box $

\vspace{1ex}
If $\phi \in L^2(S_T)^3$ for some $T<\infty $, then $\mbox{$\mathfrak V$} ^{( \tau )} ( \phi )$
is smooth on $\mathbb{R}^3 \times (T, \infty )$:
\begin{lemma} \label{lemmaL10.120}
Let $T \in (0, \infty ) $ and $\phi \in L^2(S_T)^3. $
Then $\mbox{$\mathfrak V$} ^{( \tau )} ( \phi )| \mathbb{R}^3 \times (T, \infty )$
belongs to
$C^1 \bigl(\, \mathbb{R}^3 \times  (T, \infty ) \,\bigr) ^3$, and
\begin{eqnarray*}
\partial _t^l \partial _x^{\alpha } \mbox{$\mathfrak V$} ^{( \tau )} ( \phi )(x,t) 
=
\int_{ 0}^T \int_{ \partial \Omega } \partial _t^l \partial _x ^{\alpha }
\Lambda (x-y,t-\sigma , \tau ) \cdot  \phi (y, \sigma ) \,
do_y\, d \sigma
\end{eqnarray*}
for
$
x \in \mathbb{R}^3 ,\; t \in (T, \infty ),\; \alpha \in \mathbb{N} _0^3,\; l \in \mathbb{N} _0$
with $| \alpha |+l\le 1,$
where the preceding integral exists as Lebesgue integral.
\end{lemma}
{\bf Proof:}
If 
$x \in \mathbb{R}^3 ,\; t \in (T, \infty )$
in the situation of Lemma \ref{lemmaL10.120},
the time integral in the definition of
$\mbox{$\mathfrak V$} ^{( \tau )} ( \phi )(x,t)$ only extends from $0$ to $T$, and not from $0$
to $t$. Thus the lemma follows from Corollary \ref{corollaryC10.60} and Lebesgue's theorem
by a similar but simpler reasoning as in the proof of Lemma \ref{lemmaL10.85}, simpler in particular
because $L^2(S_T)^3 \subset L^1(S_T)^3$.
\hfill $\Box $

\vspace{1ex}
We will need the following, much more deep-lying properties of $\mbox{$\mathfrak V$} ^{( \tau )} ( \phi )$.
\begin{theorem} \label{theoremT10.20}
The relation $ \mbox{$\mathfrak V$}^{( \tau )} ( \phi )|S_{\infty   }\in H_{\infty }$ holds
for $\phi \in L^2(S_{\infty })^3$. There is a constant $c_1>0$ such that
$\|\mbox{$\mathfrak V$}^{( \tau )} ( \phi )|S_{\infty   }\|_{ H_{\infty }}\le c_1\, \| \phi \|_2$
for such $\phi  $.
\end{theorem}
{\bf Proof:}
For $\phi  \in L^2( S_{\infty })^3,$ define $\mbox{$\mathfrak V$}^{( 0 )} ( \phi  )$
in the same way as $\mbox{$\mathfrak V$}^{( \tau )} ( \phi  )$, but with the term
$\Lambda (x-y,t-\sigma , \tau )$ replaced by $\Gamma (x-y, t-\sigma )$;
compare \cite[Lemma 2.19]{DeSIMA1} and use Lemma \ref{lemmaL10.50}.
According to \cite[Theorem 2.3.3]{Shen} and the estimates in \cite[p. 362]{Shen},
we have $\mbox{$\mathfrak V$}^{( 0 )} ( \phi  )|S_{\infty } \in H_{\infty } $
for $\phi  \in L^2(S_{\infty })^3,$ and there is $c>0$ with
$\|\mbox{$\mathfrak V$}^{( 0 )} ( \phi  )|S_{\infty   }\|_{ H_{\infty }}\le c\, \| \phi \|_2$
for $\phi $ as before. The theorem follows with \cite[Theorem 4]{DeEst}.
\hfill $\Box $
\begin{theorem}[\mbox{\cite[Corollary 3]{DeEst}}] \label{theoremT10.30}\label{theoremT10.40}
There is $c_2>0$ such that $\| \phi \|_2\le \| \mbox{$\mathfrak V$} ^{( \tau )} ( \phi )|S_{\infty }\|_{H_{\infty }}$
for $\phi \in L^2_n(S_{\infty  }).$

For $\widetilde{ b} \in H_{\infty }$,
the integral equation
$\mbox{$\mathfrak V$} ^{( \tau )} ( \phi )|S_{\infty } = \widetilde{ b}$
is solved by a unique function $\phi \in L^2_n(S_{\infty })$.
(The function space $ L^2_n(S_{ \infty  })$ was introduced in (\ref{2.10})).
\end{theorem}
\begin{corollary} \label{corollaryC10.70} Let $T \in (0, \infty ) $ and $\phi \in L^2( S_T)^3$.
Then all the conclusions listed in Corollary \ref{corollaryC10.65} for the function $\Psi$
introduced there are true for $\Psi= \mbox{$\mathfrak V$}^{( \tau )} ( \phi )| S_{\infty }$ as well.
\end{corollary}
{\bf Proof:}
Lemma \ref{lemmaL10.120}, Theorem \ref{theoremT10.20} and Lemma \ref{lemmaL10.10}. \hfill $\Box $
\begin{corollary} \label{corollaryC10.80}
Let $ \phi \in L^2(S_{\infty }),\; T \in (0, \infty ) .$ Then,
with $c_1$ from Theorem \ref{theoremT10.20}, we have
$\| \mbox{$\mathfrak V$} ^{( \tau )} ( \phi )|S_{T, \infty }\|_{H_{T, \infty }}
\le
\| \mbox{$\mathfrak V$} ^{( \tau )} ( \phi )\|_{H_{\infty }}
\le c_1\,\| \phi \|_2.
$
\end{corollary}
{\bf Proof:}
Theorem \ref{theoremT10.20}, Corollary \ref{corollaryC10.40}. \hfill $\Box $
\begin{corollary} \label{corollaryC10.90}
Let $T \in (0, \infty ) $ and $\phi  \in L^2_n( S_{\infty })$ with $\phi |S_T =0$.
Then, with $c_2$ from Theorem \ref{theoremT10.30},
$\| \phi \|_2
\le
c_2\,
\| \mbox{$\mathfrak V$} ^{( \tau )} ( \phi )\|_{H_{\infty }}
=
c_2\,
\| \mbox{$\mathfrak V$} ^{( \tau )} ( \phi )|S_{T, \infty }\|_{H_{T, \infty }}.
$
\end{corollary}
{\bf Proof:}
The inequality stated  in the corollary holds according to Theorem \ref{theoremT10.30},
whereas the equation follows from Corollary \ref{corollaryC10.40}.
\hfill $\Box $

\vspace{1ex}
In view of a representation formula for solutions to (\ref{10}), (\ref{20b}), (\ref{30}),
we first state a uniqueness theorem for such solutions.
\begin{theorem} \label{theoremT10.50}
Let
$
b:S_{\infty }\mapsto \mathbb{R}^3 ,\; 
a \in L^{1}_{loc}( \overline{ \Omega }^c)^3,\;
f \in L^2_{loc} \bigl(\, [0, \infty ),\, [H_{\sigma }^1( \overline{ \Omega }^c) ] ^{\prime} ).
$
Then there is at most one function
$u \in L^2_{loc} \bigl(\, [0, \infty ) ,\, H^1( \overline{ \Omega }^c)^3 \,\bigr) $
such that
\begin{eqnarray} && \label{T10.50.9} 
u(t)| \partial \Omega  = b(t) \; \mbox{for}\;\; t \in (0, \infty ) ,
\quad 
\mbox{div}_x u =0,
\\ && \label{T10.50.10}
\int_{ 0}^{\infty } \Bigl( \int_{ \overline{ \Omega}^c} \Bigl[ - \bigl(\, u(x,t) \cdot V (x) \,\bigr) 
\, \varphi ^{\prime} (t)
\\ && \hspace{2ex} \nonumber  
+ \bigl(\, \nabla _x u(x,t) \cdot \nabla V (x) \,\bigr)  \, \varphi (t)
+ \tau \, \bigl(\,  \partial _1u(x,t) \cdot V (x) \,\bigr) \, \varphi (t) \Bigr]\, dx
- f(t)( V ) \, \varphi (t) \Bigr) \, dt
\\ && \nonumber 
=
\int_{ \overline{ \Omega}^c} a(x) \cdot V (x) \, dx \; \varphi (0)
\end{eqnarray}
for
$ \varphi \in C ^{ \infty } _0 \bigl(\, [0, \infty ) \,\bigr) ,\; V \in
C ^{ \infty } _0( \overline{ \Omega }^c)^3$
with
$\mbox{div}\, V =0.$
\end{theorem}
{\bf Proof:}
This theorem may be shown in the same way as an analogous result for the Stokes system.
We refer to \cite[Theorem 3.7]{DeJMFM} and its proof.
\hfill $\Box $

\vspace{1ex}
Now we construct a solution to (\ref{10}), (\ref{20b}), (\ref{30}) in the form
$u = \mbox{$\mathfrak R$} ^{( \tau )} (f) + \mbox{$\mathfrak I$} ^{( \tau )} (a) + \mbox{$\mathfrak V$}
^{( \tau )} ( \phi )$,
with $\phi $ being the solution to the integral equation (\ref{T10.60.10}) below.
In this way we obtain a representation formula for the velocity as a sum of three potential functions,
as announced in section 1. We consider this solution as a weak one in the sense of Theorem
\ref{theoremT10.50} because we want to range it in a uniqueness class which is as large as possible.
If $f$ and $a$ are smooth, then $u$ is the velocity part of a solution satisfying (\ref{10}) and
(\ref{30}) pointwise; see \cite[Lemma 2.14, Theorem 2.16, Lemma 2.21]{DeSIMA1}.
(In \cite[Theorem 2.16]{DeSIMA1}, there is a misprint: $\mbox{$\mathfrak R$} ^{( \tau )} (f)$
is only a $C ^{ \infty } $- and not a $C ^{ \infty } _0$-function.)
The boundary condition on
$\partial \Omega $ holds according to (\ref{T10.50.9}), and the boundary condition at infinity
is satisfied
in the sense that $u(t) \in H^1( \overline{ \Omega }^c)^3$ for $t \in (0, \infty ) $.
If the functions $f,\, a $ and $b$ decay sufficiently rapidly, this boundary condition will
be fulfilled in a stronger way, on the basis
of Theorem \ref{theoremT4.10}, \ref{theoremT4.60} and Lemma \ref{lemmaL5.30} below.
We will come back to this subject in section 8.
\begin{theorem} \label{theoremT10.60}
Let $\epsilon _0 \in (0,\, 1/2 ],\; q \in [1, \, 3/2 ),\; b \in H_{\infty },\;
a \in H ^{1/2+\epsilon _0}_{\sigma } ( \overline{ \Omega }^c),\;
f \in L^2 \bigl(\, 0, \infty ,\; L^q( \overline{ \Omega }^c)^3 \,\bigr)
\cap
L^2 \bigl(\, 0, \infty ,\; L ^{3/2} ( \overline{ \Omega }^c)^3 \,\bigr).
$
Then there is a unique function $\phi \in  L^2_n( S_{\infty })$ (see (\ref{2.10}))
with
\begin{eqnarray} \label{T10.60.10}
\mbox{$\mathfrak V$} ^{( \tau )} ( \phi )|S_{\infty }= - \mbox{$\mathfrak R$} ^{( \tau )} (f)
-\mbox{$\mathfrak I$} ^{( \tau )} (a)+b.
\end{eqnarray}
Moreover there is a unique function $u \in L^2_{loc} \bigl(\, [0, \infty ),\,
H^1( \overline{ \Omega }^c)^3 \,\bigr) $
such that (\ref{T10.50.9}) and (\ref{T10.50.10}) hold. This function is given by
\begin{eqnarray} \label{T10.60.20}
u = \mbox{$\mathfrak R$} ^{( \tau )} (f) + \mbox{$\mathfrak I$} ^{( \tau )} (a) + \mbox{$\mathfrak V$}
^{( \tau )} ( \phi )| \overline{ \Omega }^c \times (0, \infty ).
\end{eqnarray}
\end{theorem}
{\bf Proof:}
Theorem \ref{theoremT10.60} holds according to \cite[Theorem 2.26]{DeDCDS-A}.
Note that the unique solvability of (\ref{T10.60.10}) is a consequence of Theorem \ref{theoremT10.40},
\ref{theoremT10.6} and \ref{theoremT10.5}.
\hfill $\Box $

\section{Temporal decay of the potential $\mbox{$\mathfrak R$} ^{( \tau )} (f)$.}
\setcounter{equation}{0}

This section has two aims. Firstly, we want to estimate
$\| \mbox{$\mathfrak R$} ^{( \tau )} (f)|S_{T, \infty }\|_{H_{T, \infty }}$,
and secondly, we are going to determine an upper bound 
of $| \partial _x^{\alpha }\mbox{$\mathfrak R$} ^{( \tau )} (f)(x,t)|$
under the assumption $| \alpha |\le 1$.
Our starting point will be to split $f$ into a sum $f ^{(1)} +f ^{(2)} $,
with
$f ^{(1)} =\chi_{\mathbb{R}^3} \times (t/2, \infty )$
and 
$f ^{(2)} =\chi_{\mathbb{R}^3} \times (0, t/2 )$.
The decay of $\mbox{$\mathfrak R$} ^{( \tau )} (f ^{(1)} )(x,t)$ is then due to the asymptotic behaviour of $f$
for $t\to \infty $, whereas $\mbox{$\mathfrak R$} ^{( \tau )} (f ^{(2)} )(x,t)$ becomes small for large $t$
due to the decay properties of the fundamental solution $\Lambda $. The next lemma
addresses some key technical difficulties in the estimate of
$\| \mbox{$\mathfrak R$} ^{( \tau )} (f ^{(2)} )|S_{T, \infty }\|_{H_{T, \infty }}$ 
and also of
$\mbox{$\mathfrak R$} ^{( \tau )} (f ^{(2)} )(x,t)$.
\begin{lemma} \label{lemmaL6.10}
Let $\alpha \in \mathbb{N} _0^3,\; l \in \mathbb{N} _0$ with $| \alpha |+l\le 1.$
Let $q \in [1,4],\; s \in [1, \infty )$ with $3/(2\, q)+1/s >1-| \alpha |/2-l/2.$
Take $T \in (0, \infty ) $ and $f \in L^s \bigl(\, 0,T,\, L^q ( \mathbb{R}^3 )^3 \,\bigr) .$
Then
\begin{eqnarray} \label{L6.10.10} &&
| \partial _t^l \partial _x^{\alpha } \mbox{$\mathfrak R$} ^{( \tau )} (f)(x,t)|
\le
\mbox{$\mathfrak C$} \,\|f\|_{q,s;T}
\,
\bigl(\, (t-T)^{-3/(2\, q)-1/s +1-| \alpha |/2-l}
\\&&\nonumber \hspace{4em}
+\delta _{1l}\, (t-T)^{-3/(2\, q)-1/s +1/2}\,\bigr)
\quad \mbox{for}\;\;
x \in \mathbb{R}^3 ,\; t \in (T, \infty ) .
\end{eqnarray}
(Recall that
$\mbox{$\mathfrak R$} ^{( \tau )} (f)| \mathbb{R}^3 \times (T, \infty ) \in C^1 \bigl(\, \mathbb{R}^3 \times
(T, \infty ) \,\bigr) ^3$
by Lemma \ref{lemmaL10.85}.)
\end{lemma}
{\bf Proof:}
Take $x$ and $t$ as in the lemma. In the case $l=0$, hence $| \alpha |\le 1$,
the exponents $q$ and $s$ verify the condition
$3/(2\, q)+1/s >1-| \alpha |/2.$
This and the inequality $t-\sigma >(t-T)/2$ for $\sigma \in (0,T)$ allow us to apply
Theorem \ref{theoremT4.30} with $M=(t-T)/2$ and $ \varrho = \infty .$
Lemma \ref{lemmaL10.85} and Theorem \ref{theoremT4.30}
then yield (\ref{L6.10.10}) in the case $l=0$.

Let us suppose for the rest of this proof that $l=1$. This means in particular that
$\alpha =0$ and $3/(2\, q)+1/s>1/2.$
Consider the case $s>1$ and $q>1$. Then by Lemma \ref{lemmaL10.85}, H\"older's inequality
and Lemma \ref{lemmaL10.50},
\begin{eqnarray} \label{L6.10.20}&&
| \partial _t \mbox{$\mathfrak R$} ^{( \tau )} (f)(x,t)|
\\&&\nonumber 
\le
C\, \int_{ 0}^T \int_{ \mathbb{R}^3  } \sum_{j \in \{5,\, 4\}}
\bigl(\, |x-y-\tau \,(t- \sigma )\, e_1| + (t-\sigma ) ^{1/2} \,\bigr) ^{-j}\, |f(y, \sigma )|\, dy\, d \sigma
\\&&\nonumber
\le
C\,
\|f\|_{q,s;T}\,\sum_{j \in \{5,\, 4\}}\Bigl( 
\int_{ 0}^T \Bigl[ \int_{ \mathbb{R}^3  } 
\bigl(\, |x-y-\tau \,(t- \sigma )\, e_1| + (t-\sigma ) ^{1/2} \,\bigr) ^{-j\, q ^{\prime} }\,
dy \Bigr]^{s ^{\prime} /q ^{\prime} }\, d \sigma  \Bigr) ^{1/s ^{\prime} }.
\end{eqnarray}
But for any $\sigma \in (0,T),\; j \in \{5,\, 4\},$ the preceding integral with respect to $y$ equals
$\int_{ \mathbb{R}^3 }\bigl(\, |z|+(t-\sigma ) ^{1/2} \,\bigr) ^{-j\, q ^{\prime} }\, dy,$
and is hence bounded by
$C(q)\, (t-\sigma )^{-j\, q ^{\prime} /2+3/2}.$
Therefore from (\ref{L6.10.20}),
\begin{eqnarray*} 
| \partial _t \mbox{$\mathfrak R$} ^{( \tau )} (f)(x,t)|
\le
C(q)\,
\|f\|_{q,s;T}\,\sum_{j \in \{5,\, 4\}}\Bigl( 
\int_{ 0}^T (t-\sigma )^{[-j/2+3/(2q ^{\prime} )]\,s ^{\prime} }\, d \sigma \Bigr) ^{1/s ^{\prime} }.
\end{eqnarray*}
But
$-j/2+3/(2q ^{\prime} )=(-j+3)/2-3/(2\,q)$ for $j \in \{5,\, 4\} .$
Since $3/(2\, q)+1/s >1/2,$ as mentioned above, we have
$[(-j+3)/2-3/(2\,q)]\, s ^{\prime} <-1$
for $j$ as before, so we may conclude that
\begin{eqnarray*} 
| \partial _t \mbox{$\mathfrak R$} ^{( \tau )} (f)(x,t)|
\le
C(q,s)\,
\|f\|_{q,s;T}\,\sum_{j \in \{5,\, 4\}}
(t-T )^{(-j+3)/2-3/(2q)+1/s ^{\prime} }.
\end{eqnarray*}
Thus we have shown (\ref{L6.10.10}) in the case $s>1,\; q>1$, under the assumption $l=1$, as
we may recall. If $s=1,\; q>1$, it follows as in (\ref{L6.10.20}) that
$| \partial _t \mbox{$\mathfrak R$} ^{( \tau )} (f)(x,t)|$
is bounded by
\begin{eqnarray*}
C\,
\sum_{j \in \{5,\, 4\}}\int_{ 0}^T\Bigl( 
\int_{ \mathbb{R}^3  } 
\bigl(\, |x-y-\tau \,(t- \sigma )\, e_1| + (t-\sigma ) ^{1/2} \,\bigr) ^{-j\, q ^{\prime} }\,
dy \Bigr)^{1 /q ^{\prime} }\, \|f(\, \cdot \, , \sigma )\|_q\, d \sigma .
\end{eqnarray*}
Since $(t- \sigma ) ^{1/2} \ge (t-T) ^{1/2} $ for $\sigma \in (0,T),$
we may conclude that the term $| \partial _t \mbox{$\mathfrak R$} ^{( \tau )} (f)(x,t)|$
may be estimated by
$C(q)\, \|f\|_{q,1; \infty }\sum_{j \in \{5,\, 4\}}(t-T)^{-j/2+3/(2 q ^{\prime} )},$
which is the looked-for result in the case $s=1,\; q>1$.
If $s>1,\; q=1,$ inequality (\ref{L6.10.20}) is replaced by the estimate
\begin{eqnarray*} &&
| \partial _t \mbox{$\mathfrak R$} ^{( \tau )} (f)(x,t)|
\le
C\,\sum_{j \in \{5,\, 4\}} \int_{ 0}^T \int_{ \mathbb{R}^3  } 
(t-\sigma ) ^{-j/2}  |f(y, \sigma )|\, dy\, d \sigma
,
\end{eqnarray*}
so that
$
| \partial _t \mbox{$\mathfrak R$} ^{( \tau )} (f)(x,t)|
\le
C\,\|f\|_{1,s;T}\,\sum_{j \in \{5,\, 4\}}\bigl( 
\int_{ 0}^T (t-\sigma )^{-j \, s ^{\prime} /2 }\, d \sigma  \bigr) ^{1/s ^{\prime} }.
$
Inequality (\ref{L6.10.10}) with $q=1,\; s>1$ follows.
It is obvious how to evaluate
$| \partial _t \mbox{$\mathfrak R$} ^{( \tau )} (f)(x,t)|$ if $s=q=1$.
\hfill $\Box $

\vspace{1ex}
Further below (Corollary \ref{corollaryC6.20}), we will estimate
$\| \mbox{$\mathfrak R$} ^{( \tau )} (f)|S_{T+\mu,\,  \infty } \|_{H_{T+\mu,\, \infty }}$
in terms of negative powers of $\mu$, under the assumption that $f| \mathbb{R}^3 \times (T, \infty ) =0.$
In the next three lemmas, we derive this type of bound.
\begin{lemma} \label{lemmaL6.15}
Let $q,\,s,\, T,\, f$ be given as in Lemma \ref{lemmaL6.10}. Suppose in addition
that $3/(2\, q)+1/s> 3/2.$
Let $\mu \in (0, \infty ) $. Then
\begin{eqnarray*}&&
\Bigl( \int_{ T+\mu } ^{ \infty }
\Bigl[ \| \mbox{$\mathfrak R$} ^{( \tau )} (f)(t)| \partial \Omega \|^2_{H^1( \partial \Omega )^3}
+
\int_{ \partial \Omega }
|n^{ (\Omega )} \cdot \partial _t \mbox{$\mathfrak R$} ^{( \tau )} (f)(x,t)|^2\, do_x
\Bigr]\, dt \Bigr) ^{1/2}
\\&&
\le
\mbox{$\mathfrak C$} \,
\|f\|_{q,s;T}\, (\mu^{-3/(2\, q)-1/s+3/2}+\mu^{-3/(2\, q)-1/s+1/2}).
\end{eqnarray*}
(Recall that $\mbox{$\mathfrak R$} ^{( \tau )} (f)| \mathbb{R}^3 \times (T, \infty )$
is a $C^1$-function by Lemma \ref{lemmaL10.85}.)
\end{lemma}
{\bf Proof:}
Since $3/(2 \, q)+1/s > 3/2,$ the term
$\bigl(  \int_{ T+\mu}^{ \infty } \bigl[ (t-T)^{-3/q-2/s+2}+(t-T)^{-3/ q-2/s} \bigr] \, dt\, \bigr) ^{1/2} $
is bounded by
$\mbox{$\mathfrak C$} \, (\mu^{-3/(2\, q)-1/s+3/2}+\mu^{-3/(2\, q)-1/s+1/2})$.
Thus the lemma follows from Lemma \ref{lemmaL6.10}. Note that the largest exponent $-3/(2 \, q)-1/s+3/2$
arises in the case $\alpha =0,\; l=0$ in Lemma \ref{lemmaL6.10}, and the smallest one
$-3/(2 \, q)-1/s+1/2$ if $\alpha =0,\; l=1$.
\hfill $\Box $

\vspace{1ex}
In order to deal with $\partial _t ^{1/2} \mbox{$\mathfrak R$} ^{( \tau )} (f)(x,t)$,
we need a preparatory result:
\begin{lemma} \label{lemmaL6.20}
Let $q \in (1,2],\; F \in L^q( \mathbb{R}^3 )^3$ and $ d \in (0,3]$.
Then, for $r \in (0, \infty ) $,
\begin{eqnarray*}
\Bigl( \int_{ \partial \Omega }\Bigl[ \int_{ \mathbb{R}^3 }| \Lambda (x-y,r, \tau ) \cdot F(y)|\, dy \Bigr]^2
\, do_x \Bigr) ^{1/2}
\le
\mbox{$\mathfrak C$} \,
\|F\|_q\,
\max\{1,\: r^{-d/(2q ^{\prime} )-3/(2\, q)+1/2}\}.
\end{eqnarray*}
\end{lemma}
{\bf Proof:}
Choose $R_0>0$ so large that $\overline{ \Omega }\subset B_{R_0/2}$.
This means that $B_{R_0}^c \subset B_{R_0/2}(x)^c$ and $B_{R_0} \subset B_{2R_0}(x) $
for $x \in \partial \Omega $.

Let \mbox{$\mathfrak A$} denote the left-hand side of the estimate stated in the lemma,
but with the integral over $\mathbb{R}^3 $ replaced by one over $B_{R_0}$.
Take $r \in (0, \infty ) $. Then, by H\"older's inequality we see that $\mbox{$\mathfrak A$}^2$
is bounded by
$$
\mbox{$\mathfrak C$} \, \sum_{j,k = 1}^3 
\int_{ \partial \Omega } \Bigl[ \Bigl( \int_{ B_{R_0}}
|\Lambda _{jk} (x-y,r, \tau )|\, dy \Bigr) ^{1/q ^{\prime} }
\Bigl( \int_{ B_{R_0}}  |\Lambda _{jk} (x-y,r, \tau ) |\, |F(y)|^q\, dy \Bigr) ^{1/q  }\Bigr]^2\, do_x
.
$$
But for
$x \in \partial \Omega ,\; 1\le j,k\le 3 $,
by Corollary \ref{corollaryC10.60} with $K=2\, R_0$
and because $B_{R_0}\subset B_{2\,R_0}(x)$,
\begin{eqnarray*}&&
\int_{ B_{R_0}}|\Lambda _{jk} (x-y,r, \tau )|\, dy
\le
C(R_0)\,
\int_{ B_{R_0}} (|x-y|+r ^{1/2} ) ^{-3} \, dy
\\&&
\le
C(R_0)\, r^{-d/2}\, \int_{ B_{2R_0}(x)} (|x-y|+r ^{1/2} ) ^{-3+d} \, dy
\le
C(R_0,d)\, r^{-d/2}.
\end{eqnarray*}
Therefore
$
\mbox{$\mathfrak A$}
\le
\mbox{$\mathfrak C$} \,
 r^{-d/(2q ^{\prime} )}\,\sum_{j,k = 1}^3 
\bigl(\,
\int_{ \partial \Omega } 
\bigl[\, 
\int_{ B_{R_0}}  |\Lambda _{jk} (x-y,r, \tau ) |\, |F(y)|^q\, dy
\,\bigr]
^{2/q } do_x
\,\bigr) 
^{1/2} .
$
Since $2/q\ge 1$,
we may conclude with Theorem \ref{theoremT1.10} that
\begin{eqnarray*}
\mbox{$\mathfrak A$}
\le
\mbox{$\mathfrak C$} \,
 r^{-d/(2 q ^{\prime} )}\, \sum_{j,k = 1}^3 
\biggl(\, \int_{ B_{R_0}}
\Bigl( \int_{ \partial \Omega }   |\Lambda _{jk} (x-y,r, \tau ) |^{2/q}\, do_x \Bigr) ^{q/2 }
\, |F(y)|^q\, dy  \,\biggr) ^{1/q} .
\end{eqnarray*}
On the other hand, again by Corollary \ref{corollaryC10.60} with $K=2\, R_0$,
for $y \in B_{R_0},\; 1\le j,k\le 3$,
\begin{eqnarray*}
\int_{ \partial \Omega }   |\Lambda _{jk} (x-y,r, \tau ) |^{2/q}\, do_x
\le
\mbox{$\mathfrak C$}  \,
\int_{ \partial \Omega }   (|x-y|+r ^{1/2} )^{-6/q}\, do_x
\le
\mbox{$\mathfrak C$} \,r^{-3/q+1}.
\end{eqnarray*}
As a consequence,
$\mbox{$\mathfrak A$} \le \mbox{$\mathfrak C$} \,r^{-d/(2q ^{\prime} )-3/(2\, q)+1/2}\, \|F\|_q$.
Let \mbox{$\mathfrak B$} be defined also by the left-hand side of the estimate
stated in the lemma, but this time with the integral over $\mathbb{R}^3 $ replaced
by one over $B_{R_0}^c$. We get
\begin{eqnarray} \label{L6.20.30}
\mbox{$\mathfrak B$}
\le
\sum_{j,k = 1}^3 \biggl( \int_{ \partial \Omega } \Bigl[ \Bigl( \int_{ B_{R_0}^c}
|\Lambda _{jk} (x-y,r, \tau ) |^{q ^{\prime} }\, dy \Bigr) ^{1/q ^{\prime}  }\, 
\|F\|_q \Bigr]^2\, do_x \biggr) ^{1/2}.
\end{eqnarray}
In order to estimate the integral over $B_{R_0}^c$ in (\ref{L6.20.30}), we recall that
$B_{R_0}^c \subset B_{R_0/2}(x)^c$ for $x \in \partial \Omega $,
so we may use Corollary \ref{corollaryC10.60} with $K=R_0/2$ to obtain
\begin{eqnarray*}&&
\int_{ B_{R_0}^c}
|\Lambda _{jk} (x-y,r, \tau ) |^{q ^{\prime} }\, dy
\le
C(R_0,q)\,
\int_{ B_{R_0/2}(x)^c} \bigl(\, |x-y|\, \nu(x-y) \,\bigr)^{-3\, q ^{\prime} /2} \, dy
\\&&
\le
C(R_0,q)\,
\int_{ B_{R_0/2}^c} \bigl(\, |z|\, \nu(z) \,\bigr)^{-3\, q ^{\prime} /2} \, dz
\quad \mbox{for}\;\;
x \in \partial \Omega ,\; 1\le j,k\le 3.
\end{eqnarray*}
But $q\le 2$, in particular $q<4$, so $3\, q ^{\prime} /2 >2,$
hence with Corollary \ref{corollaryC1.10} we arrive at the inequality
$
\int_{ B_{R_0}^c}
|\Lambda _{jk} (x-y,r, \tau ) | ^{q ^{\prime}  }\, dy
\le
C(R_0,q).
$
This estimate is inserted into (\ref{L6.20.30}). We then obtain
$\mbox{$\mathfrak B$} \le \mbox{$\mathfrak C$} \,\|F\|_q.$
The lemma now follows with the estimate of \mbox{$\mathfrak A$} obtained above.
\hfill $\Box $
\begin{lemma} \label{lemmaL6.30}
Let $s \in (1, \infty ),\; q \in (1,\, 3/2)$
with $3/(2\, q)+1/s >3/2.$
Let
$
T, \mu \in (0, \infty ) ,\;
f \in L^s \bigl(\, 0,T,\, L^q( \mathbb{R}^3 )^3 \,\bigr) .
$
Then
$$
\Bigl( \int_{ T+\mu } ^{ \infty } \int_{ \partial \Omega }
| \partial _t ^{1/2} \mbox{$\mathfrak R$} ^{( \tau )} (f)(x,t)|^2\, do_x\, dt \Bigr) ^{1/2}
\le
\mbox{$\mathfrak C$} \,
\|f\|_{q,s;T}\, \Bigl( \sum_{j \in \{2,\, 3\}} \mu^{-3/(2\, q)-1/s+j/2}
+\mu^{-1/s}\Bigr) .
$$
(According to Corollary \ref{corollaryC10.65}, the fractional derivative
$\partial _t ^{1/2} \mbox{$\mathfrak R$} ^{( \tau )} (f)(x,t)$ is well defined for any
$t \in (T, \infty )$ and a. e. $x \in \partial \Omega $.)
\end{lemma}
{\bf Proof:}
For brevity we set $\Psi:= \mbox{$\mathfrak R$} ^{( \tau )} (f)|S_{\infty }.$
According to Corollary \ref{corollaryC10.65}, for a. e. $x \in \partial \Omega $
equation (\ref{L10.10.10}) is verified with $\Psi (x,\, \cdot \,)$ in the role of $\phi $.
We will estimate each term on the right-hand side of (\ref{L10.10.10}) in the norm
of $L^2(S_{T+\mu, \infty })^3.$ To begin with, we use Lemma \ref{lemmaL6.10} to obtain
\begin{eqnarray} \label{L6.30.10}&&
\Bigl( \int_{ T+\mu} ^{ \infty } \int_{ \partial \Omega }| (t-T) ^{-1/2} \,
\Psi \bigl(\, x,\, (t+T)/2 \,\bigr) |^2\, do_x\, dt \Bigr) ^{1/2}
\\&&\nonumber
\le
\mbox{$\mathfrak C$} \,\|f \|_{q,s;T}\,
\Bigl( \int_{ T+\mu} ^{ \infty } \int_{ \partial \Omega } \bigl[\, (t-T) ^{-1/2} \,
\bigl(\,(t+T)/2 -T \,\bigr) ^{-3/(2\, q)-1/s+1}\,\bigr]^2\, do_x\, dt \Bigr) ^{1/2}
\\&&\nonumber
\le
\mbox{$\mathfrak C$} \,\|f \|_{q,s;T}\,
\Bigl( \int_{ T+\mu} ^{ \infty } 
(t -T) ^{-3/ q-2/s+1}\, dt \Bigr) ^{1/2}
\le
\mbox{$\mathfrak C$} \,\|f \|_{q,s;T}\,
\mu ^{-3/(2\, q)-1/s+1},
\end{eqnarray}
where the last inequality holds because $-3/q-2/s+1<-1$ due to our
assumptions on $q$ and $s$. Again by referring to Lemma \ref{lemmaL6.10}, we get for
$t \in (T+\mu ,\, \infty ),\; x \in \partial \Omega $ that the term
$| \int_{ (t+T)/2}^t (t-r) ^{-1/2} \, \partial _rv(x,r)\, dr |$ is bounded by
\begin{eqnarray*}
\mbox{$\mathfrak C$} \,\|f\|_{q,s;T}\,
\int_{ (t+T)/2}^t (t-r) ^{-1/2} \, \sum_{j \in \{0,\, 1\} }(t-T)^{-3/(2\, q)-1/s+j/2}\, dr
,
\end{eqnarray*}
and thus by $\mbox{$\mathfrak C$} \,\|f\|_{q,s;T}\, \sum_{j \in \{0,\, 1\} }(t-T)^{-3/(2\, q)-1/s+1/2+j/2}.$
As a consequence, since $-3/ q-2/s+2<-1$ in view of our assumptions on $q$ and $s$,
\begin{eqnarray} \label{L6.30.20}
&&
\Bigl( \int_{ T+\mu}^{ \infty } \int_{ \partial \Omega }
\Bigl| \int_{ (t+T)/2}^t (t-r) ^{-1/2} \, \partial _r\Psi (x,r)\, dr \Bigr|^2
\,do_x\, dt \Bigr) ^{1/2}
\\&&\nonumber
\le
\mbox{$\mathfrak C$} \,\|f\|_{q,s;T}\,
\sum_{j \in \{0,\, 1\} }\mu ^{-3/(2 q)-1/s+1+j/2}.
\end{eqnarray}
Moreover,
\begin{eqnarray} \label{L6.30.30}
\Bigl( \int_{ T+\mu}^{ \infty } \int_{ \partial \Omega }
\Bigl| \int^{ (t+T)/2}_0 (t-r) ^{-3/2} \, v(x,r)\, dr \Bigr|^2\,do_x\,dt \Bigr) ^{1/2} 
\le
\mbox{$\mathfrak A$} ^{(1)} +\mbox{$\mathfrak A$} ^{(2)} ,
\end{eqnarray}
where $\mbox{$\mathfrak A$} ^{(j)} $ for $j \in \{1,\, 2\} $
is defined as the left-hand side of (\ref{L6.30.30})
with $\Psi (x,r)$ replaced by
\begin{eqnarray*}
\mbox{$\mathfrak R$} ^{(1)} (x,r):= \int_{ 0}^{\min\{r,T\}}\int_{ \mathbb{R}^3 }\chi_{(0,1]}(r-\sigma )\,
\Lambda (x-y,r-\sigma , \tau ) \cdot f(y, \sigma )\, dy\, d \sigma
\end{eqnarray*}
$(x \in \partial \Omega ,\; r \in (0, \infty ) )$
in the case $j=1$, and by a term $\mbox{$\mathfrak R$} ^{(2)} (x,r)$ differing from
$\mbox{$\mathfrak R$} ^{(1)} (x,r)$ insofar as 
$\chi_{(1, \infty )}(r-\sigma )$ substitutes for $\chi_{(0,1]}(r-\sigma )$  
in the case $j=2$. In order to estimate $\mbox{$\mathfrak A$} ^{(1)} $, we
exploit the integration over $\partial \Omega $ (Lemma \ref{lemmaL6.20}) in order to
reduce the singularity of $\mbox{$\mathfrak R$} ^{(1)} (x,r)$ when $r$ tends to zero.
In a first step, we use Theorem \ref{theoremT1.10} to obtain
$
\mbox{$\mathfrak A$} ^{(1)}
\le
\bigl(\,  \int_{ T+\mu}^{ \infty } \mbox{$\mathfrak B$} (t)^2\, dt \,\bigr) ^{1/2} ,
$
with
\begin{eqnarray*} 
\mbox{$\mathfrak B$} (t):=
\int_{ 0}^{(t+T)/2}(t-r) ^{-3/2} \,
\Bigl( \int_{ \partial \Omega }| \mbox{$\mathfrak R$} ^{(1)} (x,r)|^2\, do_x \Bigr) ^{1/2}
\, dr
\end{eqnarray*}
Let $d$ be any number from $(0,3)$, for example $d=3/2$.
Theorem \ref{theoremT1.10} applied once more yields
that
$\bigl(\,  \int_{ \partial \Omega }| \mbox{$\mathfrak R$} ^{(1)} (x,r)|^2\, do_x \,\bigr)  ^{1/2}$
is bounded by
\begin{eqnarray*}
\int_{ 0}^{\min\{r,T\}}\chi_{(0,1]}(r-\sigma )\, \Bigl( \int_{ \partial \Omega } \Bigl[
\int_{ \mathbb{R}^3 }
|\Lambda (x-y,t-\sigma , \tau ) \cdot f(y, \sigma )|\, dy \Bigr]^2\, do_x \Bigr) ^{1/2} \, d \sigma
,
\end{eqnarray*}
and hence by
$
\mbox{$\mathfrak C$} \, \int_{ 0}^{\min\{r,T\}}  \chi_{(0,1]}(r-\sigma )\,\|f( \sigma )\|_q\,
(r-\sigma )^{-d/(2q ^{\prime} )-3/(2q)+1/2}\, d \sigma 
$
according to Lemma \ref{lemmaL6.20}.
This estimate and H\"older's inequality imply for $t \in (T+\mu ,\, \infty )$ that
\begin{eqnarray*}&&
\mbox{$\mathfrak B$} (t)
\le
\Bigl( \int_{ 0}^{(t+T)/2}(t-r) ^{-3 s ^{\prime} /2}\, dr \Bigr) ^{1/s ^{\prime} }
\\&&\hspace{2.5em}
\cdot 
\Bigl(  \int_{ 0}^{(t+T)/2} \Bigl[ \int_{ 0}^{\min\{r,T\}}\chi_{(0,1]}(r-\sigma )\,
(r-\sigma )^{-d/(2q ^{\prime} )-3/(2q)+1/2}\,\|f( \sigma )\|_q\, d \sigma \Bigr]^s\, dr \Bigr) ^{1/s}
,
\end{eqnarray*}
hence
with Young's inequality,
\begin{eqnarray*}
\mbox{$\mathfrak B$} (t)
\le
\mbox{$\mathfrak C$} \,\|f\|_{q,s;T}\,
(t-T)^{-3/2+1/s ^{\prime} }\,\int_{ \mathbb{R} } \chi_{(0,1]}( \sigma )\,
\sigma ^{-d/(2 q ^{\prime} )-3/(2q)+1/2}\,d \sigma 
.
\end{eqnarray*}
But $-d/(2\, q ^{\prime} )-3/(2q)+1/2>-1$
because $d<3$ and $q>1$, so
\begin{eqnarray*}
\mbox{$\mathfrak B$} (t)
\le
\mbox{$\mathfrak C$} \,\|f\|_{q,s;T}\,(t-T)^{-3/2+1/s ^{\prime} }
\le
\mbox{$\mathfrak C$} \,\|f\|_{q,s;T}\,(t-T)^{-1/2-1/s },
\end{eqnarray*}
for $t$ as before. Recalling that
$
\mbox{$\mathfrak A$} ^{(1)}
\le
\bigl(\,  \int_{ T+\mu}^{ \infty } \mbox{$\mathfrak B$} (t)^2\, dt \,\bigr) ^{1/2} ,
$
we now obtain
\begin{eqnarray} \label{L6.30.50}
\mbox{$\mathfrak A$}^{(1)} 
\le
\mbox{$\mathfrak C$} \,\|f\|_{q,s;T}\, \Bigl( \int_{ T+\mu}(t-T)^{-1-2/s }\, dt \Bigr) ^{1/2}
\le
\mbox{$\mathfrak C$} \,\|f\|_{q,s;T}\,\mu^{-1/s}.
\end{eqnarray}
Turning to $\mbox{$\mathfrak A$} ^{(2)} $, we deduce from H\"older's inequality that
\begin{eqnarray} \label{L6.30.60}&&\hspace{-3em}
\int_{ 0}^{(t+T)/2}(t-r) ^{-3/2} \,| \mbox{$\mathfrak R$} ^{(2)} (x,r)|\, dr
\le
\Bigl( \int_{ 0}^{(t+T)/2}(t-r)^{-3s ^{\prime} /2}\, dr \Bigr) ^{1/s ^{\prime} }\,\| \mbox{$\mathfrak R$} ^{(2)}
(x,\, \cdot \,)\|_s
\\&&\nonumber \hspace{-3em}
\le
\mbox{$\mathfrak C$} \, (t-T)^{-3/2+1/s ^{\prime} }\,\| \mbox{$\mathfrak R$} ^{(2)} (x,\, \cdot \, )\|_s
\le
\mbox{$\mathfrak C$} \, (t-T)^{-1/2-1/s }\,\| \mbox{$\mathfrak R$} ^{(2)} (x,\, \cdot \, )\|_s
\end{eqnarray}
for
$x \in \partial \Omega ,\; t \in (T+\mu ,\, \infty ).$
But $q<3/2$, hence $1-3/(2q)< 0$, so we get by Theorem \ref{theoremT4.30} with $\alpha =0,\; M=1,\;
\varrho =s$ that
$
\|\mbox{$\mathfrak R$} ^{(2)} (x,\, \cdot \, )\|_s\le C(q,s)\|f\|_{q,s;T}
$
for any $x \in \mathbb{R}^3 $.
This choice of $\varrho $ is possible since $s>1$.
Thus from (\ref{L6.30.60}),
$
\mbox{$\mathfrak A$} ^{(2)}
\le
\mbox{$\mathfrak C$} \, \bigl(\,  \int_{ T+\mu}^{ \infty } (t-r)^{-1-2/s}\, dr
\,\bigr) ^{1/2} \, \|f\|_{q,s;T}
\le
\mbox{$\mathfrak C$} \, \|f\|_{q,s;T}\,\mu^{-1/s}.
$
The lemma follows from the preceding inequality, (\ref{L6.30.10}) -- (\ref{L6.30.30}), (\ref{L6.30.50}),
and from equation (\ref{L10.10.10}) as indicated at the beginning of this proof.
\hfill $\Box $
\begin{corollary} \label{corollaryC6.20}
Let $q \in (1,\, 3/2),\; \hat{q}\in [1,\,3/2),\, s \in (1, \infty )$ with
$
3/(2\, q)+1/s >3/2,\;
f \in L^c \bigl(\, 0, \infty ,\; L^d( \mathbb{R}^3 )^3 \,\bigr)
$
for
$
(c,d) \in \{(s,q),\, (2,\hat{q}),\, (2,\,3/2)\}.
$
Let $T \in [1, \infty ).$ Then
\begin{eqnarray*}&&
\| \mbox{$\mathfrak R$} ^{( \tau )} (f)|S_{T, \infty }\|_{H_{T, \infty }}
\le
\mbox{$\mathfrak C$} \,\|f\|_{q,s; \infty }
\, (T^{-3/(2q)-1/s+3/2}+T^{-1/s})
\\&&\hspace{3em}
+
\mbox{$\mathfrak C$} \,
\bigl(\, \|f| \mathbb{R}^3 \times (T/2,\, \infty )\|_{\hat{q},2; \infty }
+ \|f| \mathbb{R}^3 \times (T/2,\, \infty )\|_{3/2,\, 2; \infty }\,\bigr) .
\end{eqnarray*}
(Recall that by Theorem \ref{theoremT10.5}, we have $\mbox{$\mathfrak R$} ^{( \tau )} (f)|S_{\infty }
\in H_{\infty }$.)
\end{corollary}
{\bf Proof:}
Put $f ^{(1)} :=\chi_{\mathbb{R}^3 \times (T/2,\, \infty )}\, f,
\;
f ^{(2)} :=\chi_{\mathbb{R}^3 \times (0,\, T/2)}\, f.
$
Then by Corollary \ref{corollaryC10.40} and Theorem \ref{theoremT10.5},
\begin{eqnarray} \label{C6.20.20}
\| \mbox{$\mathfrak R$} ^{( \tau )} (f ^{(1)} )|S_{T, \infty }\|_{H_{T, \infty }}
\le
\| \mbox{$\mathfrak R$} ^{( \tau )} (f ^{(1)} )\|_{H_{\infty }}
\le
\bigl(\, \|f ^{(1)} \|_{\hat{q},2; \infty }
+ \|f ^{(1)} \|_{3/2,\, 2; \infty }\,\bigr) .
\end{eqnarray}
Concerning $\mbox{$\mathfrak R$} ^{( \tau )} (f ^{(2)} ),$ we use equation (\ref{C10.65.10}) with
$f ^{(2)} $ in the place of $f$, and then Lemma \ref{lemmaL6.15} and \ref{lemmaL6.30} with $T$
replaced by $T/2$ and $\mu=T/2$. Taking account of the assumption $T\ge 1$, we find that
\begin{eqnarray} \label{C6.20.30}
\| \mbox{$\mathfrak R$} ^{( \tau )} (f ^{(2)} )|S_{T, \infty }\|_{H_{T, \infty }}
\le
\mbox{$\mathfrak C$} \,\|f\|_{q,s; \infty }
\, (T^{-3/(2q)-1/s+3/2}+T^{-1/s}).
\end{eqnarray}
Corollary \ref{corollaryC6.20} follows from (\ref{C6.20.20}) and (\ref{C6.20.30}).
\hfill $\Box $

\vspace{1ex}
In section 8, we will need a pointwise estimate on spatial and temporal decay of $\mbox{$\mathfrak R$}
^{( \tau )} (f)$, based on Theorem \ref{theoremT4.10} and Lemma \ref{lemmaL6.10}. Contrary to the
situation in Theorem \ref{theoremT4.10}, we want to avoid the assumption
$f|B_{R_0}\times (0, \infty ) \in L^2 \bigl(\, B_{R_0}\times (0, \infty ) \,\bigr) ^3$
because it is inconvenient in the nonlinear case. The ensuing lemma allows us to replace this
condition by the weaker one
$f|B_{R_0}\times (0, \infty ) \in L^2 \bigl(\, 0, \infty ,\, L^1(B_{R_0})^3\,\bigr) $.
We will apply this lemma only for $\varrho =2$, but still consider the case $\varrho \in [1, \infty )$
because one may hope that in future work,
the decay rates obtained in the special situation of Lemma \ref{lemmaL6.40} may be recovered
in a more general setting.
\begin{lemma} \label{lemmaL6.40}
Let
$\varrho \in [1, \infty ),\; R_0 \in (0, \infty ) $
and
$f \in L ^{\varrho } \bigl(\, 0, \infty ,\, L^1(B_{R_0})^3\,\bigr) $.
Take $R \in (R_0, \infty ).$
Then, for $\alpha \in \mathbb{N} _0^3$ with $ | \alpha |\le 1,\; \epsilon \in [0,1],\, x \in B_R^c$
and $t \in (0, \infty ) $,
\begin{eqnarray*}&&
| \partial _x^{\alpha }\mbox{$\mathfrak R$} ^{( \tau )} (f)(x,t)|
\\&&
\le
\mbox{$\mathfrak C$} \,
\bigl [\,
\bigl(\,   |x|\, \nu (x) \,\bigr)^{-1/2-1/\varrho -| \alpha |/2}\,
\bigl(\, (1+t) ^{-2} \, \|f\|_{1, \varrho ; \infty }
+
\|f| B_{R_0}\times (t/2,\,  \infty )\|_{1, \varrho ; \infty })
\\&&\hspace{3em}
+
\bigl(\,   |x|\, \nu (x) \,\bigr)^{(-1/2 -1/\varrho -| \alpha |/2)\, (1-\epsilon )} \,
(1+t)^{(-1/2-1/\varrho -| \alpha |/2)\, \epsilon }\,\|f\|_{1, \varrho ; \infty }
\,\bigr] .
\end{eqnarray*}
\end{lemma}
{\bf Proof:}
Take $\alpha ,\, \epsilon ,\,x,\, t$ as above, and put
$f ^{(1)} :=\chi_{B_{R_0} \times (0,\,t/2)}\, f,
\;
f ^{(2)} :=\chi_{B_{R_0} \times (t/2,\, \infty )}\, f.
$
For $y \in B_{R_0}$, we get
\begin{eqnarray} \label{L6.40.10}
|x-y|\ge |x|-|y|= |x|\, (1-R_0/R) +|x|\, R_0/R-|y|
\ge |x|(1-R/R_0)\ge R-R_0,
\end{eqnarray}
where $R-R_0>0$.
Thus we may apply Corollary \ref{corollaryC10.60} with $K=R-R_0$.
Suppose that $\varrho >1$. If $\varrho =1$, the ensuing argument
remains valid with some small modifications.
Due Lemma \ref{lemmaL10.85},
Corollary \ref{corollaryC10.60} used as indicated,
(\ref{L6.40.10}) and Lemma \ref{lemmaL10.1}, we get
\begin{eqnarray} \label{L6.40.20}&&
| \partial _x^{\alpha }\mbox{$\mathfrak R$} ^{( \tau )} (f ^{(1)} )(x,t)|
\\&&\nonumber 
\le
\mbox{$\mathfrak C$} \,
\int_{ 0}^{t/2}\int_{ B_{R_0}}\bigl(\, |x-y|\, \nu (x-y) +t-\sigma \,\bigr) ^{-3/2-| \alpha |/2}\,
|f(y, \sigma )|\, dy\, d \sigma
\\&&\nonumber 
\le
\mbox{$\mathfrak C$} \,
\int_{ 0}^{t/2}\int_{ B_{R_0}}\bigl(\, |x|\, \nu (x) +t-\sigma \,\bigr) ^{-3/2-| \alpha |/2}\,
|f(y, \sigma )|\, dy\, d \sigma
\\&&\nonumber 
\le
\mbox{$\mathfrak C$} \,
\Bigl( \int_{ 0}^{t/2}\bigl(\, |x|\, \nu (x) +t-\sigma \,\bigr) ^{(-3/2-| \alpha |/2)\, \varrho ^{\prime} }
\,d \sigma \Bigr) ^{1/ \varrho ^{\prime} } 
\|f\|_{1, \varrho ; \infty }
\\&&\nonumber 
\le
\mbox{$\mathfrak C$} \,
\bigl(\, |x|\, \nu (x) +t \,\bigr) ^{-3/2 -| \alpha |/2+1/ \varrho ^{\prime}} \,\|f\|_{1, \varrho ; \infty }.
\end{eqnarray}
In the case $t\le 1$, we thus get
$
| \partial _x^{\alpha }\mbox{$\mathfrak R$} ^{( \tau )} (f ^{(1)} ) (x,t)|
\le
\mbox{$\mathfrak C$} \,
\bigl(\, |x|\, \nu (x)\,\bigr) ^{-3/2+1/ \varrho ^{\prime} -| \alpha |/2}\,\|f\|_{1, \varrho ; \infty }
$,
and then multiply the right-hand side by
$4\, (1+t) ^{-2}$.
If $t\ge 1$, it follows from (\ref{L6.40.20}) that
$
| \partial _x^{\alpha }\mbox{$\mathfrak R$} ^{( \tau )} (f ^{(1)} )(x,t)|
\le
\mbox{$\mathfrak C$} \,
\bigl(\, |x|\, \nu (x)\,\bigr) ^{(-3/2+1/ \varrho ^{\prime} -| \alpha |/2)\,(1-\epsilon )}\,
(1+t)^{(-3/2+1/ \varrho ^{\prime} -| \alpha |/2)\,\epsilon }\,
\|f\|_{1, \varrho ; \infty }.
$
As in (\ref{L6.40.20}), we obtain that
$
| \partial _x^{\alpha }\mbox{$\mathfrak R$} ^{( \tau )} (f ^{(2)} )(x,t)|
$
is bounded by
\begin{eqnarray*} 
\mbox{$\mathfrak C$} \,
\Bigl( \int_{ t/2}^{t}\bigl(\, |x|\, \nu (x) +t-\sigma \,\bigr) ^{(-3/2-| \alpha |/2)\, \varrho ^{\prime} }\,
d \sigma \Bigr) ^{1/ \varrho ^{\prime} }\,
\|f|B_{R_0}\times (t/2,\, \infty )\|_{1, \varrho ; \infty },
\end{eqnarray*}
and hence by
$
\mbox{$\mathfrak C$} \,
\bigl(\, |x|\, \nu (x)\,\bigr) ^{-3/2+1/ \varrho ^{\prime} -| \alpha |/2}\,
\|f|B_{R_0}\times (t/2,\, \infty )\|_{1, \varrho ; \infty }.
$
The lemma follows from the preceding estimates and because $-3/2+1/ \varrho ^{\prime} =-1/2-1/ \varrho .$
\hfill $\Box $
\begin{lemma} \label{lemmaL6.50}
Let $\alpha \in \mathbb{N} _0^3$ with
$| \alpha |\le 1,\; q \in [1,4],\; \widetilde{ q},\, \overline{ q},\,
s \in [1, \infty ),\; \widetilde{ s},\, \overline{ s} \in [1, \infty ]$
with
$
3/(2\, d)+1/c>1-| \alpha |/2
$
for $(c,d) \in \{(s,q),\, ( \widetilde{ s}, \widetilde{ q})\}$,
and
$
3/(2\, \overline{ q})+1/ \overline{ s}<1-| \alpha |/2.
$
Suppose that $R_0 \in (0, \infty ) $ and
$f \in L^c \bigl(\, 0, \infty ,\, L^q( B_{R_0}^c)^3 \,\bigr) $
for $(c,d) \in \{(s,q),\, ( \widetilde{ s}, \widetilde{ q}), ( \overline{  s}, \overline{  q})\}.$
Then, for $x \in \mathbb{R}^3 ,\; t \in [1, \infty ),$
\begin{eqnarray*}&&
| \partial _x^{\alpha }\mbox{$\mathfrak R$} ^{( \tau )} (f )(x,t)|
\le
\mbox{$\mathfrak C$} \, \bigl(\,
\|f\|_{q,s; \infty }\, (1+t) ^{-3/(2q)-1/s+1-| \alpha |/2}
\\&&\hspace{3em}
+
\|f|B_{R_0}^c\times (t/2,\, \infty )\|_{\widetilde{ q}, \widetilde{ s}; \infty }
+
\|f|B_{R_0}^c\times (t/2,\, \infty )\|_{\overline{  q}, \overline{  s}; \infty }\,\bigr) .
\end{eqnarray*}
\end{lemma}
{\bf Proof:}
Let $x,\, t$ as in the lemma, and put
$f ^{(1)} :=\chi_{B_{R_0}^c \times (0,\,t/2)}\, f,
\;
f ^{(2)} :=\chi_{B_{R_0}^c \times (t/2,\, \infty )}\, f.
$
Lemma \ref{lemmaL6.10} with $f$ replaced by $f ^{(1)} $ and with $T=t/2$, and the assumptions
$t\ge 1,\; 3/(2\, q)+1/s>1-| \alpha |/2$ imply
$
| \partial _x^{\alpha }\mbox{$\mathfrak R$} ^{( \tau )} (f ^{(1)} )(x,t)|
\le
\mbox{$\mathfrak C$} \,\|f\|_{q,s; \infty }\, (1+t)^{-3/(2q)-1/s+1-| \alpha |/2}.
$
We further note that according to Lemma \ref{lemmaL10.70}, we have
$
| \partial _x^{\alpha }\mbox{$\mathfrak R$} ^{( \tau )} (f ^{(2)} )(x,t)|
\le
\mbox{$\mathfrak A$} _1+\mbox{$\mathfrak A$} _2
$,
with
$
\mbox{$\mathfrak A$} _1:=
\int_{ t/2}^t \int_{ B_{R_0}^c}\chi_{(0,1]}(t-\sigma )\,| \partial _x^{\alpha }
\Lambda (x-y,t-\sigma , \tau ) \cdot f(y, \sigma )|\,
dy\, d \sigma ,
$
and $\mbox{$\mathfrak A$} _2$ defined as $\mbox{$\mathfrak A$} _1, $ but with the term
$\chi_{(0,1]}(t-\sigma )$ replaced by $\chi_{(1, \infty )}(t-\sigma )$.
Since
$
3/(2\, \overline{ q})+1/ \overline{ s}<1-| \alpha |/2,
$
Theorem \ref{theoremT4.30} with $\varrho =\infty ,\; M=1$ yields
$
\mbox{$\mathfrak A$} _1\le \mbox{$\mathfrak C$} \,
\|f|B_{R_0}^c\times (t/2,\, \infty )\|_{\overline{  q}, \overline{  s}; \infty }.
$
Similarly,
the assumption
$
3/(2\, \widetilde{  q})+1/ \widetilde{  s}>1-| \alpha |/2
$
and
Theorem \ref{theoremT4.30}, again with $\varrho =\infty ,\; M=1$, imply
$
\mbox{$\mathfrak A$} _2\le \mbox{$\mathfrak C$} \,
\|f|B_{R_0}^c\times (t/2,\, \infty )\|_{\widetilde{   q}, \widetilde{   s}; \infty }.
$
Lemma \ref{lemmaL6.50} follows by combining the preceding estimates.
\hfill $\Box $

\vspace{1ex}
Now we collect what we found in this section about temporal pointwise decay
of $\mbox{$\mathfrak R$} ^{( \tau )} (f)$, and combine it with a pointwise spatial decay estimate.
\begin{corollary} \label{corollaryC6.30}
Let $\alpha \in \mathbb{N} _0^3$ with
$| \alpha |\le 1,\; \varrho \in [1, \infty ),\; q \in [1,4],\; \widetilde{ q},\, \overline{ q},\,
s \in [1, \infty ),\; \widetilde{ s},\, \overline{ s} \in [1, \infty ]$
with
$
3/(2\, d)+1/c>1-| \alpha |/2
$
for $(c,d) \in \{(s,q),\, ( \widetilde{ s}, \widetilde{ q})\}$,
and
$
3/(2\, \overline{ q})+1/ \overline{ s}<1-| \alpha |/2.
$
Suppose that $R_0 \in (0, \infty ) $ and
the function $f: \mathbb{R}^3 \times (0, \infty )\mapsto \mathbb{R}^3$
satisfies
$
f|B_{R_0} \times (0, \infty ) \in L^{\varrho } \bigl(\, 0, \infty ,\, L^1( B_{R_0})^3 \,\bigr),\;
f|B_{R_0}^c \times (0, \infty ) \in L^c \bigl(\, 0, \infty ,\, L^q( B_{R_0}^c)^3 \,\bigr) $
with $(c,d)$ being anyone of the pairs
$(s,q),\, ( \widetilde{ s}, \widetilde{ q})$ and $ ( \overline{  s}, \overline{  q})$.
Let $R \in (R_0, \infty )$ and further suppose there is $\mbox{$\mathfrak D$} _0>0$ such that
\begin{eqnarray} \label{C6.30.10}
| \partial _x^{\alpha }\mbox{$\mathfrak R$} ^{( \tau )} \bigl(\, f | B_{R_0}^c \times (0, \infty ) \,\bigr) (x,t)|
\le
\mbox{$\mathfrak D$} _0\, \bigl(\, |x|\,\nu(x) \,\bigr) ^{-1/2-1/ \varrho -| \alpha |/2}
\; \mbox{for}\;
x \in B_R^c,\; t>0.
\end{eqnarray}
(This condition is satisfied if, for example, $\varrho =2$ and $\chi_{B_{R_0}^c \times (0, \infty )}\, f$ fulfills the
assumptions of Theorem \ref{theoremT4.10} in the place of $f$.)
Then, for $x \in B_R^c,\; t  \in (0, \infty ) ,\; \epsilon \in [0,1],$
\begin{eqnarray} \label{C6.30.20}&&\hspace{-3em}
| \partial _x^{\alpha }\mbox{$\mathfrak R$} ^{( \tau )} (f )(x,t)|
\le
\mbox{$\mathfrak C$} \, \bigl[\, 
\bigl(\, |x|\,\nu(x) \,\bigr) ^{-1/2-1/ \varrho -| \alpha |/2}\,
\bigl(\, ( \mbox{$\mathfrak D$} _0+K_1)\, (1+t) ^{-2} +K_2(t) \,\bigr)
\\&&\nonumber 
+
( \mbox{$\mathfrak D$} _0^{1-\epsilon }+K_1)\, \bigl(\, |x|\,\nu(x) \,\bigr) ^{(-1/2-1/ \varrho -| \alpha |/2)\,(1-\epsilon )}
\\&&\nonumber \hspace{4em}
\bigl(\, (1+t)^{-1/2-1/ \varrho -| \alpha |/2} + K_3\, (1+t)^{-3/(2q)-1/s+1-| \alpha |/2} +K_4(t) \,\bigr) ^{\epsilon }
\,\bigr],
\end{eqnarray}
with
$
K_1:=\|f|B_{R_0}\times (0, \infty )\|_{1, \varrho ; \infty },\;
K_2(t):=\|f|B_{R_0}\times (t/2,\, \infty )\|_{1, \varrho ; \infty },\;
K_3:=\|f|B_{R_0}^c\times (0, \infty )\|_{q,s; \infty },\;
K_4(t):=\|f|B_{R_0}^c\times (t/2,\, \infty )\|_{\widetilde{ q}, \widetilde{ s}; \infty }
+
\|f|B_{R_0}^c\times (t/2,\, \infty )\|_{\overline{  q}, \overline{  s}; \infty }.
$
\end{corollary}
{\bf Proof:} Take $x,\,t,\, \epsilon $ as in the lemma.
We apply Lemma \ref{lemmaL6.40} with $f$
replaced by $f ^{(1)} :=f|B_{R_0}\times (0, \infty )$. It follows that
$
| \partial _x^{\alpha }\mbox{$\mathfrak R$} ^{( \tau )} (f ^{(1)}  )(x,t)|
$
is bounded by the right-hand side of (\ref{C6.30.20}).
Suppose that $t\ge 1$.
Starting from the equation
$d=d^{\epsilon} \,d^{1-\epsilon}$
for $d>0$,
we use (\ref{C6.30.10}) as well as Lemma \ref{lemmaL6.50} with $f$ replaced by
$f ^{(2)} :=\chi_{B_{R_0}^c\times (0, \infty )}\, f$.
It follows that 
$
| \partial _x^{\alpha }\mbox{$\mathfrak R$} ^{( \tau )} (f ^{(2)}  )(x,t)|
$
is also bounded by the right-hand side of (\ref{C6.30.20}).
In the case $t\le 1$, inequality (\ref{C6.30.10}) multiplied by $4\,(1+t) ^{-2} $
yields a suitable estimate of
$
| \partial _x^{\alpha }\mbox{$\mathfrak R$} ^{( \tau )} (f ^{(2)}   )(x,t)|
.
$
\hfill $\Box $

\section{Temporal decay of the potential $\mbox{$\mathfrak I$} ^{( \tau )} (a)$.}
\setcounter{equation}{0}

First we consider the behaviour of $\| \mbox{$\mathfrak I$} ^{( \tau )} (a)|S_{T, \infty }\|_{H_{T, \infty }}$
when $T $ tends to infinity.
\begin{theorem} \label{theoremT7.10}
Let $\epsilon _0 \in (0,\, 1/2],\; p \in (1,2],\; a \in H_{\sigma }^{1/2+ \epsilon _0}( \overline{ \Omega }^c)
\cap L^p( \overline{ \Omega }^c)^3$ and $T \in [1, \infty )$.
Then 
$
\| \mbox{$\mathfrak I$} ^{( \tau )} (a)|S_{T, \infty }\|_{H_{T, \infty }}
\le
\mbox{$\mathfrak C$} \,\|a\|_p\, T ^{-3/(2p)+1/2}
$.
(Note that
$\mbox{$\mathfrak I$} ^{( \tau )} (a)|S_{\infty } \in H_{\infty }$
due to the assumption
$a \in H_{\sigma }^{1/2+ \epsilon _0}( \overline{ \Omega }^c)$; see Theorem \ref{theoremT10.6}.)
\end{theorem}
{\bf Proof:}
Recall that $\mbox{$\mathfrak I$} ^{( \tau )} (a) \in C^1 \bigl(\, \mathbb{R}^3 \times (0, \infty ) \,\bigr) ^3$
according to Lemma \ref{lemmaL10.100}. From (\ref{L10.100.10}) and Lemma \ref{lemmaL10.90} we
get for $\alpha \in \mathbb{N} _0^3$ and $l \in \mathbb{N} _0$ with $| \alpha |+l\le 1$
that
\begin{eqnarray*} \label{T7.10.10}&&\hspace{-2em}
\Bigl( \int_{ T} ^{ \infty } \int_{ \partial \Omega }| \partial _t^l \partial _x^{\alpha }
\mbox{$\mathfrak I$} ^{( \tau )} (a)(x,t)|^2\, do_x\,dt \Bigr) ^{1/2}
\le
\mbox{$\mathfrak C$} \,
\|a\|_p\,
\Bigl( \int_{ T}^{ \infty } (t^{-3/p-| \alpha |-2l}+\delta _{1l}\, t^{-3/p-1})\, dt \Bigr) ^{1/2}
\\&&\nonumber \hspace{-2em}
\le
\mbox{$\mathfrak C$} \,
\|a\|_p\,
(T^{-3/(2p)-| \alpha |/2-l+1/2}+\delta _{1l}\, T^{-3/(2p)})
.
\end{eqnarray*}
Since $T\ge 1$,
the right-hand side of this estimate is bounded by $\mbox{$\mathfrak C$} \,
\|a\|_p\,T^{-3/(2p)+1/2}$,
so
\begin{eqnarray} \label{T7.10.20}
\int_{ T} ^{ \infty }\bigl(\, \|\mbox{$\mathfrak I$} ^{( \tau )} (a)(t)| \partial \Omega
\|_{H^1( \partial \Omega )^3}^2 +\| n^{ (\Omega )} \cdot \partial _4 \mbox{$\mathfrak I$} ^{( \tau )} (a)(t)
\|_2^2\,\bigr) \, dt 
\le
\mbox{$\mathfrak C$} \, \|a\|_p^2 \,T^{-3/p+1}.
\end{eqnarray}
Turning to an estimate of $\partial _4 ^{1/2} \bigl(\, \mbox{$\mathfrak I$} ^{( \tau )} (a)|S_{\infty }\,\bigr) ,$
we note that by Corollary \ref{corollaryC10.66}, this fractional derivative exists, and
equation (\ref{L10.10.10}) holds with $\mbox{$\mathfrak I$} ^{( \tau )} (a)(x,\, \cdot \, )$
in the role of $\phi $ and with $T$ replaced by $T/2$,
for a. e. $x \in \partial \Omega $.
With this in mind, we find with
Lemma \ref{lemmaL10.90} that
\begin{eqnarray}\label{T7.10.30}&&
\Bigl( \int_{ T} ^{ \infty } \int_{ \partial \Omega }| (t-T/2) ^{-1/2} \,
\mbox{$\mathfrak I$} ^{( \tau )} (a) \bigl(\, x,\, (t+T/2)/2 \,\bigr) | ^2\, do_x\, dt \Bigr)  ^{1/2}
\\&&\nonumber 
\le
C(p)\, \|a\|_p\, \Bigl( \int_{ T} ^{ \infty }  \bigl(\, t ^{-1/2} \, (t+T/2)^{-3/(2p)} \,\bigr)  ^2\,  dt \Bigr)  ^{1/2}
\le
C(p)\, \|a\|_p\,  T ^{-3/(2p)}.
\end{eqnarray}
We further find for $t \in (T, \infty ),\; x \in \partial \Omega $ that
by (\ref{L10.100.10}) and Lemma \ref{lemmaL10.90},
\begin{eqnarray*} &&\hspace{-2em}
\int_{ (t+T/2)/2} ^{ t } (t-r) ^{-1/2} \,| \partial _4\mbox{$\mathfrak I$} ^{( \tau )} (a)(x,r)|\, dr
\\&&\hspace{-2em}
\le
C(p) \, \|a\|_p\, \int_{(t+T/2)/2}^t (t-r) ^{-1/2} \,\sum_{j \in \{0,\, 1\}  } r^{-3/(2p)-1/2-j/2}\, dr
\\&&\hspace{-2em}
\le
C(p) \, \|a\|_p\, \sum_{j \in \{0,\, 1\}  } (t+T/2)^{-3/(2p)-1/2-j/2}\,\int_{(t+T/2)/2}^t (t-r) ^{-1/2} \, dr
\le
C(p) \, \|a\|_p\,  t^{-3/(2p)},
\end{eqnarray*}
where the last inequality is valid because $t\ge T\ge 1$. As a consequence,
since $p<3,$ hence $-3/p< -1$,
\begin{eqnarray} \label{T7.10.40}
\Bigl( \int_{ T} ^{ \infty } \int_{ \partial \Omega }\Bigl| 
\int_{ (t+T/2)/2} ^{ t } (t-r) ^{-1/2}  \partial _4\mbox{$\mathfrak I$} ^{( \tau )} (a)(x,r)\, dr \Bigr| 
do_x\, dt \Bigr) ^{1/2}
\le
C(p)  \|a\|_p  T^{-3/(2p)+1/2}.
\end{eqnarray}
Next we consider the term
\begin{eqnarray} \label{T7.10.50}
\mbox{$\mathfrak A$}
:=
\Bigl( \int_{ T} ^{ \infty } \int_{ \partial \Omega }\Bigl| 
\int^{ (t+T/2)/2} _{ 0 } (t-r) ^{-3/2}\,  \mbox{$\mathfrak I$} ^{( \tau )} (a)(x,r)\, dr \Bigr| \,
do_x\, dt \Bigr) ^{1/2}.
\end{eqnarray}
In the case $p>3/2$, we have $3/(2\,p)< 1$, so we get with Lemma \ref{lemmaL10.90} that
\begin{eqnarray} \label{T7.10.60}&&
\mbox{$\mathfrak A$}
\le
C(p)\,\|a\|_p\,
\Bigl( \int_{ T} ^{ \infty } \Bigl[
\int^{ (t+T/2)/2} _{ 0 } (t-r) ^{-3/2} \,r^{-3/(2p)}\, dr \Bigr]^2\, dt \Bigr) ^{1/2} 
\\&& \nonumber
\le
C(p)\,\|a\|_p\,
\Bigl( \int_{ T} ^{ \infty } t ^{-3} \,\Bigl[
\int^{ (t+T/2)/2} _{ 0 } r^{-3/(2p)}\, dr \Bigr]^2\, dt \Bigr) ^{1/2} 
\le
C(p)\,\|a\|_p\, T ^{-3/(2p)}.
\end{eqnarray}
Now suppose that $p\le 3/2$ so that $3/(2\, p)\ge 1.$ Then we use the splitting
$\mbox{$\mathfrak A$} \le \mbox{$\mathfrak A$} _1+\mbox{$\mathfrak A$} _2,$
where $\mbox{$\mathfrak A$} _1$ is defined in the same way as \mbox{$\mathfrak A$} (see (\ref{T7.10.50})),
but with the lower bound $0$ of the integral with respect to $r$ replaced by $3/4$.
Similarly, the definition of $\mbox{$\mathfrak A$} _2$ follows that of \mbox{$\mathfrak A$},
but the integral with respect to $r$ is to extend from $0$ to $3/4$.
If $p<3/2,$  we again use the first two inequalities in (\ref{T7.10.60}), with the lower bound of
the integral with respect to $r$ being $3/4$ instead of $0$.
Since $3/(2\, p)>1$, we have $\int_{ 3/4}^{(t+T/2)/2}t^{-3/(2p)}\, dr\le C(p)$, so we get
that
$
\mbox{$\mathfrak A$} _1
\le
C(p)\, \|a\|_p\, \bigl(\,  \int_{ T}^{ \infty } t ^{-3} \, dt \,\bigr)  ^{1/2} \le C(p)\, \|a\|_p\, T ^{-1} .
$

In the case $3/(2\, p)=1$, we observe that
$
\int_{ 3/4}^{(t+T/2)/2}t^{-3/(2p)}\, dr
\le
\int_{ 3/4}^{3t/4}r ^{-1} \, dr
=
\ln t.
$
Therefore, again starting as in (\ref{T7.10.60}),
we obtain
$
\mbox{$\mathfrak A$} _1\le C(p)\, \|a\|_p\, T ^{-1} \, (1+\ln T)
$
for $p=3/2$.
By what we have seen further above, the preceding estimate is valid in the case $p<3/(2\, p)$ as well,
and thus holds if $p\le 3/(2\, p)$. In order to deal with $\mbox{$\mathfrak A$} _2$, we reduce
the singularity of the variable $r$ in $\mbox{$\mathfrak I$} ^{( \tau )} (a)(x,r)$
for $r\downarrow 0$ by exploiting the integration on $\partial \Omega $. To this end, we use H\"older's
inequality and Theorem \ref{theoremT1.10}, as well as Lemma \ref{lemmaL10.40}, obtaining
with the abbreviation $F(x,y,r):=\mbox{$\mathfrak H$} (x-y- \tau \, r\, e_1,\, r)$ that
\begin{eqnarray*} &&\hspace{-2em}
\int_{ \partial \Omega }|\mbox{$\mathfrak I$} ^{( \tau )} (a)(x,r)|^2\, do_x 
\\&&\hspace{-2em}
\le
\int_{ \partial \Omega }\Bigl[ \Bigl( \int_{ \mathbb{R}^3 }F(x,y,r)
\, dy \Bigr) ^{1/p ^{\prime} }
\Bigl( \int_{ \mathbb{R}^3 }F(x,y,r)
\, |a(y)|^p\, dy \Bigr) ^{1/p}
\Bigr]^2\, do_x  
\\&&\hspace{-2em}
=
\int_{ \partial \Omega }\
\Bigl( \int_{ \mathbb{R}^3 }F(x,y,r)
\, |a(y)|^p\, dy \Bigr) ^{2/p}
\, do_x 
\le
\Bigl( \int_{ \mathbb{R}^3 }\Bigl[\int_{ \partial \Omega }F(x,y,r)^{2/p}
\, do_x \Bigr]^{p/2}\, |a(y)|^p\, dy \Bigr) ^{2/p}
\end{eqnarray*}
for $r \in (0, \infty ) $. Note that $2/p\ge 1$ because $p \in (1,2]$.
Now we apply Lemma \ref{lemmaL10.40} again to deduce from the preceding estimate that
$
\bigl(\,  \int_{ \partial \Omega }|\mbox{$\mathfrak I$} ^{( \tau )} (a)(x,r)|^2\, do_x \,\bigr)  ^{1/2}
$
is bounded by
\begin{eqnarray*}
C \,
\Bigl( \int_{ \mathbb{R}^3 }\Bigl[\int_{ \partial \Omega } (|x-y- \tau \, r\, e_1|+ r ^{1/2} )^{-6/p}
\, do_x \Bigr]^{p/2}\, |a(y)|^p\, dy \Bigr) ^{1/p},
\end{eqnarray*}
and hence by $\mbox{$\mathfrak C$} \, r^{-3/(2p)+1/2}\,\|a\|_p$ for $r \in (0, \infty ) $.
Therefore, once more by Theorem \ref{theoremT1.10},
\begin{eqnarray*} &&
\mbox{$\mathfrak A$} _2
\le
C \Bigl( \int_{ T} ^{ \infty } t ^{-3} \, \int_{ \partial \Omega }\Bigl[ \int_{ 0}^{3/4}| \mbox{$\mathfrak I$}
^{( \tau )} (a)(x,r)|\, dr \Bigr]^2\, do_x\, dt \Bigr) ^{1/2}
\\&&
\le
C \biggl( \int_{ T} ^{ \infty } t ^{-3} \Bigl[ \int_{ 0}^{3/4} \Bigl(
\int_{ \partial \Omega } | \mbox{$\mathfrak I$}^{( \tau )} (a)(x,r)|^2\, do_x \Bigr) ^{1/2} \,dr \Bigr]^2\, dt
\biggr)  ^{1/2}
\\&&
\le
\mbox{$\mathfrak C$} \,\|a\|_p\,
\Bigl(  \int_{ T} ^{ \infty } t ^{-3} \Bigl[ \int_{ 0}^{3/4} r^{-3/(2p)+1/2}\, dr \Bigr]^2\, dt \Bigr) ^{1/2} .
\end{eqnarray*}
Since $p>1$, we have $-3/(2\, p)+1/2>-1,$ so we arrive at the inequality
$\mbox{$\mathfrak A$} _2\le \mbox{$\mathfrak C$} \, \|a\|_p\, T ^{-1} .$
On combining the definition of \mbox{$\mathfrak A$} in (\ref{T7.10.50}), the estimate of \mbox{$\mathfrak A$}
in the case $p >3/2$ in (\ref{T7.10.60}), the inequalities $\mbox{$\mathfrak A$} \le \mbox{$\mathfrak A$} _1+
\mbox{$\mathfrak A$} _2$ and $\mbox{$\mathfrak A$} _1\le \mbox{$\mathfrak C$} \,\|a\|_p\, T ^{-1} \,
(1+\ln T)$ (see above), and the preceding estimate of $\mbox{$\mathfrak A$} _2$,
we find
$\mbox{$\mathfrak A$} \le \mbox{$\mathfrak C$} \,\|a\|_p\, T ^{-\varrho (p)}\, \sigma (p,T),$
with
$\varrho (p):=3/(2\, p),\; \sigma (p,T):=1$ if $p>3/2,$
and
$\varrho (p):=1,\; \sigma (p,T):=1+\ln T$ in the case $p\le3/2$.
But $T\ge 1$, so we get in any case that
$
\mbox{$\mathfrak A$} \le \mbox{$\mathfrak C$} \,\|a\|_p\, T ^{-3/(2p)+1/2}.
$
Since equation (\ref{L10.10.10}) is valid with
$
\mbox{$\mathfrak I$}^{( \tau )} (a)(x,\, \cdot \, )
$
in the role of $\phi $ and with $T$ replaced by $T/2$, for a. e. $x \in \partial \Omega $
(Corollary \ref{corollaryC10.66}), we may deduce from the preceding estimate and
(\ref{T7.10.30}) -- (\ref{T7.10.50}) that
$
\bigl( \int_{ T}^{ \infty }  \int_{ \partial \Omega }| \partial _4 ^{1/2} \mbox{$\mathfrak I$}^{( \tau )} (a)(x,t)|^2
\, do_x\, dt \bigr)  ^{1/2}
\le
\mbox{$\mathfrak C$} \,\|a\|_p\, T ^{-3/(2p)+1/2}.
$
But by Corollary \ref{corollaryC10.66}, equation (\ref{C10.65.10}) holds for
$\Psi := \mbox{$\mathfrak I$} ^{( \tau )} (a)|S_{\infty }$,
so Theorem \ref{theoremT7.10} follows from
from (\ref{T7.10.20}) and the estimate of
$\partial _4 ^{1/2} \bigl(\, \mbox{$\mathfrak I$}^{( \tau )} (a)|S_{\infty }\,\bigr) $ just proved.
\hfill $\Box $

\vspace{1ex}
A pointwise estimate with respect to the asymptotics of $\mbox{$\mathfrak I$}^{( \tau )} (a)$
in time and in space is readily available, due to Theorem \ref{theoremT4.60} and Lemma
\ref{lemmaL10.90}.
\begin{corollary} \label{corollaryC7.10}
Take $R_0,\, \delta _0,\, \kappa _0$ and $a$ as in Theorem \ref{theoremT4.60}.
Let $p \in [1, \infty ] $ and suppose that $a \in L^p( \overline{ \Omega }^c)^3$.
Let $R \in (R_0, \infty ),\; \alpha \in \mathbb{N} _0^3$ with $| \alpha |\le 1.$
Then, for $\epsilon \in [0,1],\; x \in B_R^c,\; t \in (0, \infty ) $,
\begin{eqnarray*}&&\hspace{-2em}
|\mbox{$\mathfrak I$}^{( \tau )} (a)(x,t)|
\le
\mbox{$\mathfrak C$} \,
\bigl[\, 
( \delta _0+\| a| \Omega _{R_0}\|_1)\,
\bigl(\, |x|\,\nu (x) \,\bigr) ^{(-1-| \alpha |/2)} \, (1+t)^{-2 }
\\&&\hspace{2em}
+
( \delta _0+\| a| \Omega _{R_0}\|_1)^{1-\epsilon }\, \|a\|_p^{\epsilon }\,
\bigl(\, |x|\,\nu (x) \,\bigr) ^{(-1-| \alpha |/2)\,(1-\epsilon )}\, (1+t)^{(-3/(2p)-| \alpha |/2)\, \epsilon }
\,\bigr]
.
\end{eqnarray*}
\end{corollary}
{\bf Proof:}
In the case $t\le 1$, Theorem \ref{theoremT4.60} immediately yields
that $|\mbox{$\mathfrak I$}^{( \tau )} (a)(x,t)|$ is bounded by
$
\mbox{$\mathfrak C$} \,
( \delta _0+\| a| \Omega _{R_0}\|_1)
\bigl(\, |x|\,\nu (x) \,\bigr) ^{-1-| \alpha |/2} \, (1+t)^{-2}.
$
If
$t\ge 1$,
the looked-for estimate holds due to Theorem \ref{theoremT4.60}, (\ref{L10.100.10}),
Lemma \ref{lemmaL10.90} and because $d=d^{\epsilon }\, d^{1-\epsilon }$ for $d>0$.
\hfill $\Box $

\section{Temporal decay of the potential $\mbox{$\mathfrak V$}^{( \tau )} ( \phi )$ and
of the solution of the integral equation (\ref{T10.60.10}).}
\setcounter{equation}{0}

A key element of our theory is a decay estimate of the solution to the integral equation
(\ref{T10.60.10}). This element will be presented in Theorem \ref{theoremT3.10} below.
Its proof depends on certain features of the potential $\mbox{$\mathfrak V$} ^{( \tau )} ( \phi ),$
which we establish in the ensuing two lemmas.
\begin{lemma} \label{lemmaL5.10}
Let $\mu,\, T \in (0, \infty ) $ and $\phi \in L^2(S_T)^3$. Then
\begin{eqnarray*} 
\Bigl( \int_{ T+\mu} \int_{ \partial \Omega }| \partial _t ^{1/2} \mbox{$\mathfrak V$} ^{( \tau )} ( \phi )(x,t)|^2
\, do_x \, dt \Bigr) ^{1/2}
\le
\mbox{$\mathfrak C$} \,(\mu^{-3/2}+\mu ^{-1} )\, T ^{1/2} \,\| \phi \|_2.
\end{eqnarray*}
(By Corollary \ref{corollaryC10.70}, the fractional derivative
$\partial _t ^{1/2} \mbox{$\mathfrak V$} ^{( \tau )} ( \phi )(x,t)$
exists for a. e. $x \in \partial \Omega $ and for $t \in (T, \infty )$).
\end{lemma}
{\bf Proof:}
For brevity we set $v:=\mbox{$\mathfrak V$} ^{( \tau )} ( \phi )$.
By Corollary \ref{corollaryC10.70}, equation (\ref{L10.10.10}) holds with
$v(x,\, \cdot \,)$ in the place of the function $\phi $ in Lemma \ref{lemmaL10.10}, for a. e. $x \in \partial
\Omega $. With this reference in mind, we observe that by Lemma \ref{lemmaL10.120} and \ref{lemmaL10.50},
for $r \in (T, \infty ),\; x \in \partial \Omega ,\; l \in \{0,\, 1\} ,$
\begin{eqnarray*} \label{L5.10.30}&&
| \partial _r^l v(x,r)|
\le
\int_{ 0}^T \int_{ \partial \Omega } \bigl(\, (r-\sigma )^{-3/2-l}+(r-\sigma )^{-3/2-l/2}\,\bigr) \,
| \phi (y, \sigma )|\, do_y\, d \sigma
\\&&\nonumber
\le
\mbox{$\mathfrak C$} \, \bigl(\, (r-T )^{-3/2-l}+(r-T )^{-3/2-l/2}\,\bigr) \,
T ^{1/2}  \,\| \phi\|_2,
\end{eqnarray*}
where we used that $\| \phi \|_1\le C \, T ^{1/2} \,\| \phi \|_2.$
The preceding inequality implies
\begin{eqnarray} \label{L5.10.40}&&
\Bigl( \int_{ T+\mu}^{ \infty } \int_{ \partial \Omega }| (t-T) ^{-1/2} \,v \bigl(\, x,\,(t+T)/2 \,\bigr) |^2
\, do_x\, dt \Bigr) ^{1/2}
\\&&\nonumber
\le
\mbox{$\mathfrak C$} \, T ^{1/2} \,\| \phi \|_2\,
\Bigl( \int_{ T+\mu}^{ \infty }  (t-T) ^{-4}\, dt \Bigr) ^{1/2}  
\le
\mbox{$\mathfrak C$} \, T ^{1/2} \,\| \phi \|_2\,\mu^{-3/2},
\\ \label{L5.10.50}&&\hspace{-3em}
\Bigl( \int_{ T+\mu}^{ \infty } \int_{ \partial \Omega }\Bigl|  \int_{ (t+T)/2}^t (t-r) ^{-1/2} \,
\partial _rv ( x, r) \, dr \Bigr| ^2 \, do_x\, dt \Bigr) ^{1/2}
\\&&\nonumber \hspace{-3em}
\le
\mbox{$\mathfrak C$} \, T ^{1/2} \,\| \phi \|_2\,
\Bigl( \int_{ T+\mu}^{ \infty } \Bigl[  \int_{ (t+T)/2}^t (t-r) ^{-1/2} \,
\bigl(\, (t-T) ^{-5/2} + (t-T) ^{-2} \,\bigr) \, dr \Bigr]^2\, dt \Bigr) ^{1/2} 
\\&&\nonumber \hspace{-3em}
\le
\mbox{$\mathfrak C$} \, T ^{1/2} \,\| \phi \|_2\,(\mu^{-3/2}+\mu ^{-2} ).
\end{eqnarray}
For $t \in [T+\mu ,\, \infty ),\; r \in \bigl(\, 0,\, (t+T)/2 \,\bigr],$
the inequality $t-r\ge (t-T)/2$ holds. Therefore
\begin{eqnarray}\label{L5.10.55} &&
\mbox{$\mathfrak B$}
:=
\Bigl( \int_{ T+\mu}^{ \infty } \int_{ \partial \Omega }\Bigl|  \int^{ (t+T)/2}_0 (t-r) ^{-3/2} \,
v ( x, r) \, dr \Bigr| ^2 \, do_x\, dt \Bigr) ^{1/2}
\\&&\nonumber
\le
\mbox{$\mathfrak C$} \,
\Bigl( \int_{ T+\mu}^{ \infty } (t-T) ^{-3}  \, \int_{ \partial \Omega }
\Bigl[ \int^{ (t+T)/2}_0 |v ( x, r)| \, dr \Bigr] ^2\, do_x\, dt \Bigr) ^{1/2} .
\end{eqnarray}
By extending the domain of integration of the variable $r$ to $(0, \infty ) $,
we may separate the integration  with respect to $r$ and $t$. In this way we get
\begin{eqnarray} \label{L5.10.60}
\mbox{$\mathfrak B$}
\le
\mbox{$\mathfrak C$} \, \mu ^{-1} \,
\Bigl( \int_{ \partial \Omega }\Bigl[ \int_{ 0}^{ \infty } |v(x,r)|\, dr \Bigr]^2\, do_x \Bigr) ^{1/2} .
\end{eqnarray}
But by Corollary \ref{corollaryC10.60} with $K=diam\, \Omega $, 
the term
$
\int_{ 0 }^{ \infty } |v(x,r)|\, dr
$
for $x \in \partial \Omega $
is bounded by
$
\mbox{$\mathfrak C$} \,
\int_{ 0}^{ \infty } \int_{ 0}^T \int_{ \partial \Omega }(|x-y|^2+r-\sigma )^{-3/2}\, | \phi (y, \sigma )|\, do_y\,
d \sigma \, dr.
$
Integrating with respect to $r$, this triple integral, in turn, may be estimated by 
$
\mbox{$\mathfrak C$} \, \int_{ 0}^T \int_{ \partial \Omega }|x-y| ^{-1} \,| \phi (y, \sigma )|\, do_y\, d \sigma,
$
and hence with H\"older's inequality by
$
\mbox{$\mathfrak C$} \, T ^{1/2} \, \int_{ \partial \Omega }|x-y| ^{-1} \,\| \phi (y,\, \cdot \,)\|_2\, do_y.
$
On applying H\"older's inequality once more, we may thus deduce from (\ref{L5.10.60}) that
\begin{eqnarray*} &&
\mbox{$\mathfrak B$} 
\le
\mbox{$\mathfrak C$} \, \mu ^{-1} \,T ^{1/2} \,
\Bigl( \int_{ \partial \Omega } \Bigl[ \int_{ \partial \Omega }|x-y| ^{-1} \,do_y \Bigr]\,
\Bigl[ \int_{ \partial \Omega }|x-y| ^{-1} \,\| \phi (y,\, \cdot \,)\|_2^2\,do_y \Bigr]\, do_x \Bigr) ^{1/2}
\\&&
\le
\mbox{$\mathfrak C$} \, \mu ^{-1} \, T ^{1/2} \,
\Bigl( \int_{ \partial \Omega } 
\int_{ \partial \Omega }|x-y| ^{-1} \, do_x \,\| \phi (y,\, \cdot \,)\|_2^2\,do_y  \Bigr) ^{1/2}
\le
\mbox{$\mathfrak C$} \, \mu ^{-1} \, T ^{1/2} \,\| \phi \|_2.
\end{eqnarray*}
The lemma follows from this inequality, (\ref{L5.10.40}) -- (\ref{L5.10.55}) and, as explained
at the beginning of this proof, from (\ref{L10.10.10}).
\hfill $\Box $
\begin{lemma} \label{lemmaL5.20}
As in Lemma \ref{lemmaL5.10}, let $\mu,\, T \in (0, \infty ) $ and $\phi \in L^2(S_T)^3$. Then
\begin{eqnarray} \label{L5.20.10}
&&
\Bigl( \int_{ T+\mu}^{ \infty } \Bigl[\|\mbox{$\mathfrak V$} ^{( \tau )} ( \phi )(t)| \partial \Omega
\|_{H^1( \partial \Omega )^3}^2
+
\| n^{ (\Omega )} \cdot 
\partial _t  \mbox{$\mathfrak V$} ^{( \tau )} ( \phi )(t)\|_2^2
\Bigr] \, dt \Bigr) ^{1/2}
\\&&\nonumber
\le
\mbox{$\mathfrak C$} \,(\mu^{-2}+\mu ^{-1} )\, T ^{1/2} \,\| \phi \|_2.
\end{eqnarray}
(Recall that by Lemma \ref{lemmaL10.120}, we have
$\mbox{$\mathfrak V$} ^{( \tau )} ( \phi )| \mathbb{R}^3 \times (T, \infty )
\in C^1 \bigl(\, \mathbb{R}^3 \times (T, \infty ) \,\bigr) .$)
\end{lemma}
{\bf Proof:}
We again set $v:=\mbox{$\mathfrak V$} ^{( \tau )} ( \phi )$ for brevity.
Lemma \ref{lemmaL10.50} and \ref{lemmaL10.120} yield for $x \in \partial \Omega $ and
$t \in [T+\mu,\, \infty )$
that
\begin{eqnarray} \label{L5.20.20}&&
|v(x,t)| + \sum_{j = 1}^3 | \partial _jv(x,t)|+| \partial _tv(x,t)|
\\&&\nonumber
\le
\mbox{$\mathfrak C$} \, \int_{ 0}^T \int_{ \partial \Omega }
\bigl(\, (t-\sigma )^{-3/2} +(t-\sigma )^{-2}+(t-\sigma )^{-5/2}\,\bigr)  \,| \phi (y, \sigma )|\,do_y\, d \sigma
\\&&\nonumber
\le
\mbox{$\mathfrak C$} \, 
\bigl(\, (t-T )^{-3/2} +(t-T )^{-5/2}\,\bigr)  \,\| \phi \|_1
\le
\mbox{$\mathfrak C$} \, T ^{1/2} \, 
\bigl(\, (t-T )^{-3/2} +(t-T )^{-5/2}\,\bigr)  \,\| \phi \|_2.
\end{eqnarray}
Denote the left-hand side of (\ref{L5.20.10}) by \mbox{$\mathfrak B$}.
We may deduce from (\ref{L5.20.20}) that \mbox{$\mathfrak B$} is bounded by
$
\mbox{$\mathfrak C$} \, T ^{1/2} \, \| \phi \|_2
\bigl(\, \int_{ T+\mu}^{ \infty } 
\bigl[\, (t-T )^{-3} +(t-T )^{-5}\,\bigr]  \,dt \,\bigr) ^{1/2},
$
and hence by
$
\mbox{$\mathfrak C$} \, T ^{1/2} \, \| \phi \|_2\, ( \mu ^{-1} +\mu ^{-2} ).
$
\hfill $\Box $
\begin{corollary} \label{corollaryC5.10}
There is a constant $c_3>0$ such that
$
\|\mbox{$\mathfrak V$} ^{( \tau )} ( \phi )|S_{T+\mu, \infty }\|_{H_{T+\mu, \infty }}
$
is bounded by
$
c_3\, (1+T) ^{1/2} \, (1+\mu ) ^{-1} \, \| \phi \|_2
$,
for $T,\,\mu \in (0, \infty ) $ and $\phi \in L^2(S_T)^3$.
\end{corollary}
{\bf Proof:}
First suppose that $\mu \in (0,1]$.
Then Corollary \ref{corollaryC10.40} and Theorem \ref{theoremT10.20} yield that
$
\|\mbox{$\mathfrak V$} ^{( \tau )} ( \phi )|S_{T+\mu, \infty }\|_{H_{T+\mu, \infty }}
\le
\|\mbox{$\mathfrak V$} ^{( \tau )} ( \phi )|S_{\infty }\|_{H_{\infty }}
\le
\mbox{$\mathfrak C$} \, \| \phi \|_2
\le
\mbox{$\mathfrak C$} \, \| \phi \|_2\,(1+\mu) ^{-1} \, (1+T) ^{1/2} .
$
If $\mu\ge 1$, we use Lemma \ref{lemmaL5.10}, \ref{lemmaL5.20} and Corollary \ref{corollaryC10.70} to obtain
that
$
\|\mbox{$\mathfrak V$} ^{( \tau )} ( \phi )|S_{T+\mu, \infty }\|_{H_{T+\mu, \infty }}
$
is bounded by
$
\mbox{$\mathfrak C$} \, T ^{1/2} \, \| \phi \|_2\, ( \mu ^{-1} +\mu ^{-2} ),
$
and thus again by
$
\mbox{$\mathfrak C$} \, \| \phi \|_2\,(1+\mu) ^{-1}  \, (1+T) ^{1/2} .
$
\hfill $\Box $
\begin{theorem} \label{theoremT3.10}
Let $\widetilde{ b} \in H_{\infty }$. Suppose there are numbers $\delta \in (0, \infty ) ,\; \zeta \in (0,1)$
such that
$
\| \widetilde{ b}|S_{T, \infty }\|_{H_{T, \infty }}\le \delta \, T^{-\zeta } \quad \mbox{for}\;\; T \in (1, \infty ).
$
Let
$\phi $
be the unique function from
$L^2_n(S_{\infty })$
such that
$\mbox{$\mathfrak V$} ^{( \tau )} ( \phi )|S_{\infty }=\widetilde{ b}$ (Theorem \ref{theoremT10.40}).
Then $\| \phi |S_{T, \infty }\|_2
\le
\mbox{$\mathfrak C$}\, (1+T)^{-\zeta }$ for $T \in (0, \infty ) $.
\end{theorem}
{\bf Proof:}
For brevity, set
$\epsilon :=1-\zeta $. Let $n \in \mathbb{N} $ with
$n\ge 400\, c_1^2\, c_2^2$,
where $c_1$ was introduced in Theorem \ref{theoremT10.20}, and $c_2$ in Theorem \ref{theoremT10.30}.
Let $k $ be the unique number from $\mathbb{N} $ such that $k \, \epsilon < 1\le (k+1)\, \epsilon $.
Define $\varphi ( \epsilon )$ as in Lemma \ref{lemmaL1.10}. In view of a later proof by induction,
suppose that $T,\, B \in (0, \infty ) $ are such that
\begin{eqnarray} \label{T3.10.25}\label{T3.10.30}&&\hspace{-3em}
T
\ge
\max \{ 4,\; (160\, c_2\, c_3\, n)^2,\; (320\, c_2\, c_3 k)^{1/ \varphi ( \epsilon )}\},
\\&&\label{T3.10.x}\hspace{-3em}
B
\ge
\max\{ 40\, c_2\, c_3\, (\| \phi |S_1 \|_2+1),\; 20\, \delta \},
\quad 
\| \phi |S_{t, \infty }\|_2\le B\,t^{-1+\epsilon }\; \mbox{for}\;\; t \in [1,T].
\end{eqnarray}
We want to show that $\| \phi |S_{T, \infty }\|_2\le (B/2)\,T^{-1+\epsilon }$
in this situation, hence
$\| \phi |S_{t, \infty }\|_2\le B\,t^{-1+\epsilon }$ for $t \in   [1,\,2\, T].$
(Note that $S_{t, \infty }\subset S_{T, \infty }$ for $t \in [T,\, 2\,T]$.)
With this aim in mind, we refer to Lemma \ref{lemmaL1.20}, choosing a number
$i_0 \in \{0,\, ...,\, n-1\} $ such that the inequality
$
\| \phi | \partial \Omega \times (t_0-T/(2\, n),\, t_0)\|\le \| \phi |S_T \backslash S_{T/2 }\|_2\,n ^{-1/2} 
$
holds
with $t_0:=(T/2)\, \bigl(\, 1+(i_0+1)/n \,\bigr) .$
Since $t_0\le T$ and $\phi |S_{t_0, \infty }\in L^2_n(S_{t_0, \infty })$,
we obtain with Corollary \ref{corollaryC10.90} that 
\begin{eqnarray} \label{T3.10.50} &&
\| \phi |S_{T, \infty }\|_2
\le
\| \phi | S_{t_0, \infty }\|_2
\le
c_2\,\|\mbox{$\mathfrak V$} ^{( \tau )} ( \phi |S_{t_0,  \infty } )|S_{t_0,\infty }\|_{H_{t_0, \infty }}
\\&& \nonumber
\le
c_2\, \bigl(\, \|\mbox{$\mathfrak V$} ^{( \tau )} ( \phi )|S_{t_0,\infty }\|_{H_{t_0, \infty }}
+
\|\mbox{$\mathfrak V$} ^{( \tau )} ( \phi |S_{t_0})|S_{t_0,\infty }\|_{H_{t_0, \infty }}\,\bigr) .
\end{eqnarray}
where the last inequality holds because
$
\mbox{$\mathfrak V$} ^{( \tau )} ( \phi )=
\mbox{$\mathfrak V$} ^{( \tau )} ( \phi |S_{t_0} )
+
\mbox{$\mathfrak V$} ^{( \tau )} ( \phi |S_{t_0,\infty }).
$
Recalling that $\mbox{$\mathfrak V$} ^{( \tau )} ( \phi )|S_{\infty }=\widetilde{ b}$,
and splitting $\mbox{$\mathfrak V$} ^{( \tau )} ( \phi |S_{t_0} )$
into a suitable sum, we may conclude from (\ref{T3.10.50}) that
$
\| \phi |S_{T,\infty }\|_2
\le
c_2\, (\| \widetilde{ b}|S_{t_0, \infty }\|_{H_{t_0, \infty }} + \mbox{$\mathfrak A$} _1+\mbox{$\mathfrak A$} _2
+ \sum_{j = 0}^{k-1}\mbox{$\mathfrak B$} _j+\mbox{$\mathfrak A$} _3),
$
where the terms
$
\mbox{$\mathfrak A$} _1,\, \mbox{$\mathfrak A$} _2,\,
\mbox{$\mathfrak B$} _0,\, ...,\, \mbox{$\mathfrak B$} _{k-1},\, \mbox{$\mathfrak A$} _3
$
all have the form
$\| \mbox{$\mathfrak V$} ^{( \tau )} (  \phi | A)|S_{t_0,  \infty }\|_{H_{t_0, \infty }},$
with $A$ being defined by
$
A=\partial \Omega \times \bigl(\, t_0-T/(2\, n),\, t_0 \,\bigr),\;
A=S_{t_0-T/(2\, n)} \backslash S_{T/4}
$
and
$
A=S_1
$
in the case of
$
\mbox{$\mathfrak A$} _1,\,\mbox{$\mathfrak A$} _2
$
and
$
\mbox{$\mathfrak A$} _3,
$
respectively, and
$
A=S_{(T/4)^{1-j/k}} \backslash S_{(T/4)^{1-(j+1)/k}}
$
in the case of
$\mbox{$\mathfrak B$} _j$, for $0\le j\le k-1$.
The definition of
$\mbox{$\mathfrak B$} _j$
makes sense because
$T\ge 4$ (see (\ref{T3.10.25})),
so
$
(T/4)^{1-j/k}> (T/4)^{1-(j+1)/k}
$
for $j$ as before.
But $t_0\ge T/2$, so by Corollary \ref{corollaryC10.40} and
the preceding estimate of $\| \phi |S_{T,\infty }\|_2$,
\begin{eqnarray} \label{T3.10.65}&&
\| \phi |S_{T,\infty }\|_2
\le
c_2\, (\| \widetilde{ b}|S_{T/2,\, \infty }\|_{H_{T/2,\, \infty }} + \mbox{$\mathfrak A$} _1+\mbox{$\mathfrak A$} _2
+ \sum_{j = 0}^{k-1}\widetilde{ \mbox{$\mathfrak B$}} _j+\widetilde{ \mbox{$\mathfrak A$}} _3),
\end{eqnarray}
where
$
\widetilde{ \mbox{$\mathfrak B$}} _0,\, ...,\,
\widetilde{ \mbox{$\mathfrak B$}} _{k-1},\,
\widetilde{ \mbox{$\mathfrak A$}} _3
$
are defined in the same way as
$
\mbox{$\mathfrak B$} _0,\, ...,\,
\mbox{$\mathfrak B$} _{k-1},\,
\mbox{$\mathfrak A$} _3,
$
respectively, except that the restriction to $S_{t_0, \infty }$ is replaced  by a restriction
to $S_{T/2, \,\infty  }$, and the mapping
$\|\; \|_{H_{t_0 ,\infty }}$ by $\|\;\|_{H_{T/2,\, \infty }}$.
Corollary \ref{corollaryC10.80} and
the choice of
$t_0$
imply
\begin{eqnarray*} &&
c_2\, \mbox{$\mathfrak A$} _1
\le
c_1\,c_2 \,\| \phi | \partial \Omega \times \bigl(\, t_0-T/(2\, n),\, t_0 \,\bigr) \|_2
\le
c_1\,c_2 \,\| \phi | S_T \backslash S_{T/2} \|_2\, n ^{-1/2}
\\&&
\le
c_1\,c_2 \,\| \phi |S_{T/2,\, \infty } \|_2\, n ^{-1/2}
,
\end{eqnarray*}
so by (\ref{T3.10.x}) and the choice of $n$ we get
$
c_2\, \mbox{$\mathfrak A$} _1
\le
c_1\, c_2\, B\, (T/2)^{-1+\epsilon }\, n ^{-1/2}
\le
(B/10)\, T^{-1+\epsilon }.
$
Concerning
$\mbox{$\mathfrak A$} _2$, we use Corollary \ref{corollaryC5.10} with $T$ replaced by
$t_0-T/(2\, n)$ and $\mu $ by $T/(2\, n)$,
to obtain
\begin{eqnarray*}
&&
c_2\, \mbox{$\mathfrak A$} _2
\le
c_2\, c_3\, \bigl(\, 1+T/(2\, n) \,\bigr) ^{-1} \, \bigl(\, 1+t_0-T/(2\, n) \,\bigr) ^{1/2}
\, \| \phi  | S_{t_0-T/(2\, n)}\backslash S_{T/4}\|_2.
\\&&
\le
2\,c_2\, c_3\, n\, T ^{-1} \, ( 1+T) ^{1/2}
\, \| \phi  | S_{T/4,\, \infty }\|_2
.
\end{eqnarray*}
But
$T\ge 1$,
hence $(1+T) ^{1/2} \le 2 \, T ^{1/2} $, so it follows from
the preceding estimate of $c_2\, \mbox{$\mathfrak A$} _2$ and from (\ref{T3.10.x}) with $t=T/4$ that
$
c_2\, \mbox{$\mathfrak A$} _2
\le
4\, c_2\, c_3\, n\, T ^{-1/2} \, \| \phi  | S_{T/4,\, \infty }\|_2
\le
16\,c_2\, c_3\, n\, T ^{-3/2+\epsilon } \, B.
$
By the choice of $T$ in (\ref{T3.10.30}), we may conclude that
$c_2 \, \mbox{$\mathfrak A$} _2\le (B/10)\,T^{-1+\epsilon }$.
Moreover, by Corollary \ref{corollaryC5.10}
with $ T,\; \mu $ replaced by $T/2-1$ and $1$, respectively, we get
$
c_2\, \widetilde{ \mbox{$\mathfrak A$}} _3\le 4\, c_2\, c_3\, T ^{-1} \,\| \phi |S_1\|_2,
$
hence 
$
c_2\, \widetilde{ \mbox{$\mathfrak A$}} _3\le (B/10)\, T ^{-1} \le (B/10)\, T^{-1+\epsilon }
$
by the choice of $B$ in (\ref{T3.10.x}) and because $T\ge1$.
Due to
the assumptions on $\widetilde{ b}$, the definition of $\epsilon $
and because $T>2$, we obtain
$
c_2\,\| \widetilde{ b}|S_{T/2,\, \infty }\|_{H_{T/2,\, \infty }}
\le
\delta (T/2)^{-1+\epsilon }
,
$
so that
$
c_2\,\| \widetilde{ b}|S_{T/2,\, \infty }\|_{H_{T/2,\, \infty }}\le (B/10)\, T^{-1+\epsilon }
$
by the choice of $B$ in (\ref{T3.10.x}).
This leaves us to estimate $c_2\, \sum_{j = 0}^{k-1}\widetilde{ \mbox{$\mathfrak B$} }_j$.
To this end let $j \in \{0,\, ...,\, k-1\} .$ Then Corollary \ref{corollaryC5.10} with
$(T/4)^{1-j/k}$ in the role of $T$ and $T/2-(T/4)^{1-j/k}$ in that of $\mu$ implies that
$c_2\, \widetilde{ \mbox{$\mathfrak B$} }_j$ is bounded by
\begin{eqnarray*}
c_2\,c_3 \,\bigl(\, 1+T/2-(T/4)^{1-j/k} \,\bigr) ^{-1} \, \bigl(\, 1+(T/4)^{1-j/k} \,\bigr) ^{1/2}
\,\| \phi |S_{(T/4)^{1-j/k}} \backslash S_{(T/4)^{1-(j+1)/k} }\|_2.
\end{eqnarray*}
But $T\ge 4$, so 
$
\bigl(\, 1+(T/4)^{1-j/k} \,\bigr) ^{1/2}
\le
2 ^{1/2} \, (T/4)^{1/2-j/(2k)}
$
and
$
1+T/2-(T/4)^{1-j/k}\ge 1+T/4\ge T/4.
$
Therefore we get
$
c_2\, \widetilde{ \mbox{$\mathfrak B$} }_j
\le
8\,c_2\,c_3  \, T^{-1/2-j/(2k)} \,\| \phi |S_{(T/4)^{1-(j+1)/k},\, \infty  }\|_2
$,
so it follows with (\ref{T3.10.x}) that
\begin{eqnarray*}
c_2\, \widetilde{ \mbox{$\mathfrak B$} }_j
\le
8\,c_2\,c_3\, T^{-1/2-j/(2k)} \, B\, (T/4)^{(1-(j+1)/k)\,(-1 +\epsilon )}
\le
32\,c_2\,c_3\,B  \, T^{Z(j)} \,T^{-1+\epsilon }.
\end{eqnarray*}
where $Z(j):=-1/2+j/(2k)+1/k -\epsilon \,(j+1)/k$.
But $Z(j)\le -\varphi ( \epsilon )$
by Lemma \ref{lemmaL1.10}, with $\varphi ( \epsilon )$ defined there.
Thus
$
c_2\, \widetilde{ \mbox{$\mathfrak B$} }_j
\le
32\,c_2\,c_3\,B  \, T^{-\varphi ( \epsilon )} \,T^{-1+\epsilon },
$
hence by the choice of $T$ in (\ref{T3.10.25}),
$
c_2\, \widetilde{ \mbox{$\mathfrak B$} }_j
\le
T^{-1+\epsilon }\,B/(10\, k) .
$
Since this holds for any $j \in \{0,\, ...,\, k-1\} ,$
we thus get
$
c_2\, \sum_{j = 0}^{k-1} \widetilde{ \mbox{$\mathfrak B$} }_j
\le
(B/10) \,T^{-1+\epsilon }.
$
Combining (\ref{T3.10.65}) with the preceding estimates of
$
\| \widetilde{ b}|S_{T/2,\, \infty }\|_{H_{T/2,\, \infty }}
,\; \mbox{$\mathfrak A$} _1,\; \mbox{$\mathfrak A$} _2,\;
\sum_{j = 0}^{k-1} \widetilde{ \mbox{$\mathfrak B$} }_j
$
and $\widetilde{ \mbox{$\mathfrak A$} }_3$,
we get
$
\| \phi |S_{T,  \infty }\|_2 \le (B/2)\, T^{-1+\epsilon }
,
$
so we have in fact shown that if
$\| \phi |S_{t, \infty }\|_2\le B\,t^{-1+\epsilon }$
for
$
t \in [1,T]
$
(see (\ref{T3.10.x})),
then
$
\| \phi |S_{T,  \infty }\|_2 \le (B/2)\, T^{-1+\epsilon }
.
$

Now, in view of a proof by induction, put
$
T_0
:=
\max \{ 4,\; (160\, c_2\, c_3\, n)^2,\; (320\, c_2\, c_3 k)^{1/ \varphi ( \epsilon )}\}
$
and
$
B_0
:=
\max\{ 40\, c_2\, c_3\, (\| \phi |S_1 \|_2+1),\; 20\, \delta ,\; \| \phi \|_2\, T_0^{1- \epsilon },\;
2\,\| \phi \|_2\}.
$
Then for $t \in [1, T_0],$ we find
$
\| \phi |S_{t, \infty }\|_2 \le \| \phi \|_2
\le
\| \phi \|_2\, T_0^{1-\epsilon }\, t^{-1+\epsilon } \le B_0\, t^{-1+\epsilon } .
$
Suppose that $m \in \mathbb{N} _0 $ and $\| \phi |S_{t, \infty }\|_2 \le B_0\, t^{-1+\epsilon }$
for $t \in [1,\, 2^m\, T_0].$ Since $B_0$ fulfills the condition on $B$ in (\ref{T3.10.x}),
as does $2^m\,T_0$ the one on $T$ in (\ref{T3.10.25}), we may apply
the first part of this proof
with $2^m\,T_0$ and $B_0$ in the place of $T$ and $B$, respectively. It follows that
$
\| \phi |S_{2^m\,T_0,\,  \infty }\|_2 \le (B_0/2)\, (2^m\,T_0)^{-1+\epsilon }.
$
Thus for $t \in [2^m\,T_0,\, 2^{m+1}\,T_0]$,
\begin{eqnarray*} 
\| \phi |S_{t, \infty }\|_2
\le
\| \phi |S_{2^m\,T_0,\, \infty }\|_2
\le
(B_0/2)\,(2^m\,T_0)^{-1+\epsilon }
\le
B_0\,(2^{m+1}\,T_0)^{-1+\epsilon }
\le
B_0\,t^{-1+\epsilon }.
\end{eqnarray*}
Thus we have shown by induction that
$\| \phi |S_{t, \infty }\|_2 \le B_0\, t^{-1+\epsilon }$
for any $t \in [1,\, 2^m\,T_0]$
and any $m \in \mathbb{N} _0$,
so the preceding inequality is valid for any $t \in [1, \infty )$.
Hence, for such $t$ we find
$
\| \phi |S_{t, \infty }\|_2
\le 2\, B_0\, (1+t)^{-1+\epsilon }.
$
But for $t \in (0,1],$ it is obvious that
$
\| \phi |S_{t, \infty }\|_2 \le \| \phi \|_2 
\le
2\,\| \phi \|_2 \, (1+t)^{-1+\epsilon }
\le
B_0\,(1+t)^{-1+\epsilon }.
$
This completes the proof of Theorem \ref{theoremT3.10}. \hfill $\Box $

\vspace{1ex}
We will combine Theorem \ref{theoremT3.10} with the following pointwise temporal and spatial
estimate of $\mbox{$\mathfrak V$} ^{( \tau )} (  \phi )$.
\begin{lemma} \label{lemmaL5.30}
Let $R \in (0, \infty ) $ with $\overline{ \Omega }\subset B_R,\;
\phi \in L^2(S_{\infty })^3$ and $\alpha \in \mathbb{N} _0^3$ with $| \alpha |\le 1$.
Then, for $x \in B_R^c,\; t \in (0, \infty ) $,
\begin{eqnarray*}
| \partial _x^{\alpha }\mbox{$\mathfrak V$} ^{( \tau )} (  \phi )(x,t)|
\le
\mbox{$\mathfrak C$} \,
\bigl[\,
\bigl(\, |x|\,\nu (x) \,\bigr) ^{-1-| \alpha |/2}\,\| \phi |S_{t/2,\, \infty }\|_2
+
\bigl(\, |x|\,\nu (x) +t\,\bigr) ^{-1-| \alpha |/2}\,\| \phi \|_2
 \,\bigr].
\end{eqnarray*}
\end{lemma}
{\bf Proof:}
Choose $R_0 \in (0,R)$ with $\overline{ \Omega }\subset B_{R_0}.$
Take $x,\, t$ as in the lemma. For $y \in \partial \Omega $, we have $|y|\le R_0$, so we find as in
(\ref{L6.40.10}) that
$
|x-y|\ge |x|\,(1-R_0/R)\ge R-R_0>0,
$
and by Lemma \ref{lemmaL10.1},
$\nu (x-y) ^{-1} \le C\,|y|\, \nu (x) ^{-1} \le C\, R_0\,\nu (x) ^{-1} .$
Hence by (\ref{T10.10.10}) and Corollary \ref{corollaryC10.60} with $K=R-R_0$,
\begin{eqnarray*} &&
| \partial _x^{\alpha }\mbox{$\mathfrak V$} ^{( \tau )} (  \phi )(x,t)|
\le
C(R,R_0, \tau )\, \int_{ 0} ^t \int_{ \partial \Omega }
\bigl(\, |x|\,\nu (x) +t-\sigma \,\bigr) ^{-3/2-| \alpha |/2}\,| \phi (y, \sigma )|\, do_y\, d \sigma
,
\end{eqnarray*}
so
$
| \partial _x^{\alpha }\mbox{$\mathfrak V$} ^{( \tau )} (  \phi )(x,t)|
\le
\mbox{$\mathfrak C$} \,
\int_{ 0} ^t \bigl(\, |x|\,\nu (x) +t-\sigma \,\bigr) ^{-3/2-| \alpha |/2}\,\| \phi (\sigma )\|_2\,  d \sigma .
$
By H\"older's inequality,
\begin{eqnarray*} &&\hspace{-2em}
\int_{ 0} ^{t/2}
\bigl(\, |x|\,\nu (x) +t-\sigma \,\bigr) ^{-3/2-| \alpha |/2}\,\| \phi (\sigma )\|_2\,  d \sigma
\le
\Bigl(
\int_{ 0} ^{t/2}
\bigl(\, |x|\,\nu (x) +t-\sigma \,\bigr) ^{-3-| \alpha |}\,d \sigma
\Bigr)  ^{1/2} \,\| \phi \|_2
\\&&\hspace{-2em}
\le
C\, \bigl(\, |x|\,\nu (x) +t \,\bigr) ^{-1-| \alpha |/2} \,\| \phi \|_2.
\end{eqnarray*}
If the integral from $0$ to $t/2$ in the preceding estimate is replaced by one from $t/2$ to $t$,
the same type of estimate yields the upper bound
$C\, \bigl(\, |x|\,\nu (x) \,\bigr) ^{-1-| \alpha |/2} \,\| \phi |S_{t/2,\, \infty }\|_2$
for this modified integral. Lemma \ref{lemmaL5.30} follows from these estimates.
\hfill $\Box $
\begin{corollary} \label{corollaryC5.20}
Let $\widetilde{ b} \in H_{\infty  },\; \delta \in (0, \infty ) ,\; \zeta \in (0,1)$ with
$\| \widetilde{ b}|S_{T, \infty }\|_{H_{T, \infty }}\le \delta \, T^{-\zeta }$
for $T \in (1, \infty )$.
Let $\phi \in L^2_n(S_{\infty })$ with $\mbox{$\mathfrak V$}  ^{( \tau )} ( \phi )|S_{\infty }=\widetilde{ b}$
(Theorem \ref{theoremT10.40}), $ R_0,\, R \in (0, \infty ) $ with $R_0<R$ and
$\overline{ \Omega  }\subset B_{R_0},$ and $ \alpha \in \mathbb{N} _0^3$ with $| \alpha |\le 1$.
Then, for $x \in B_R^c, \; t \in (0, \infty ) ,\; \epsilon \in [0,1]$,
\begin{eqnarray*}&&
| \partial _x^{\alpha }\mbox{$\mathfrak V$} ^{( \tau )} (  \phi )(x,t)|
\le
\mbox{$\mathfrak C$} \,
\bigl[\,
\bigl(\, |x|\,\nu (x) \,\bigr) ^{(-1-| \alpha |/2)}\, (1+t)^{-\zeta }
\\&&\hspace{3em}
+
\bigl(\, |x|\,\nu (x) +t \,\bigr) ^{(-1-| \alpha |/2)\,(1-\epsilon )} \,
(1+t)^{(-1-| \alpha |/2)\,\epsilon }\, \| \phi \|_2 \,\bigr]
.
\end{eqnarray*}
\end{corollary}
{\bf Proof:}
Lemma \ref{lemmaL5.30} and Theorem \ref{theoremT3.10} yield that the estimate in the corollary
holds under the additional assumption $t\ge 1$. In the case $t \in (0,1]$, we deduce from
Lemma \ref{lemmaL5.30} that
$
| \partial _x^{\alpha } \mbox{$\mathfrak V$} ^{( \tau )} (  \phi )(x,t)|
\le
\mbox{$\mathfrak C$} \,
\bigl(\, |x|\,\nu (x) \,\bigr) ^{(-1-| \alpha |/2)} \, \| \phi \|_2,
$
so we again obtain the estimate stated in the corollary.
\hfill $\Box $

\section{Main results.}
\setcounter{equation}{0}

We begin by collecting our assumptions on the right-hand side $f$
in the differential equations (\ref{10}), (\ref{180}) and (\ref{100}),
the initial data $a$ in (\ref{30}) and the Dirichlet boundary data $b$
in (\ref{20b}).

Let
$
A \in (2, \infty ),\; B \in [0,\, 3/2]$ with $A+\min\{1,B\} > 3,\; A+B\ge 7/2,\;
\varrho _0 \in (2, \infty ),\; R_1 \in (0, \infty )
$
with
$
\overline{ \Omega }\subset B_{R_1}
,\; \gamma \in L^2 \bigl(\, (0, \infty ) \,\bigr)
\cap L^{\varrho _0} \bigl(\, (0, \infty ) \,\bigr) ,\; f: \mathbb{R}^3 \times (0,  \infty ) \mapsto \mathbb{R}^3
$
measurable with
$
f|B_{R_1} \times (0, \infty ) \in L^2 \bigl(\, 0, \infty ,\, L^1( B_{R_0})^3 \,\bigr) 
$
and
\begin{eqnarray} \label{8.10a}
|f(y, \sigma )|\le \gamma ( \sigma ) \, |y| ^{-A}\, \nu (y) ^{-B}
\quad \mbox{for}\;\; 
y \in B_{R_0}^c,\; \sigma \in (0, \infty ) .
\end{eqnarray}
Moreover, suppose there are numbers
$\hat{q}_0 \in [1,\, 3/2),\; q_0 \in (1,\, 3/2),\; \overline{ q}_0 \in (1, \infty ),\;
s_0 \in (1, \infty ),\; \overline{ s}_0 \in [1, \infty ]$
such that
$3/(2\, q_0)+1/s_0 >3/2,\; 3/(2\, \overline{ q})+1/ \overline{ s}< 1,$
$$
f \in L^c \bigl(\, 0, \infty ,\, L^d( \mathbb{R}^3 )^3 \,\bigr)
\quad \mbox{for}\;\; 
(c,d) \in \{ (2,\, 3/2),\, (2,\hat{q}_0),\, (s_0,q_0)\}
$$
and
$f|B_{R_1}^c \times (0, \infty ) \in L^{\overline{ s}_0}\bigl(\, 0, \infty ,\, L^{\overline{ q}_0}(B_{R_1}^c)^3 \,\bigr) .$
Further suppose there are numbers $\delta _1 , \zeta _1\in (0, \infty ) $
such that
\begin{eqnarray} \label{8.10}
\|f| \mathbb{R}^3 \times (t/2,\, \infty )\|_{\hat{q}_0,2; \infty }
+
\|f| \mathbb{R}^3 \times (t/2,\, \infty )\|_{3/2,\, 2; \infty }
\le
\delta _1\, t^{-\zeta _1}
\quad \mbox{for}\;\;
t \in (1, \infty ).
\end{eqnarray}
Note that (\ref{8.10}) implies
$
\|f| B_{R_1} \times (t/2,\, \infty )\|_{1, 2; \infty }
\le
C(R_1) \,\delta _1\, t^{-\zeta _1}
$
for
$
t \in (1, \infty ).
$
Moreover, let $\epsilon _0 \in (0,1/2],\; p_0 \in (1,2],$
and assume that $a \in H^{1/2+\epsilon _0}_{\sigma }( \overline{ \Omega }^c)\cap L^{p_0}( \overline{ \Omega }^c)^3$.
In addition, suppose there are numbers $R_2,\, \delta _2 \in (0, \infty ) ,\; \kappa _0 \in (0,1]$
such that $\overline{ \Omega }\subset B_{R_2},\;
a| \overline{ B_{R_2}}^c \in W^{1,1}_{loc}( \overline{ B_{R_2} }^c)^3,$
\begin{eqnarray} \label{8.10b}
| \partial ^{\alpha }_ya(y)|\le \delta _2\, \bigl(\, |y|\, \nu(y) \,\bigr) ^{-1-| \alpha |/2- \kappa _0}
\quad \mbox{for}\;\;
y \in \overline{ B_{R_2}}^c,\; \alpha \in \mathbb{N} _0^3\; \mbox{with}\; | \alpha |\le 1.
\end{eqnarray}
Finally let $b \in H_{\infty }$, and suppose there are numbers $\delta _3 \in (0, \infty ) ,\; \zeta _2
\in (0,1)$ such that
\begin{eqnarray} \label{8.5}
\|b|S_{T, \infty }\|_{H_{T, \infty }}\le \delta _3\, T^{-\zeta _3}
\quad \mbox{for}\;\;
T \in (1, \infty ).
\end{eqnarray}
Now we turn to our main result on spatial and temporal pointwise decay of solutions to (\ref{10}), (\ref{20b}),
(\ref{30}).
In view of later applications in the nonlinear case, we state this result in the form of a theorem and
a corollary.
\begin{theorem} \label{theoremT8.10}
Suppose that $f,\, a$ and $b$ satisfy the assumptions listed above.
Then there is a unique function $\phi \in L^2_n(S_{\infty })$ verifying the integral
equation (\ref{T10.60.10}). Put
$u:= \mbox{$\mathfrak R$} ^{( \tau )} (f)+ \mbox{$\mathfrak I$} ^{( \tau )} (a)
+ \mbox{$\mathfrak V$} ^{( \tau )} ( \phi )| \overline{ \Omega }^c \times (0, \infty )$
(see section 4 for the definition of the preceding potential functions).
Let $\alpha \in \mathbb{N} _0^3$ with
$| \alpha |\le 1$. In addition to the assumptions above, suppose that
$3/(2\, \overline{ q})+1/ \overline{ s}< 1-| \alpha |/2$. Put
\begin{eqnarray*}&&
\varrho _1:=\min\{ \zeta _1,\, \zeta _2,\; 3/(2\, q_0)+1/s_0-3/2,\; 1/s_0,\; 3/(2\, p_0)-1/2\},
\\&&
\varrho _2:=\min\{ \zeta _1,\;
1+| \alpha |/2,\; 3/(2\, q_0)+1/s_0-1+| \alpha |/2,\; 3/(2\, p_0)+| \alpha |/2\}.
\end{eqnarray*}
Let $R \in \bigl(\, \max\{R_1,\,R_2\},\, \infty \,\bigr) .$
(The parameters $\overline{ q}_0,\, \overline{ s}_0,\, q_0,\, s_0,\, \zeta _1,\, \zeta_2,\, R_1,\, R_2,\, p_0$
were introduced at the beginning of this section.)
Then, for $x \in B_R^c,\, t \in (0, \infty ) ,\, \epsilon \in [0,1],$
\begin{eqnarray} \label{T8.10.10}&&
| \partial _x^{\alpha }u(x,t)|
\le
\mbox{$\mathfrak C$} \,
\bigl[\, \bigl(\, |x|\,\nu (x) \,\bigr) ^{-1-| \alpha |/2}\,(1+t)^{-\varrho _1}
\\&&\nonumber \hspace{3em}
+
\bigl(\, |x|\,\nu (x) \,\bigr) ^{(-1-| \alpha |/2)\,(1-\epsilon )}
\,
\bigl(\, (1+t)^{-\varrho _2}+\|f|B_{R_1}^c \times (t/2,\, \infty )\|_{ \overline{ q}_0, \overline{ s}_0; \infty }
\,\bigr)^{\epsilon } \,\bigr].
\end{eqnarray}
\end{theorem}
{\bf Proof:}
Concerning existence and uniqueness of $\phi $, we refer to Theorem \ref{theoremT10.60}.
In order to prove (\ref{T8.10.10}), take $\epsilon ,\, x,\, t$ as in the theorem. By Corollary
\ref{corollaryC6.20} and assumption (\ref{8.10}), we have
$
\|\mbox{$\mathfrak R$} ^{( \tau )} (f)|S_{T,\infty }\|_{H_{T, \infty }}
\le
\mbox{$\mathfrak C$} \, T^{-\min\{3/(2q_0)+1/s_0-3/2,\,1/s_0,\, \zeta _1\}}
$
for
$
T \in (1, \infty ).
$
Moreover, by Theorem \ref{theoremT7.10},
$
\|  \mbox{$\mathfrak I$} ^{( \tau )} (a)|S_{T,\infty }\|_{H_{T, \infty }}
\le
\mbox{$\mathfrak C$} \, T^{-3/(2p_0)+1/2}
$
equally for
$
T \in (1, \infty ).
$
These estimates, our assumptions on $b$, equation (\ref{T10.60.10}) and Corollary \ref{corollaryC5.20}
yield
\begin{eqnarray} \label{T8.10.20} &&
| \partial _x^{\alpha }\mbox{$\mathfrak V$} ^{( \tau )} ( \phi )(x,t)|
\le
\mbox{$\mathfrak C$} \,
\bigl[\, \bigl(\, |x|\,\nu (x) \,\bigr) ^{-1-| \alpha |/2}\,(1+t)^{-\zeta  _1}
\\&&\nonumber \hspace{3em}
+
\bigl(\, |x|\,\nu (x) \,\bigr) ^{(-1-| \alpha |/2)\,(1-\epsilon )}
\,
(1+t)^{(-1-| \alpha |/2) \, \epsilon } \,\bigr],
\end{eqnarray}
with $\varrho _1$ defined in the theorem.
Theorem \ref{theoremT4.10} with $R_0,\, f$ replaced by
$R_1,\;\chi_{B_{R_1}^c \times (0, \infty )}\, f$, respectively, and (\ref{8.10a}) imply
$
| \partial _x^{\alpha }\mbox{$\mathfrak R$} ^{( \tau )} \bigl(\, f| B_{R_1}^c \times (0, \infty ) \,\bigr) (y,r)|
\le
\mbox{$\mathfrak C$}\, \bigl(\, |y|\, \nu (y) \,\bigr) ^{-1-| \alpha |/2}
$
for
$
y \in B_R^c,\; r \in (0, \infty )
$.
This estimate, Corollary \ref{corollaryC6.30} with
$\varrho =2,\; ( \widetilde{ s}, \widetilde{ q})= (2,\hat{q}_0)$
and assumption (\ref{8.10}) yield
\begin{eqnarray} \label{T8.10.30}&& 
| \partial _x^{\alpha }\mbox{$\mathfrak R$} ^{( \tau )} ( f)(x,t)|
\le
\mbox{$\mathfrak C$} \,
\bigl[\, \bigl(\, |x|\,\nu (x) \,\bigr) ^{-1-| \alpha |/2}\,(1+t)^{-\zeta  _1}
\\&&\nonumber \hspace{3em}
+
\bigl(\, |x|\,\nu (x) \,\bigr) ^{(-1-| \alpha |/2)\,(1-\epsilon )}
\,
\bigl(\, (1+t)^{-\varrho _2}
+
\|f|B_{R_1}^c \times (t/2,\, \infty )\|_{ \overline{ q}_0, \overline{ s}_0; \infty } \,\bigr)^{\epsilon }
\,\bigr].
\end{eqnarray}
Corollary \ref{corollaryC7.10} provides the estimate
\begin{eqnarray} \label{T8.10.40} &&
| \partial _x^{\alpha }\mbox{$\mathfrak I$} ^{( \tau )} ( a )(x,t)|
\le
\mbox{$\mathfrak C$} \,
\bigl[\, \bigl(\, |x|\,\nu (x) \,\bigr) ^{-1-| \alpha |/2}\,(1+t)^{-2}
\\&&\nonumber \hspace{3em}
+
\bigl(\, |x|\,\nu (x) \,\bigr) ^{(-1-| \alpha |/2)\,(1-\epsilon )}
\,
(1+t)^{(-3/(2p_0)-| \alpha |/2) \, \epsilon } \,\bigr].
\end{eqnarray}
Inequality \ref{T8.10.10} follows from (\ref{T8.10.20}) - (\ref{T8.10.40}).
\hfill $\Box $

\vspace{1ex}
Theorem \ref{theoremT10.50}, \ref{theoremT10.60} and \ref{theoremT8.10}
taken together yield the following
\begin{corollary} \label{corollaryC8.10}
Under the assumptions on $f,\,a$ and $b$ listed at the beginning of this section,
there is a unique function $u \in L^2_{loc}\bigl(\, 0, \infty ,\, H^1( \overline{ \Omega }^c)^3 \,\bigr)$
satisfying (\ref{T10.50.9}) and (\ref{T10.50.10}) ($u$ velocity part of a solution to (\ref{10}), (\ref{20b}),
(\ref{30})). This function is given by
$
u= \mbox{$\mathfrak R$} ^{( \tau )} (f)+ \mbox{$\mathfrak I$} ^{( \tau )} (a)
+ \mbox{$\mathfrak V$} ^{( \tau )} ( \phi )| \overline{ \Omega }^c \times (0, \infty )
$,
with $\phi $ from (\ref{T10.60.10}), and it fulfills (\ref{T8.10.10}).
\end{corollary}
In order to exhibit the best possible rate of temporal decay inherent to inequality (\ref{T8.10.10}),
we consider the case that $f$ and $a$ are bounded functions of compact support:
\begin{corollary} \label{corollaryC8.20}
Suppose that
$f :  \mathbb{R}^3 \times (0, \infty ) \mapsto \mathbb{R}^3  $
is measurable, bounded and with compact support,
and $a \in H^1_{\sigma }( \overline{ \Omega }^c)$ also bounded and with compact support.
Let $\zeta \in [1/2,\, 1)$ and suppose that (\ref{8.5}) is satisfied with $\zeta _2=\zeta $.
Let $u$ be the solution to (\ref{T10.50.9}), (\ref{T10.50.10}). Then there is $R \in (0, \infty ) $
with $\overline{ \Omega }\subset B_R$ such that
$| \partial _x^{\alpha }u(x,t)|$ is bounded by
\begin{eqnarray*}&&
\mbox{$\mathfrak C$} \,
\bigl[\, \bigl(\, |x|\,\nu (x) \,\bigr) ^{-1-| \alpha |/2}\,(1+t)^{-\zeta }
+
\bigl(\, |x|\,\nu (x) \,\bigr) ^{(-1-| \alpha |/2)\,(1-\epsilon )}
\,
(1+t)^{(-1-| \alpha |/2) \, \epsilon } \,\bigr]
\end{eqnarray*}
for $\alpha \in \mathbb{N} _0^3$ with $| \alpha |\le 1,\; \epsilon \in [0,1],\; x \in B_R^c$
and $t \in (0, \infty ) $.
\end{corollary}
{\bf Proof:}
There are numbers $q \in (1,\, 3/2),\; p \in (1,2]$ and $s \in (1, \infty )$ so close to $1$
that $3/(2\, q)+1/s-3/2 \ge \zeta $ (in particular $3/(2\, q)+1/s -1+| \alpha |/2 \ge 1+| \alpha  |/2$),
$1/s\ge \zeta $ and $3/(2\, p)-1/2 \ge \zeta  $
(in particular $3/(2\, p) +| \alpha |/2 \ge 1+| \alpha  |/2$).
Moreover we may choose $R_0,\, T_0 \in (0, \infty ) $ such that
$\overline{ \Omega }\subset B_{R_0},\; supp(f) \subset B_{R_0} \times [0,T_0) $ and
$supp(a) \subset B_{R_0}.$
Then $f| \mathbb{R}^3 \times (t/2,\, \infty ) =0$ for $t \in [2\, T_0,\, \infty )$,
so there is $c>0$ with
\begin{eqnarray}\label{C8.20.10}
\| f| \mathbb{R}^3 \times (t/2,\, \infty )\|_{1,2; \infty }
+
\| f| \mathbb{R}^3 \times (t/2,\, \infty )\|_{3/2,\,2; \infty }
\le
c\, (1+t) ^{-2}
\quad \mbox{for}\;\;
t \in (0, \infty ) .
\end{eqnarray}
Thus the assumptions on $f$ and $a$ stated at the beginning of the chapter are verified if we suppose
that
$
R_1=R_0,\; \gamma =0,\; \hat{q}=1,\; q_0=q,\; s_0=s,\; p_0=p,
$
with $q,\,s,\,p$ as chosen above,
$\delta _1=c,\; \zeta _1=2,\; \epsilon _0=1/2,\; R_2=R_0$.
The parameters  $A,\; B,\; \varrho _0,\; \overline{ q}_0,\, \overline{ s}_0,\, \delta _2$ and $\kappa _0$
are irrelevant due to the choice of $R_0,\; \gamma $ and $T_0$;
they may be chosen as a matter of form in any way corresponding
to the specifications given at the beginning of this section. According to the assumptions
in the corollary, inequality (\ref{8.5}) holds with $\zeta _2=\zeta .$

Now take $R \in (R_0, \infty )$. Then we may conclude with Theorem \ref{theoremT8.10} that
(\ref{T8.10.10}) holds. But due to (\ref{C8.20.10}) and our choice of $q,\, s,\, p,\, \zeta _1 $
and $\zeta _2$, the parameters $\varrho _1$ and $\varrho _2$ in (\ref{T8.10.10}) equal
$\zeta $ and $1+| \alpha |/2$, respectively. This completes the proof.
\hfill $\Box $

\vspace{1ex}
Turning to the nonlinear systems (\ref{180}) and (\ref{100}), we first specify which type of
solution to the stationary problem (\ref{170}) will be considered.
\begin{theorem} \label{theoremT8.20}
Let $B \in H ^{1/2} ( \partial \Omega )^3$ with $\int_{ \partial \Omega }B \cdot n^{ (\Omega )} \, do_x =0,\;
\Psi \in L^{6/5}( \overline{ \Omega }^c)^3,\; c,\, R \in (0, \infty ) ,\; \sigma \in (4, \infty )$
such that $|\Psi(y)|\le c\, |x|^{-\sigma }$ for $x \in B_R^c$.

Then there is $U \in L^{6}( \overline{ \Omega }^c)^3\cap W^{1,1}_{loc}( \overline{ \Omega }^c)^3$
with $
\nabla U \in L^2( \overline{ \Omega }^c)^9,\; \mbox{div}\, U =0,\; U| \partial \Omega =-e_1+B,
\;
\int_{ \overline{ \Omega }^c }\bigl[\, \nabla U \cdot \nabla V + \bigl(\, \tau \, \partial _1U 
+\tau \,(U \cdot \nabla )U -\Psi \,\bigr) \cdot V \,\bigr] \, dx = 0
$
for any
$
V \in C ^{ \infty } _0( \overline{ \Omega })^3
$
with
$
\mbox{div}\, V=0.
$
In addition, there are numbers $R_3,\, c \in (0, \infty ) $ such that $\overline{ \Omega } \subset B_{R_3}$ and
\begin{eqnarray} \label{T8.20.10}
| \partial _x^{\alpha }U(x)|
\le
c\, \bigl(\, |x|\, \nu(x) \,\bigr) ^{-1-| \alpha |/2}
\quad \mbox{for}\;\;
x \in B_{R_3}^c,\; \alpha \in \mathbb{N} _0^3\; \mbox{with}\; | \alpha |\le 1.
\end{eqnarray}
\end{theorem}
{\bf Proof:}
For the first part of this theorem, up to but excluding inequality (\ref{T8.20.10}),
we refer to \cite[Theorem IX.4.1]{Ga2}, \cite[Theorem II.5.1]{Ga1}.
As for (\ref{T8.20.10}), a slightly different version of this estimate was proved in \cite{BV}
in the case $\Psi=0$;
see \cite[p. 658 above and p. 661, (3.8)]{BV}. In \cite{DeKorfu}, inequality (\ref{T8.20.10}) was deduced
from the theory in \cite{Ga2}.
\hfill $\Box $

\vspace{1ex}
In order to deal simultaneously with initial-boundary value problem (\ref{180}), (\ref{20b}), (\ref{30})
on the one hand
and (\ref{100}), (\ref{20b}), (\ref{30})
on the other,
we introduce the new parameter $\widetilde{ \tau }$,
which is to take the values $0$ or $\tau $: $\widetilde{ \tau } \in \{0,\, \tau \} .$
Then we consider a function $u$ 
with the following properties:
\begin{eqnarray} \label{8.50} &&
u \in L ^{ \infty } \bigl(\, 0, \infty ,\, H^1 ( \overline{ \Omega }^c)^3 \,\bigr) ,\;
\nabla _xu \in L^2 \bigl(\,  \overline{ \Omega }^c \times (0, \infty ) \,\bigr) ^9,
\\&&\nonumber \hspace{4em}
u(t)| \partial \Omega = b(t)\; \mbox{for}\; t \in (0, \infty ) ,\;
\mbox{div}_xu = 0,
\\[1ex] && \label{8.60}
\int_{ 0}^{ \infty }  \int_{ \overline{ \Omega}^c}
\Bigl(
- u(x,t) \cdot V (x)  \, \varphi ^{\prime} (t)
+ \nabla _x u(x,t) \cdot \nabla V (x)   \, \varphi (t)
\\&&\nonumber \hspace{1em}
+ \bigl[\,
\tau \,  \partial _1u(x,t) 
+\tau \,  \bigl( u(x,t) \cdot \nabla _x \bigr) u(x,t)
+ \widetilde{ \tau }\, \bigl(\, U(x) \cdot \nabla _x \,\bigr) u(x,t)
+ \widetilde{ \tau }\, \bigl( u(x,t) \cdot \nabla  \bigr) U(x)
\\ && \nonumber \hspace{2em}
-f(x,t)
\,\bigr]\cdot V(x) \, \varphi (t)
\Bigr)\; dx\; dt
=
\int_{ \overline{ \Omega}^c} a(x) \cdot V (x) \; dx \, \varphi (0)
\end{eqnarray}
for
$ \varphi \in C ^{ \infty } _0 \bigl(\, [0, \infty ) \,\bigr) ,\; V\in C ^{ \infty } _0( \overline{ \Omega }^c)^3$
with $ \mbox{div}\, V  =0$.
This means that $u$ is an $L^2$-strong solution to (\ref{180}), (\ref{20b}), (\ref{30}) if
$\widetilde{ \tau }=1,$ and 
an $L^2$-strong solution to (\ref{100}), (\ref{20b}), (\ref{30}) in the case $\widetilde{ \tau }=0$.
Of course, strictly speaking, $u$ is merely the velocity part of such a solution. But as usual
in this context, the preceding definition involves only this velocity part, which we call
``solution'' without any further qualification.
Results on existence of this type of solutions are due to Heywood \cite[p. 674]{Hey2},
Solonnikov \cite[Theorem 10.1, Remark 10.1]{Sol3},
and Neustupa \cite[Theorem 1]{Neu2009}, \cite[Theorem 4.1]{Neu2016}
under smallness conditions on the data.
Solutions in $L^p$-spaces with $p\ne 2$ were constructed by
Solonnikov \cite{Sol3}, Miyakawa \cite{Miya}, Shibata \cite{Sh} and Enomoto, Shibata \cite{EnSh2}.
If a solution as in (\ref{8.50}), (\ref{8.60}) is considered as given,
and if the data $f$ and $a$ decay in space as specified in (\ref{8.10a}) and (\ref{8.10b}), then it can
be shown without any smallness conditions that $u$, too, decays in space. This is the result
whose essential point is stated in (\ref{130}),
and which we now formulate in detail, choosing a version which is suitable
for what follows.
\begin{theorem} \label{theoremT8.40}
Suppose that $f,\, a$ and $b$ satisfy the assumptions listed above (some of which are not relevant
here because they are related to temporal decay). Let $U$ be the function from
Theorem \ref{theoremT8.20}, and assume that $u$ satisfies (\ref{8.50})
and (\ref{8.60}). Put
$
F(x,t):=-\tau \,  \bigl( u(x,t) \cdot \nabla _x \bigr) u(x,t),
\;
G(x,t)
:=
-\widetilde{ \tau }\, \bigl(\, U(x) \cdot \nabla _x \,\bigr) u(x,t)
- \widetilde{ \tau }\, \bigl( u(x,t) \cdot \nabla  \bigr) U(x)
$
for
$
x \in \overline{ \Omega  }^c,\; t \in (0, \infty ) .
$
Then there are constants $c,\, R_4 \in (0, \infty ) $
such that $\overline{ \Omega }\subset B_{R_4},\;\;
| \partial _x^{\alpha }\mbox{$\mathfrak R$} ^{( \tau )} (F+G)(x,t)|
\le
c \,
\bigl(\, |x|\,\nu (x) \,\bigr) ^{-1-| \alpha |/2}
$
and
$
| \partial _x^{\alpha }u(x,t)|
\le
c \,
\bigl(\, |x|\,\nu (x) \,\bigr) ^{-1-| \alpha |/2}
$
for
$
x \in B_{R_4}^c,\; t \in (0, \infty ) ,\; \alpha \in \mathbb{N} _0^3
$
with
$
| \alpha |\le 1.
$
\end{theorem}
{\bf Proof:}
First consider the case $\widetilde{ \tau }=1$. We refer to \cite[Corollary 3.5, Theorem 4.6, 4.8]{DeSIMA2}
for the first inequality, which actually is an intermediate result in the proof of the second
(\cite[Theorem 1.2]{DeSIMA2}). Note that in \cite{DeSIMA2}, the function $f$ is supposed to fulfill
the relation $f|B_{R_1}\times (0, \infty ) \in L^2 \bigl(\, B_{R_1}\times (0, \infty ) \,\bigr) ^3,$
instead of the weaker condition
$f|B_{R_1}\times (0, \infty ) \in L^2 \bigl(\, 0, \infty ,\, L^1( B_{R_1})^3 \,\bigr)$
required here.
However, the former condition enters into \cite{DeSIMA2} only via
\cite[Theorem 2.16]{DeSIMA2},
which restates \cite[Theorem 3.1]{DeDCDS-A} and is reproduced as Theorem \ref{theoremT4.10} in the work
at hand. But the conclusion of this theorem remains valid if its assumptions are modified as indicated
above. In fact, Theorem \ref{theoremT4.10} as it stands yields (\ref{130})
in the case that $f$ is replaced by
$\chi_{B_{R_0}^c \times (0, \infty )}\, f$,
and Lemma \ref{lemmaL6.40} implies (\ref{130}) if
$f  \in L^2 \bigl(\, 0, \infty ,\, L^1( B_{R_0})^3 \,\bigr) $. (In these two references,
the parameter $R_0$ is in the role of $R_1$ here.)
Combining these two results, we end up with inequality (\ref{130}) under the modified assumptions
on $f$.
We further remark that in \cite{DeSIMA2}, we supposed $\Psi =0$ and $B=0$
in the analogue of Theorem \ref{theoremT8.20} (\cite[Theorem 1.1]{DeSIMA2}). However, these
relations are not used anywhere in \cite{DeSIMA2};
only the conclusions on $U$ stated in Theorem \ref{theoremT8.20} are relevant.
Concerning the case $\widetilde{ \tau }=0$, we indicate that all the proofs in \cite{DeSIMA2}
also carry through, some of them in a much simpler form, if the function $U$ in that reference
vanishes. Therefore the preceding arguments remain valid if $\widetilde{ \tau }=0$.
\hfill $\Box $

\vspace{1ex}
In the ensuing lemmas, we collect some integrability properties of the functions $F$ (nonlinearity)
and $G$ defined in Theorem \ref{theoremT8.40}.
\begin{lemma} \label{lemmaL8.15}\label{lemmaL8.20}
Let $u$ be a function satisfying (\ref{8.50}), and define $F$ and $G$ as in Theorem \ref{theoremT8.40}.
Let $H \in \{F,\, G\}$.
Then $H \in L^2 \bigl(\, 0, \infty ,\, L^{\varrho }( \overline{ \Omega }^c)^3 \,\bigr) $
and
$
\|H| \overline{ \Omega }^c \times (t/2,\, \infty ) \|_{\varrho ,2; \infty }
\le
\mbox{$\mathfrak C$} \, 
\| \nabla _xu| \overline{ \Omega }^c \times (t/2,\, \infty ) \|_{2 }
$
for
$
t \in (0, \infty ) ,\; \varrho \in [1,\, 3/2]
$
in the case of $H=F$, and $\varrho \in [11/10,\, 3/2]$ if $H=G$.
\end{lemma}
{\bf Proof:}
First suppose that $H=F$. Then the statements of the lemma are well known.
By the proof of \cite[Lemma 3.2, estimate (3.1) with $s_0=2,\; r_0=6$]{DeSIMA2}
for example, we have
$
\|F(r)\|_{\varrho }\le \mbox{$\mathfrak C$} \, (\|u(r)\|_2+\|u(r)\|_6)\,\| \nabla _xu(r)\|_2
$
for
$r \in (0, \infty ) ,\; \varrho =1$ and $\varrho =3/2.$
The claims in Lemma \ref{lemmaL8.15} with respect to $F$ follow with (\ref{8.50}) and
Theorem \ref{theoremT1.101}.

Now consider the case $H=G$. Let $r \in (0, \infty ) $.
By H\"older's inequality and Theorem \ref{theoremT1.101}, we get
$
\| (U \cdot \nabla _x)u(r)\|_{3/2}
\le
C\,\|U\|_6\, \| \nabla _xu(r)\|_2
$
and
$
\| \bigl(\, u(r) \cdot \nabla \,\bigr) U\|_{3/2}\le C\,\|u(r)\|_6\, \| \nabla U\|_2
\le
\mbox{$\mathfrak C$} (U)\, \| \nabla _xu(r)\|_2,
$
so that
$
\|G(r)\|_{3/2}\le \mbox{$\mathfrak C$}(U) \,\| \nabla _xu(r)\|_2.
$
Moreover, as in the proof of \cite[Lemma 3.2]{DeSIMA2}, 
inequality (\ref{T8.20.10}) and Corollary \ref{corollaryC1.10} imply
$U \in L^{22/9}( \overline{ \Omega }^c)^3,\; \nabla U \in L^{66/49}( \overline{ \Omega }^c)^9$,
so by H\"older's inequality,
$
\| (U \cdot \nabla _x)u(r)\|_{11/10}
\le
C\,\|U\|_{22/9}\, \| \nabla _xu(r)\|_2
\le
\mbox{$\mathfrak C$} (U)\,\| \nabla _xu(r)\|_2
$
and
$
\| \bigl(\, u(r) \cdot \nabla \,\bigr) U\|_{11/10}\le \|u(r)\|_6\, \| \nabla U\|_{66/49}
\le
\mbox{$\mathfrak C$} (U)\, \| \nabla _xu(r)\|_2.
$
Therefore
$
\|G(r)\|_{11/10}
$
is bounded by
$
\mbox{$\mathfrak C$}(U) \,\| \nabla _xu(r)\|_2.
$
The statements in the lemma relating to $G$ now follow by interpolation.
Note that in \cite[(3.4)]{DeSIMA2}, it should read $\|u\|_{6,2;T_0}$
instead of $\|u\|_{2,6;T_0}$.
\hfill $\Box $
\begin{corollary} \label{corollaryC8.30}
In the situation of Theorem \ref{theoremT8.40}, let $\overline{ R} \in  [R_4, \infty )$,
and take $R \in ( \overline{ R}, \infty ).$
Then the estimate
$
| \partial _x^{\alpha }\mbox{$\mathfrak R$} ^{( \tau )} \bigl(\, F+G|B_{\overline{ R}}^c
\times (0, \infty ) \,\bigr) (x,t)|
\le
\mbox{$\mathfrak C$} \, \bigl(\, |x|\, \nu(x) \,\bigr) ^{-1-| \alpha |/2}
$
holds for
$
x \in B_R^c,\; t \in (0, \infty ) $ and $ \alpha \in \mathbb{N} _0^3
$
with
$
| \alpha |\le 1.
$

\end{corollary}
{\bf Proof:}
By Lemma \ref{lemmaL8.20}, we have
$F+G \in L^2 \bigl(\, 0, \infty ,\, L^{11/10}( \overline{ \Omega } ^c)^3 \,\bigr) ,$
so Lemma \ref{lemmaL6.40} with $\epsilon =0$ and
$F+G|B_{\overline{ R}} \times (0, \infty )$
in the place of $f$ yields for $x,\,t,\, \alpha $ as above that
$
| \partial _x^{\alpha }\mbox{$\mathfrak R$} ^{( \tau )} \bigl(\, F+G|B_{\overline{ R}}
\times (0, \infty ) \,\bigr) (x,t)|
\le
\mbox{$\mathfrak C$} \, \bigl(\, |x|\, \nu(x) \,\bigr) ^{-1-| \alpha |/2}$.
Thus the corollary follows from the first inequality in Theorem \ref{theoremT8.40}.
\hfill $\Box $
\begin{lemma} \label{lemmaL8.30}
Let $f,\,a,\, b,\, U,\, u, \, F,\; G$ be given as in Theorem \ref{theoremT8.40},
and take $p \in [2,6],\;
\varrho \in [p, \infty )$. Put $\overline{ R}:=\max\{R_3,\, R_4\}$,
with $R_3$ from Theorem \ref{theoremT8.20} and $R_4$ from Theorem \ref{theoremT8.40}.
Then, for $H \in \{F,\, G\}$, the function $H|B_{\overline{ R}}^c \times (0, \infty )$ belongs to
$L^{\varrho }\bigl(\, 0, \infty ,\, L^p( B_{\overline{ R}}^c)^3 \,\bigr) $, and
$
\|H|B_{\overline{ R}}^c \times (t, \infty )\|_{p, \varrho ; \infty }\le \mbox{$\mathfrak C$} \,
\| \nabla _xu| \overline{ \Omega }^c \times (t, \infty )\|_{2 }^{2/ \varrho }
\quad \mbox{for}\;\; t \in (0, \infty ) .
$
\end{lemma}
{\bf Proof:}
We note that $(1/p-1/6) ^{-1} \in (4/3,\, \infty ],$
so $\| \nabla U| B_{\overline{ R}}^c\|_{(1/p-1/6 ) ^{-1} }< \infty $
by (\ref{T8.20.10}) and Corollary \ref{corollaryC1.10}.
Inequality (\ref{T8.20.10}) further yields that
$\| U| B_{\overline{ R}}^c\|_{\infty  }< \infty $.
Define
$
G_1(x,r):= \bigl(\, u(x,r) \cdot \nabla \,\bigr) U(x)
$
for $x \in \overline{ \Omega }^c,\; r \in (0, \infty ) $.
Then we may conclude with Theorem \ref{theoremT1.101} that
$
\|G_1(r)|B_{\overline{ R}}^c\|_p
\le
C\, \|u(r)\|_6\,\| \nabla U|B_{\overline{ R}}^c\|_{(1/p-1/6) ^{-1} }
\le
\mbox{$\mathfrak C$} (U)\, \| \nabla u(r)\|_2
$
for
$
r \in (0, \infty ) .
$
Let $t \in [0, \infty ) .$
Since $\varrho \ge p\ge 2$,
the preceding estimate implies
\begin{eqnarray*} &&
\|G_1|B_{\overline{ R}}^c \times (t, \infty )\|_{p, \varrho ; \infty }
\le
\mbox{$\mathfrak C$} (U)\,
\Bigl( \int_{ t}^{ \infty } \| \nabla _xu(r)\|_2^{\varrho }\, dr \Bigr) ^{1/ \varrho }
\\&&
\le
\mbox{$\mathfrak C$}\,
\| \nabla _xu\|_{2, \infty ; \infty  }^{1-2/ \varrho }\, 
\Bigl( \int_{ t}^{ \infty } \| \nabla _xu(r)\|_2^{2}\, dr \Bigr) ^{1/ \varrho }
\le
\mbox{$\mathfrak C$}\,
\| \nabla _xu\|_{2, \infty ; \infty  }^{1-2/ \varrho }\, 
\| \nabla _xu | \overline{ \Omega }^c \times (t, \infty )\|_2^{2/ \varrho }.
\end{eqnarray*}
Put
$
G_2(x,r):= \bigl(\, U(x) \cdot \nabla _x\,\bigr) u(x,r)
$
for
$
x \in \overline{ \Omega }^c,\; r \in (0, \infty )
$.
Then Theorem \ref{theoremT8.40} yields that
$$
\|G_2(r)|B_{\overline{ R}}^c\|_p
\le
C\,\| U|B_{\overline{ R}}^c\|_{\infty  }\, \| \nabla _xu(r)|B_{\overline{ R}}^c\|_p
\le
\mbox{$\mathfrak C$} (U)\, \Bigl( \int_{ B_{\overline{ R}}^c} | \nabla u(x,r)|^2\, |x|^{-(3/2)\,(p-2)}\, dx \Bigr) ^{1/p}
$$
for $r \in (0, \infty ) $,
so
$
\|G_2(r)|B_{\overline{ R}}^c\|_p
\le
\mbox{$\mathfrak C$} \, \| \nabla _xu(r)\|_2^{2/p}.
$
Since $\varrho \ge p$, hence $2\, \varrho /p\ge 2,$
we thus get as in the estimate of $G_1$ that
\begin{eqnarray*} &&
\|G_2|B_{\overline{ R}}^c \times (t, \infty )\|_{p, \varrho ; \infty }
\le
\mbox{$\mathfrak C$}\,
\| \nabla _xu\|_{2, \infty ; \infty  }^{2(1/p-1/ \varrho )}\, 
\| \nabla _xu | \overline{ \Omega }^c \times (t, \infty )\|_2^{2/ \varrho }.
\end{eqnarray*}
But $\| \nabla _xu\|_{2, \infty ; \infty  }< \infty $ by (\ref{8.50}),
so the claims of the lemma referring to the case $H=G$ are proved.
Since by Theorem \ref{theoremT8.40} and Corollary \ref{corollaryC1.10}, there is
$c>0$ with
$
\|u(r)|B_{\overline{ R}}^c\|_{\infty }\le c,\;\| \nabla _xu(r)|B_{\overline{ R}}^c\|_{(1/p-1/6) ^{-1} }\le c
$
for
$r \in (0, \infty ) $,
the case $H=F$ may be handled in the same way.
\hfill $\Box $

\vspace{1ex}
For a function $u$ which satisfies not only (\ref{8.50}) and (\ref{8.60}), but also
the pointwise decay estimate of
$\| \nabla _xu(t)\|_2$ in (\ref{8.30}),
the following relations hold for the nonlinearity $(u \cdot \nabla _x)u$, that is, for the function $F$
from Theorem \ref{theoremT8.40}.
\begin{lemma} \label{lemmaL8.10}
In the situation of Theorem \ref{theoremT8.40}, suppose in addition that
that $u$ satisfies (\ref{8.30}). Then
\begin{eqnarray} \label{L8.10.10}&&\hspace{-4em}
\| \nabla _xu(t)\|_2\le \mbox{$\mathfrak C$} \, (1+t)^{-\kappa _1}
\quad \mbox{for}\;\;
t >0 ;
\\&& \label{L8.10.20}\hspace{-4em}
\|F| \overline{ \Omega }^c \times (t, \infty )\|_{q ,2; \infty }
\le
\mbox{$\mathfrak C$} \,
(1+t)^{
3 \kappa _1 (1-1/ q )
}
\quad \mbox{for}\;\;
t>0 ,\; q \in [1,\, 3/2];
\\&&\label{L8.10.30}\hspace{-4em}
\|F| B_{R_4}^c \times (t, \infty )\|_{q , \infty ; \infty }
\le
\mbox{$\mathfrak C$} \, (1+t)^{-2 \kappa _1/q}
\quad \mbox{for}\;\;
t >0 ,\; q \in [2, \infty ).
\end{eqnarray}
\end{lemma}
{\bf Proof:}
Since $u \in L ^{ \infty } \bigl(\, 0, \infty ,\,H^1( \overline{ \Omega  }^c)^3 \,\bigr) $
by (\ref{8.50}),
the estimate in (\ref{L8.10.10}) follows immediately with (\ref{8.30}). Moreover,
with Theorem \ref{theoremT1.101} and (\ref{L8.10.10}), we find that
$
\|F(r)\|_{3/2}
\le
C\, \|u(r)\|_6\, \| \nabla _xu(r)\|_2
\le
\mbox{$\mathfrak C$} \, \| \nabla _xu(r)\|_2^2
\le
\mbox{$\mathfrak C$} \, \| \nabla _xu(r)\|_2\, (1+r)^{-\kappa _1}
$
for
$
r \in (0, \infty ).
$
This implies (\ref{L8.10.20}) in the case $q =3/2$.
From
(\ref{8.50})
we get
$
\|F(r)\|_{1}
\le
C\, \|u(r)\|_2\, \| \nabla _xu(r)\|_2
\le
\mbox{$\mathfrak C$} \, \| \nabla _xu(r)\|_2
\le \mbox{$\mathfrak C$} 
$
for
$
r \in (0, \infty ),
$
so that (\ref{L8.10.20}) holds for $q =1$ as well.
That latter estimate for $q \in (1,\,3/2)$ now follows by interpolation.
Let $q \in [2, \infty )$.
Then for $r$ as before, by the second estimate in Theorem \ref{theoremT8.40} and
(\ref{L8.10.10}),
$
\|F(r)|B_{R_4}^c\|_{q}
\le
\mbox{$\mathfrak C$} \, \bigl(\, \int_{ B_{R_4}^c}|x|^{-q-(3/2)(q-2)}\, | \nabla _xu(x,r)|^2\, dx \,\bigr) ^{1/q}
\le
\mbox{$\mathfrak C$} \, \| \nabla _xu(r)\|_2^{2/q}
\le
\mbox{$\mathfrak C$} \, (1+r)^{-2 \kappa _1/q}.
$
This proves (\ref{L8.10.30}).
\hfill $\Box $

\vspace{1ex}
Now we are in a position to prove our main results about the asymptotics of solutions to the
nonlinear problem (\ref{180}), (\ref{20b}), (\ref{30}) and (\ref{100}), (\ref{20b}), (\ref{30}).
\begin{theorem} \label{theoremT8.50}
Suppose that $f,\, a$ and $b$ satisfy the assumptions listed at the beginning of this section.
Let $\alpha \in \mathbb{N} _0^3$ with $| \alpha |\le 1$, and let the parameters $\overline{ q}_0,\,
\overline{ s}_0$ introduced at the beginning of this section verify the additional condition
$3/(2\, \overline{ q}_0)+1/ \overline{ s}_0<1-| \alpha |/2$. Take $U$ as in Theorem \ref{theoremT8.20},
let $\widetilde{ \tau } \in \{0,\, 1\}$, and suppose that $u$ is given as in (\ref{8.50}), (\ref{8.60}),
that is, $u$ is a $L^2$-strong solution to
(\ref{180}), (\ref{20b}), (\ref{30}) if $\widetilde{ \tau }=1,$
and to (\ref{100}), (\ref{20b}), (\ref{30}) in the case $\widetilde{ \tau }=0.$
Let $R \in (\max\{R_1,\, ...,\, R_4\},\; \infty )$, with $R_3$ from Theorem \ref{theoremT8.20}, $R_4$
from Theorem \ref{theoremT8.40}, and $R_1,\,R_2$ fixed at the beginning of this section.

Then there is a bounded function $X  :(0, \infty ) \mapsto (0, \infty ) $
with $X  (t)\downarrow 0$ for $t\to \infty $ such that
\begin{eqnarray} \label{T8.50.10}
| \partial _x^{\alpha }u(x,t)|
\le
\bigl(\, |x|\,\nu (x) \,\bigr) ^{(-1-| \alpha |/2)\,(1-\epsilon )}\, X  (t)^{\epsilon }
\quad \mbox{for}\;\;
x \in B_R^c,\; t \in (0, \infty ) ,\; \epsilon \in [0,1].
\end{eqnarray}
If $\widetilde{ \tau }=0$ and $u$ additionally satisfies (\ref{8.30}) with $\kappa _1>0$,
and if
$
q_1 ,\,\hat{q}_1 \in (1,\, 3/2)
,\;
\overline{ q}_1 \in [2, \infty )$ with $3/(2-| \alpha |)< \overline{ q}_1$,
then
\begin{eqnarray} \label{T8.50.15}&&
| \partial _x^{\alpha }u(x,t)|
\le
\mbox{$\mathfrak C$} \,
\bigl[\, \bigl(\, |x|\,\nu (x) \,\bigr) ^{(-1-| \alpha |/2)}\,
(1+t)^{-\min\{ \varrho _1,\; 3/(2\, q_1)-1,\;
3 \kappa _1\,(1-1/\hat{q}_1)
,\, \kappa _1\}}
\\&&\nonumber \hspace{2em}
+
\bigl(\, |x|\,\nu (x) \,\bigr) ^{(-1-| \alpha |/2)\,(1-\epsilon )}\,
\\&&\nonumber \hspace{4em}
\cdot \bigl(\, (1+t)^{-\min\{ \varrho _2,\; 3/(2\, q_1)-1/2+| \alpha |/2,\;  2 \kappa _1/ \overline{ q}_2\}}
+\|f|B_{R_1}^c \times (t/2,\, \infty )\|_{\overline{ q}_0, \overline{ s}_0; \infty } \,\bigr) ^{\epsilon }
\,\bigr]
\end{eqnarray}
for
$
x,\,t,\, \epsilon 
$
as in (\ref{T8.50.10}), with $\varrho _1,\, \varrho _2$ from Theorem \ref{theoremT8.10}.

Assume in addition that $f$ and $a$ are bounded and with compact support, the parameter $\zeta _2$ from
(\ref{8.5}) belongs to $[1/2,\, 1)$, and $\delta >0$ is suffiently small. Then
\begin{eqnarray} \label{T8.50.16}&&
| \partial _x^{\alpha }u(x,t)|
\le
\mbox{$\mathfrak C$} \,
\bigl[\, \bigl(\, |x|\,\nu (x) \,\bigr) ^{(-1-| \alpha |/2)}\,
(1+t)^{-\min\{ 1/2,\,  \kappa _1 \}+\delta }
\\&&\nonumber \hspace{2em}
+
\bigl(\, |x|\,\nu (x) \,\bigr) ^{(-1-| \alpha |/2)\,(1-\epsilon )}\,
(1+t)^{-(\min\{1+| \alpha |/2 ,\; k( \alpha )\}+\delta )\, \epsilon }
\,\bigr],
\end{eqnarray}
with $k( \alpha ):=\kappa _1$ in the case $\alpha =0$,
and $k( \alpha ):=2\, \kappa _1 /3$ if $| \alpha |=1$.
\end{theorem}
{\bf Proof:}
Define $F$ and $G$ as in the proof of Theorem \ref{theoremT8.40}.
By our assumptions on $f$ and by Lemma \ref{lemmaL8.15}, we have
$f+F+G \in L^2 \bigl(\, 0, \infty ,\, L^p( \overline{ \Omega }^c)^3 \,\bigr) $
for
$p \in [\max\{11/10,\, \hat{q}_0\},\: 3/2].$
From (\ref{8.50}), (\ref{8.60}), we obtain that $u$ satisfies
(\ref{T10.50.10}) with
$f+F+G$ in the place of $f$,
as well as (\ref{T10.50.9}). Thus Theorem \ref{theoremT10.60} yields there is a unique function
$\phi \in L^2_n(S_{\infty })$ such that
\begin{eqnarray} \label{T8.50.20}
\mbox{$\mathfrak V$} ^{( \tau )} ( \phi )|S_{\infty }
=
-\mbox{$\mathfrak R$} ^{( \tau )} (f+F+G)- \mbox{$\mathfrak I$} ^{( \tau )} (a)+b,
\end{eqnarray}
and
\begin{eqnarray} \label{T8.50.30}
u = \mbox{$\mathfrak R$} ^{( \tau )} (f+F+G)+\mbox{$\mathfrak I$} ^{( \tau )} (a)+\mbox{$\mathfrak V$} ^{( \tau )}
( \phi )| \overline{ \Omega }^c \times (0, \infty ).
\end{eqnarray}
On the other hand, since $f$ belongs to
$L^2 \bigl(\,  0, \infty ,\, L^{q} ( \overline{ \Omega }^c)^3 \,\bigr) $
for $q \in \{\hat{q}_0,\, 3/2\}$,
we know from Theorem \ref{theoremT10.5} that
$\mbox{$\mathfrak R$} ^{( \tau )} (f)|S_{\infty }\in H_{\infty }.$
Moreover
$\mbox{$\mathfrak I$} ^{( \tau )} (a)|S_{\infty }\in H_{\infty } $
by Theorem \ref{theoremT10.6}, and $b \in H_{\infty }$ by assumption.
Therefore Theorem \ref{theoremT10.40} states there is a unique function
$\phi ^{(1)} \in L^2_n(S_{\infty })$
such that
\begin{eqnarray} \label{T8.50.40}
\mbox{$\mathfrak V$} ^{( \tau )} ( \phi ^{(1)}  )|S_{\infty }
=
-\mbox{$\mathfrak R$} ^{( \tau )} (f)- \mbox{$\mathfrak I$} ^{( \tau )} (a)+b.
\end{eqnarray}
According to Lemma \ref{lemmaL8.15}, the function $F+G$ belongs to
$L^2 \bigl(\, 0, \infty ,\, L^p( \overline{ \Omega }^c)^3 \,\bigr) $
for
$p \in [11/10,\: 3/2]$
so that
$\mbox{$\mathfrak R$} ^{( \tau )} (F+G)|S_{\infty }\in H_{\infty }$
(Theorem \ref{theoremT10.5}). Therefore, once more by Theorem \ref{theoremT10.40},
there is 
$\phi ^{(2)} \in L^2_n(S_{\infty })$
such that
\begin{eqnarray} \label{T8.50.50}
\mbox{$\mathfrak V$} ^{( \tau )} ( \phi ^{(2)}  )|S_{\infty }
=
-\mbox{$\mathfrak R$} ^{( \tau )} (F+G)|S_{\infty }.
\end{eqnarray}
In view of the uniqueness statement in Theorem \ref{theoremT10.40}, and because of
(\ref{T8.50.20}), (\ref{T8.50.40}) and (\ref{T8.50.50}), we may conclude that
$\phi = \phi ^{(1)} +\phi ^{(2)} ,$
so by (\ref{T8.50.30})
$u = u ^{(1)} +u ^{(2)} $,
where
$
u ^{(1)}
:=
\mbox{$\mathfrak R$} ^{( \tau )} (f)+\mbox{$\mathfrak I$} ^{( \tau )} (a)+\mbox{$\mathfrak V$} ^{( \tau )}
( \phi ^{(1)} )| \overline{ \Omega }^c \times (0, \infty )
$
and
$
u ^{(2)}
:=
\mbox{$\mathfrak R$} ^{( \tau )} (F+G)+\mbox{$\mathfrak V$} ^{( \tau )}( \phi ^{(2)}  )
| \overline{ \Omega }^c \times (0, \infty ).
$
By the definition of $u ^{(1)} $, inequality (\ref{T8.10.10}) holds with $u ^{(1)} $ in the place of $u$, hence
\begin{eqnarray} \label{T8.50.55}
| \partial _x^{\alpha }u ^{(1)} (x,t)|
\le
\mbox{$\mathfrak C$} \,
\bigl(\, |x|\,\nu(x) \,\bigr) ^{(-1-| \alpha |/2)\,(1-\epsilon )}\, X  _1(t)^{\epsilon }
\end{eqnarray}
for $x,\, t,\, \epsilon $
as in (\ref{T8.50.10}), with
\begin{eqnarray*}
X  _1(t)
:=
(1+t)^{-\varrho _1}+(1+t)^{-\varrho _2}+\|f|B_{R_1}^c \times (t/2,\, \infty )\|_{\overline{ q}_0, \overline{ s}_0; \infty }.
\end{eqnarray*}
Turning to $u ^{(2)} $, we may deduce from Lemma \ref{lemmaL5.30} that
\begin{eqnarray} \label{T8.50.100}
| \partial _x^{\alpha }\mbox{$\mathfrak V$} ^{( \tau )} ( \phi ^{(2)}  )(x,t)|
\le
\mbox{$\mathfrak C$} \,
\bigl(\, |x|\,\nu(x) \,\bigr) ^{(-1-| \alpha |/2)\,(1-\epsilon )}\, X  _2(t)^{\epsilon }
\end{eqnarray}
for $x,\, \epsilon $
as in (\ref{T8.50.10}) and $t \in [1, \infty )$, where $X  _2(t)$ is defined by
\begin{eqnarray*}
X  _2(t)
:=
(1+t)^{-1-| \alpha |/2}+\| \phi ^{(2)} |S_{t/2,\, \infty }\|_2.
\end{eqnarray*}
If $t \in (0,1),$ we argue as in the proof of Corollary \ref{corollaryC5.20}, referring to Lemma
\ref{lemmaL5.30} to obtain
$
|\mbox{$\mathfrak V$} ^{( \tau )} ( \phi ^{(2)}  )(x,t)|
\le
\mbox{$\mathfrak C$} \,
\bigl(\, |x|\,\nu(x) \,\bigr) ^{-1-| \alpha |/2}
$
for $x,\, t$ as in (\ref{T8.50.10}).
Thus we obtain (\ref{T8.50.100}) again, this time for $t \in (0,1)$.
Put $\overline{ R}:=\max\{R_3,\, R_4\}.$
The term
%
$\| F+G|B_{\overline{ R}}^c \times (t/2,\,  \infty )\|_{6,8; \infty }$
is bounded by
$
\mbox{$\mathfrak C$} \,
\| \nabla _x u| \overline{ \Omega }^c \times (t/2,\, \infty )\|_2^{1/4}
$
for $t \in (0, \infty ) $
according to Lemma \ref{lemmaL8.30}.
Moreover Lemma \ref{lemmaL8.15} yields that
$
F+G
\in 
L^2 \bigl(\, 0, \infty ,\;L^p(\overline{ \Omega }^c)^3 \,\bigr) $
and
$
\| F+G|\overline{\Omega }^c \times (t/2,\,  \infty )\|_{p,2; \infty }
\le
\mbox{$\mathfrak C$} \,
\| \nabla _x u| \overline{ \Omega }^c \times (t/2,\, \infty )\|_2
$
for $t \in (0, \infty ),\; p \in [11/10,\, 3/2] $.
From this latter estimate we may conclude that
$
\| F+G|(\overline{\Omega }^c\cap B_{\overline{ R}}) \times (t/2,\,  \infty )\|_{1,2; \infty }
\le
\mbox{$\mathfrak C$} \,
\| \nabla _x u| \overline{ \Omega }^c \times (t/2,\, \infty )\|_2
$.
By Corollary \ref{corollaryC8.30}, we know there is $\mbox{$\mathfrak D$} _0>0$
with
$
| \partial _x^{\alpha } \mbox{$\mathfrak R$} ^{( \tau )} \bigl(\, F+G|B_{\overline{ R}}^c \times (0, \infty ) \,\bigr)
(x,t)|
\le
\mbox{$\mathfrak D$} _0\,
\bigl(\, |x|\,\nu(x) \,\bigr) ^{-1-| \alpha |/2}
$
for
$x \in B_R^c,\; t \in (0, \infty ) $,
hence (\ref{C6.30.10}) is valid with $F+G,\, \overline{ R}$ in the place of $f$ and $R_0$,
respectively.
Due to these estimates, we may apply Corollary \ref{corollaryC6.30} with
$\varrho =2,\; q=4/3,\; s=2,\; \widetilde{ q}=11/10,\; \widetilde{ s}=2,\; \overline{ q}=6,\;
\overline{ s}= 8
$
and $R_0,\, f$ replaced by $\overline{ R}$ and $F+G$, respectively, to obtain
\begin{eqnarray} \label{T8.50.120}
| \partial _x^{\alpha } \mbox{$\mathfrak R$} ^{( \tau )} ( F+G)(x,t)|
\le
\mbox{$\mathfrak C$} \,
\bigl(\, |x|\,\nu(x) \,\bigr) ^{(-1-| \alpha |/2)\,(1-\epsilon )}
\, X  _3(t)^{\epsilon }
\end{eqnarray}
for
$x,\, t,\, \epsilon $
as in (\ref{T8.50.10}), with
$
X  _3(t)
:=
(1+t)^{-7/8}
+
\sum_{r \in \{1,\, 1/4\}}
\| \nabla _x u| \overline{ \Omega }^c \times (t/2,\, \infty )\|_2^r.
$
As a consequence of (\ref{T8.50.100}) and (\ref{T8.50.120}) and the definition of $u ^{(2)} $,
we arrive at the estimate
\begin{eqnarray*} 
| \partial _x^{\alpha } u ^{(2)} (x,t)|
\le
\mbox{$\mathfrak C$} \,
\bigl(\, |x|\,\nu(x) \,\bigr) ^{(-1-| \alpha |/2)\,(1-\epsilon )}
\, \bigl(\, X  _2(t)+X  _3(t) \,\bigr) ^{\epsilon }
\end{eqnarray*}
for
$x,\, t,\, \epsilon $
as in (\ref{T8.50.10}).
Inequality (\ref{T8.50.10}) now follows with (\ref{T8.50.55}) and because
$u=u ^{(1)} +u ^{(2)} $.

Now suppose that $\widetilde{\tau }=0$ so that $G=0$.
Further suppose that $u$ satisfies (\ref{8.30}) with $\kappa _1 >0$.
Let $q_1,\, \hat{q}_1$ and $\overline{ q}_1$ be chosen as in the theorem.
Note that
$F \in L^2 \bigl(\, 0, \infty ,\,L^{q_1}( \overline{ \Omega }^c)^3 \,\bigr) $
according to Lemma \ref{lemmaL8.15}
By applying (\ref{L8.10.20}) with $q=\hat{q}_1$ and $q=3/2$, and
Corollary \ref{corollaryC6.20} with $q=q_1,\; s=2,\; \hat{q}=\hat{q}_1$,
we get
$
\| \mbox{$\mathfrak R$} ^{( \tau )} (F)|S_{T, \infty }\|_{H_{T, \infty }}
\le
\mbox{$\mathfrak C$} \,
T^{-\min\{3/(2q_1)-1,\; 3 \kappa _1\,(1-1/\hat{q}_1),\; \kappa _1\}}
$
for
$
T \in (1, \infty ).
$
This inequality, (\ref{T8.50.50}) and Corollary \ref{corollaryC5.20} imply
\begin{eqnarray} \label{T8.50.70}&&\hspace{-3em}
| \partial _x^{\alpha }\mbox{$\mathfrak V$} ^{( \tau )} ( \phi ^{(2)}  )(x,t)|
\le
\mbox{$\mathfrak C$} \,
\bigl[\, \bigl(\, |x|\,\nu(x) \,\bigr) ^{-1-| \alpha |/2}\,
(1+t)^{-\min\{3/(2q_1)-1,\; 3 \kappa _1\,(1-1/\hat{q}_1),\; \kappa _1\}}
\\&&\nonumber \hspace{2em}
+
\bigl(\, |x|\,\nu(x) \,\bigr) ^{(-1-| \alpha |/2)\,(1-\epsilon )}\,
(1+t)^{(-1-| \alpha |/2)\, \epsilon } \,\bigr]
\end{eqnarray}
for $x,\,t,\, \epsilon $
as in (\ref{T8.50.10}).
Next we want to apply Corollary \ref{corollaryC6.30} to $\mbox{$\mathfrak R$} ^{( \tau )} (F)$,
with $f,\, R_0$ replaced by $F$ and $\overline{ R}$, respectively.
As indicated further above,
inequality (\ref{C6.30.10}) holds for the preceding choice of $f$ and $R_0$ in that latter estimate.
Recall that by Lemma \ref{lemmaL8.15},
$
F
\in 
L^2 \bigl(\, 0, \infty ,\; L^p(\overline{ \Omega }^c)^3 \,\bigr)
$
for any $p \in [1,\,3/2]$.
On the other hand, it follows with (\ref{L8.10.20}) that
$
\|F|( \overline{ \Omega }^c\cap B_{\overline{ R}})\times (t/2,\, \infty )\|_{1,2;\infty }
\le
C( \overline{ R})\,
\|F| \overline{ \Omega }^c\times (t/2,\, \infty )\|_{3/2,\,2;\infty }
\le (1+t)^{-\kappa _1}
$.
Due to (\ref{L8.10.30}), we get
$
\|F|B_{\overline{ R}}^c\times (t/2,\, \infty )\|_{\overline{ q}_1, \infty ;\infty }
\le
\mbox{$\mathfrak C$} \,
(1+t)^{-2\kappa _1/ \overline{ q}_1}
$
for
$t \in (0, \infty ) $.
Thus Corollary \ref{corollaryC6.30} with
$\varrho =2,\; q=1,\; s=2,\; \widetilde{ q}=3/2,\; \widetilde{ s}=2,\; \overline{  q}=\overline{  q}_1,\;
\overline{  s}=\infty $
and
$\overline{ R}$ in the place of $R_0$
yields
\begin{eqnarray*}&&
| \partial _x^{\alpha } \mbox{$\mathfrak R$} ^{( \tau )} (F)(x,t)|
\le
\mbox{$\mathfrak C$} \,
\bigl[\, \bigl(\, |x|\,\nu(x) \,\bigr) ^{-1-| \alpha |/2}\, (1+t)^{-\min \{2,\kappa _1\}}
\\&&\hspace{2em}
+
\bigl(\, |x|\,\nu(x) \,\bigr) ^{(-1-| \alpha |/2)\,(1-\epsilon )}
\,
(1+t)^{-\min\{ 1+| \alpha |/2, \; 2 \kappa _1/ \overline{ q}_1\}\, \epsilon }
\,\bigr]
\end{eqnarray*}
for
$x,\, t,\, \epsilon $
as in (\ref{T8.50.10}).
The preceding estimate, (\ref{T8.50.70}) and the definition of $u ^{(2)} $
imply
\begin{eqnarray*}&&
| \partial _x^{\alpha } u ^{(2)} (x,t)|
\le
\mbox{$\mathfrak C$} \,
\bigl[\, \bigl(\, |x|\,\nu(x) \,\bigr) ^{-1-| \alpha |/2}\,
(1+t)^{-\min\{ 3/(2q_1)-1,\; 3 \kappa _1\,(1-1/\hat{q}_1),\; \kappa _1\}}
\\&&\nonumber \hspace{2em}
+
\bigl(\, |x|\,\nu(x) \,\bigr) ^{(-1-| \alpha |/2)\,(1-\epsilon )}
\,
(1+t)^{-\min\{ 1+| \alpha |/2, \; 2 \kappa _1/ \overline{ q}_1\}\, \epsilon }
\,\bigr],
\end{eqnarray*}
again for
$x,\, t,\, \epsilon $
as in (\ref{T8.50.10}).
This result combined with (\ref{T8.10.10}) with $u ^{(1)} $ in the place of $u$
and the equation $u=u ^{(1)} +u ^{(2)} $ lead to (\ref{T8.50.15}).

Suppose in addition that $f$ and $a $ are bounded with compact support, and $\zeta _2 \in [1/2,\, 1)$.
Then $\varrho _1$ may be taken as $\zeta _2$, and $\varrho _2$ as $1+| \alpha |/2$; see the proof
of Corollary \ref{corollaryC8.20}. Moreover $3/(2\, q_1)\uparrow 3/2$ for $q_1 \downarrow 1$,
and $3 \kappa _1\,(1-1/\hat{q}_1)\to \kappa _1$ for $\hat{q}_1 \uparrow 3/2$. If $\alpha =0$, we may take
$\overline{ q}_1=2$. In the case $| \alpha |=1,$ the relation $2 \kappa_1 / \overline{ q}_1\uparrow
2\, \kappa _1/3 $
holds for $\overline{ q}_1\downarrow 3$.
The term $\|f|B_{R_1}^c \times (t/2,\, \infty )\|_{\overline{ q}_0, \overline{ s}_0; \infty }$ may be estimated
by $\mbox{$\mathfrak C$} \, (1+t) ^{-2} $.
These remarks and (\ref{T8.50.15}) imply (\ref{T8.50.16}).
\hfill $\Box$

\end{document}